\documentclass{amsart}
\usepackage[T1]{fontenc}
\usepackage{amsmath,a4wide}
\usepackage{amssymb}
\usepackage{amscd}
\usepackage{epsfig}
\usepackage[T1]{fontenc}
\usepackage{epsfig}
\usepackage{amsmath}
\usepackage{graphicx}
\usepackage{subfigure}

\usepackage{mathrsfs}
\usepackage{epstopdf}
\usepackage{color}
\usepackage{booktabs}
\usepackage{ifpdf}
\usepackage{cite}
\usepackage{amsthm}
\usepackage{amsmath}
\definecolor{myblue}{rgb}{0.2,0,0.9}
\definecolor{blue-violet}{rgb}{0.54, 0.17, 0.89}
\usepackage{enumitem}
\usepackage{algorithm,algorithmic}

\usepackage{pgfplots}
\usepgfplotslibrary{groupplots}
\pgfplotsset{compat=1.12}

\usepackage{tikz}
\usetikzlibrary{arrows,intersections}
\usepackage{boldline}
\usepackage{multirow}
\usepackage[margin=1.32in]{geometry}

\usepackage{varioref}
\usepackage{xr-hyper}
\usepackage{hyperref}
\hypersetup{
	colorlinks,linkcolor=blue,citecolor=blue,filecolor=red,urlcolor=blue}

\definecolor{myblue}{rgb}{0.2,0,0.9}
\definecolor{blue_violet}{rgb}{0.54, 0.17, 0.89}
\definecolor{darkgreen}{rgb}{0,0.35,0}

\allowdisplaybreaks

\makeatletter
\DeclareRobustCommand*\cal{\@fontswitch\relax\mathcal}
\makeatother

\newtheorem{thm}{Theorem}[section]
\newtheorem{pro}[thm]{Proposition}

\newtheorem{cor}[thm]{Corollary}
\newtheorem{lem}[thm]{Lemma}
\numberwithin{equation}{section}

\theoremstyle{definition}
\newtheorem{rem}[thm]{Remark}
\newtheorem{dfn}[thm]{Definition}
\newtheorem{as}[thm]{Assumption}

\usepackage{xparse}
\makeatletter
\RenewDocumentCommand{\title}{om}{%
	\IfNoValueTF{#1}
	{\gdef\shorttitle{}}
	{\gdef\shorttitle{#1}}%
	\gdef\@title{#2}%
}
\makeatother

\makeatletter
\def\@tocline#1#2#3#4#5#6#7{\relax
	\ifnum #1>\c@tocdepth 
	\else
	\par \addpenalty\@secpenalty\addvspace{#2}%
	\begingroup \hyphenpenalty\@M
	\@ifempty{#4}{%
		\@tempdima\csname r@tocindent\number#1\endcsname\relax
	}{%
		\@tempdima#4\relax
	}%
	\parindent\z@ \leftskip#3\relax \advance\leftskip\@tempdima\relax
	\rightskip\@pnumwidth plus4em \parfillskip-\@pnumwidth
	#5\leavevmode\hskip-\@tempdima
	\ifcase #1
	\or\or \hskip 2em \or \hskip 2em \else \hskip 3em \fi%
	#6\nobreak\relax
	\hfill\hbox to\@pnumwidth{\@tocpagenum{#7}}\par
	\nobreak
	\endgroup
	\fi}
\makeatother

\title{Robust mean-field control \\ under common noise uncertainty}
\author{Mathieu Lauri{\`e}re}
\address{Shanghai Frontiers Science Center of Artificial Intelligence and Deep Learning, and NYU-ECNU Institute of Mathematical Sciences, NYU Shanghai}
\email{mathieu.lauriere@nyu.edu}

\author{Ariel Neufeld}
\address{Division of Mathematical Sciences, Nanyang Technological University}
\email{ariel.neufeld@ntu.edu.sg}

\author{Kyunghyun Park}
\address{Division of Mathematical Sciences, Nanyang Technological University}
\email{kyunghyun.park@ntu.edu.sg}

\thanks{\textit{Key words:} mean-field control, common noise uncertainty, robust optimization, propagation of chaos, Markov decision process, dynamic programming.
}
\thanks{\textit{Funding:} 
	M.\;Lauri{\`e}re acknowledges the support of the  grant ``AI-driven Initiative to Promote Research Paradigm Reform and Empower Discipline Advancement.'' Computing resources were provided by NYU Shanghai HPC. A.\;Neufeld acknowledges the support of
	the MOE AcRF Tier 2 Grant MOE-T2EP20222-0013. 
	K.\;Park acknowledges the support of the National Research
	Foundation of Korea (grant DOI: RS-2025-02633175).
}

\date{\today.}

\begin{document}

\begin{abstract}
	We propose and analyze a framework for discrete-time robust mean-field control problems under common\;noise uncertainty. In this framework, the mean-field interaction describes the collective behavior of infinitely many cooperative agents' state and action, while~the common noise---a random disturbance affecting all agents' state dynamics---is uncertain. A~social planner optimizes over open-loop controls on an infinite horizon to maximize the representative agent's worst-case expected reward, where worst-case corresponds to the most adverse probability measure among all candidates inducing the unknown true law of the common noise process. We refer to this optimization as a robust mean-field control problem under common noise uncertainty. 
	We first show that this problem arises as the asymptotic limit of a cooperative $N$-agent robust optimization problem, commonly known as propagation of chaos. We then prove the existence of an optimal open-loop control by linking the robust mean field  control problem to a lifted robust Markov decision problem on the space of probability measures and by establishing the dynamic programming principle and Bellman--Isaac fixed point theorem for the lifted robust Markov decision problem. Finally, we complement our theoretical results with numerical experiments motivated by distribution planning and systemic risk in finance, highlighting the advantages of accounting for common noise uncertainty.

\medskip
\end{abstract}

\maketitle


	\section{Introduction}\label{sec:intro}
	Mean-field control problems \cite{bensoussan2013mean,carmona2018probabilistic}, also known as optimal control of McKean--Vlasov dynamics, have emerged as a fundamental framework for optimizing the behavior of large populations of {\it cooperative} agents. By considering a social planner or central controller managing an infinite (or very large) number of homogeneous agents, mean-field control problems capture a wide range of scenarios including in economics and finance (e.g., \cite{fischer2016continuous,fu2024mean,carmona2021deep,carmona2020applications}), and robotics (e.g., \cite{lerman2004review,elamvazhuthi2019mean,khamis2015multi,cui2023scalable}. 
	
	One significant extension of the mean-field control paradigm is the inclusion of {\it common\;noise}---a random disturbance affecting the dynamics of all agents (e.g., \cite{carmona2023model,motte2022mean,djete2022extended,djete2022mckean,pham2017dynamic,motte2023quantitative}). This feature has become prominent because it captures systemic, correlated randomness (such as macroeconomic shocks or environmental disturbances) that affects the entire population simultaneously. In particular, accounting for common noise enhances the realism of mean-field control problems' applications in financial engineering, including portfolio optimization, optimal liquidation, or systemic risk (e.g., \cite{MR3325083,balata2019class,pham2016linear}), as well as in economics, including contract theory or the production of exhaustible resources (e.g., \cite{elie2021mean,graber2016linear,basei2019weak}).
	
	However, mean-field control problems with common noise 
	inevitably face a key challenge: {\it model uncertainty}.  When a social planner implements a mean-field control problem with common noise, it is likely that there is a margin for potential inaccuracies in the model parameters or distributions governing the common noise process. Crucially, because the common noise process affects all agents simultaneously, even small modeling errors in the common noise process can have widespread impact across our prediction of the system's evolution or our computation of the optimal control. This motivates the need for a {\it robust framework}---also known as the worst-case or Knightian approach (e.g., \cite{chen2002ambiguity,dow1992uncertainty,gilboa1989maxmin,garlappi2007portfolio})---in which the social planner seeks an optimal policy that performs robustly under uncertain dynamics of the common noise.
	
	In this article, we aim to propose and analyze a discrete-time robust mean-field control problem under \emph{common noise uncertainty}. The starting point for our problem is based on the two recent works by Carmona\;et\;al.\cite{carmona2023model} and Motte and Pham \cite{motte2022mean}, where infinite time-horizon discounted mean-field control problems with common noise are considered. Both two works establish the correspondence between the conditional Mckean--Vlasov dynamics for the representative agent's state (that typically appear in mean-field control problems with common noise) and the {\it lifted} Markov decision process on the space of probability measures on the state space. This correspondence enables to articulate dynamic programming Bellman fixed point equations, leading to derive optimal open-loop (and closed-loop Markov) policies for mean-field control problems. Furthermore, \cite{motte2022mean} establishes the propagation of chaos result which connects the mean-field control problem to a social planner's optimization problem with a large but finite number of cooperative agents. This  ensures that the optimal open-loop policy for the mean-field control problem can~be a useful approximation of the optimal policy for such large but finite cooperative agents problems.
	
	Building on \cite{carmona2023model,motte2022mean}, we introduce a probabilistic framework for robust mean-field control problems under {common noise uncertainty}. This framework is designed to encompass both the finite cooperative $N$-agent system and the conditional McKean--Vlasov dynamics when the common noise distribution is unknown (see Section \ref{subsec:propagation}). In contrast with the fixed probability measure setting in \cite{carmona2023model,motte2022mean} which induces a single law for the common noise, we construct a set of probability measures, allowing the common noise to have multiple laws within a prescribed uncertainty measures set (see Definition~\ref{dfn:measures}). 
	This extension is inspired by the robust Markov decision framework of \cite{neufeld2023markov,neufeld2024robust,langner2024markov}, which enables to specify a wide range of different uncertainty sets of probability measures and associated transition kernels.  
	
	Using this framework, we establish three main results. First, we prove a propagation of chaos result linking the finite $N$-agent robust control problem to its mean-field (infinite-agent) counterpart under common noise uncertainty. Under mild regularity conditions on the system and reward functions, we show that the $N$-agent robust control problem converges to the robust mean-field control problem as $N\to \infty$ (see Theorem \ref{thm:propagation_chaos}). This implies that the optimal open-loop policy obtained from the robust mean-field control problem serves as an approximately optimal policy for the finitely many $N$-agent robust control problem. The proof is based on the Wasserstein convergence rates for empirical measures \cite{fournier2015rate,boissard2014mean}. In this regard, our propagation of chaos result can be viewed as a robust analog of the results in \cite{motte2022mean,motte2023quantitative}.
	
	Second, we establish a dynamic programming principle for the robust mean-field control problem by {\it lifting} it to the space of probability measures on the state space.  To that end, we show that the conditional McKean--Vlasov state dynamics under common noise uncertainty corresponds to a lifted robust Markov decision process on the space of probability measures (see Proposition \ref{pro:lift_dynamics}). This correspondence allows us to derive the Bellman–Isaacs fixed-point equations for the value function in the lifted space of distributions. The proof relies on Berge’s maximum theorem to construct local (i.e., one time-step) optimal control and worst-case common noise measure (see Proposition \ref{pro:dpp}), and the Banach fixed-point theorem to establish the existence and uniqueness of a fixed point for the Bellman–Isaacs operator (see Proposition~\ref{pro:fixed_point}). We then construct an optimal open-loop policy for the robust mean-field control problem by aggregating the local optimizers (see Theorem \ref{thm:origin_MDP}). A crucial toolkit in this construction is the use of an extrinsic randomization source with an atomless distribution (see Assumption~\ref{as:randomize}), which also appears in \cite{carmona2023model}. 
	This randomization not only facilitates the implementation of randomized policies in a decentralized manner but also ensures that each agent's distribution of controls aligns with the law of optimal policy prescribed by the social planner. While the existence of a randomization source is not explicitly assumed in \cite{motte2022mean}, a randomization hypothesis on the initial information is imposed therein, which in turn induces a structure from which a randomization source naturally exists; see Remark 3.1 therein.
	
	Third, we introduce a closed-loop Markov policy formulation of the robust mean-field control problem. We establish the equivalence between open-loop and closed-loop formulations (Corollary~\ref{cor:origin_MDP_closed}) and obtain an 
		optimal closed-loop Markov policy. This result can be considered as a robust analog of the main results in \cite{carmona2023model}.

	
	Finally, in order to illustrate all our theoretical results, we provide two numerical examples (see Section~\ref{sec:numerics}). In the first example, inspired by Example 1 of~\cite{carmona2023model}, the central planner's goal is to steer the population distribution towards a target distribution. In the second example, inspired by the systemic risk model of~\cite{MR3325083}, the central planner's goal is to stabilize a financial system and avoid that too many institutions default. In both examples, we underscore the importance and benefits of incorporating common noise uncertainty into mean-field control frameworks.

	\vspace{0.5em}
	\noindent {\it Related literature.} Classic mean-field control problems have been described predominantly in continuous time (see, e.g., \cite{pham2017dynamic,lauriere2016dynamic,fischer2016continuous,fu2024mean,djete2022mckeanAOP,bayraktar2018randomized,soner2024viscosity,burzoni2020viscosity,bensoussan2024control,cosso2019zero,wang2017social,sanjari2020optimal}). Several works \cite{djete2022extended,djete2022mckean,lacker2017limit,fornasier2019mean} have rigorously established the connection between mean-field control and large systems of controlled processes in continuous time settings. 
	
	Notably, robust mean-field control problems in continuous-time settings, involving uncertainty in the drift or volatility of the common noise, have been investigated in \cite{huang2021social,wang2020social,wang2017socialjump}. The conceptual structure of the arguments in \cite{huang2021social} bears certain similarities to ours: in the paper, a centralized control problem under volatility uncertainty of the common noise (analogous to our lifted robust Markov decision problem) is tackled, and then decentralized strategies for the population of agents (analogous to our construction of optimal open-loop policies for the robust mean-field control problem) are obtained.
	Nevertheless, there are key differences. In particular, the continuous-time works rely on the theory of forward-backward stochastic differential equations, which are not suitable in the discrete-time setting we consider. Instead, our analysis requires a measure-theoretic construction of optimal controls and a derivation of the dynamic programming principle on the space of probability measures. Most notably, while the aforementioned works do not establish a propagation of chaos result, our article provides the first such result under common noise uncertainty. 
	
	Several works on mean-field game and control problems have introduced robustness via min–max formulations (e.g., \cite{huang2017robust,liang2022robust,carmona2020policy,carmona2021linear,zaman2024robust}). However, these models do not consider common noise but idiosyncratic noise which is uncertain. 
	In contrast, our framework explicitly accounts for common noise uncertainty, which introduces fundamentally different technical and conceptual challenges. While extending the model to include both idiosyncratic and common noise uncertainty is of clear interest, such an extension leads to significant technical obstacles that invalidate key arguments used in establishing the propagation of chaos result and the lifted dynamic programming principle. This is beyond the scope of the present paper, and we leave it for future work. 
	
	Moving away from the above continuous time settings to discrete time settings, some works \cite{pham2016discrete,gu2023dynamic,gu2021mean,gu2025mean,lauriere2016dynamic} have explored dynamic progarmming principles for discrete time mean-field control problems, but without considering common noise. More relevant to our setting, recent works---including those we benchmark against \cite{motte2022mean,carmona2023model} and others such as \cite{bauerle2023mean,motte2023quantitative,bayraktar2024infinite}---have investigated discrete-time mean-field control problems with common noise. 
	Notably, a recent work \cite{langner2024markov} by two of the authors of the present article proposes a framework for discrete time mean-field Markov games under model uncertainty.  
	 In contrast, we focus on a cooperative control setting (as opposed to a game-theoretic equilibrium) and consider the model uncertainty in the law of the common noise process. This leads to a different optimization structure and our lifted dynamic programming formulation on the space of measures is specifically tailored to this social control setting. Furthermore, our propagation of chaos result has no analogue in \cite{langner2024markov}, whose results concern approximate Nash equilibria rather than centralized performance guarantees.
	
	Finally, for completeness, we note that a substantial body of work has focused on robust Markov decision processes under model uncertainty, which also underpin our lifted dynamic programming result on the space of probability measures (see, e.g., \cite{bauerle2022distributionally,bayraktar2023nonparametric,el2005robust,liu2022distributionally,neufeld2024robust,neufeld2023markov,neufeld2024non,wiesemann2013robust,xu2012distributionally,yang2017convex,li2023policy,lu2025distributionally} in the optimal control literature, and~\cite{hansen2008robustness} in the economics literature).
	
	\section{Main results}\label{sec:main}
	\subsection{Notation and preliminaries}\label{subsec:notation} 
	Throughout this article, we work with Polish spaces. If $X$ is such a space with corresponding metric $d_X$, we denote by ${\cal B}_X$ its Borel $\sigma$-algebra and by ${\cal P}(X)$ the set of all Borel probability measures on~$X$.  Let $C_{b}(X;\mathbb{R})$ be the set of all bounded and continuous functions $f:X\to \mathbb{R}$, endowed with the supremum norm $\|f\|_{\infty}:= \sup_{x \in X} |f(x)|$ where $|\cdot|$ denotes the Euclidean norm. For any $L\geq0$, we denote by $\operatorname{Lip}_{b,L}(X;\mathbb{R})\subset C_{b}(X;\mathbb{R})$ the set of all $L$-Lipschitz continuous~functions. 

	We equip ${\cal P}({X})$ with the topology induced by weak convergence, i.e., for any $\mu\in {\cal P}({X})$ and any $(\mu^n)_{n\in\mathbb{N}}\subseteq {\cal P}({X})$, we have
	\begin{align}\label{eq:topology_w}
		\mu^n \rightharpoonup \mu\;\; \mbox{as $n\rightarrow \infty$}\; \Leftrightarrow \;\; \lim_{n\rightarrow \infty} \int_X f(x) \mu^n(dx) = \int_X f(x) \mu(dx)\;\;\mbox{for any $f \in C_{b}(X;\mathbb{R})$}.
	\end{align}

	If $X$ is compact, then the weak topology given in \eqref{eq:topology_w} is equivalent to the topology induced by the 1-Wasserstein distance ${\cal W}_{{\cal P}(X)}(\cdot,\cdot)$ which we recall to be the following:  For any $\mu,\hat\mu \in {\cal P}(X)$, denote by $\operatorname{Cpl}_{X\times X}(\mu,\hat{\mu})\subset {\cal P}(X\times X)$ the subset of all couplings with marginals~$\mu,\hat{\mu}$. Then the 1-Wasserstein distance between $\mu$ and $\hat{\mu}$ is defined~by
	\[
		{\cal W}_{{\cal P}(X)}(\mu,\hat{\mu}):= \inf_{\Gamma \in  \operatorname{Cpl}_{X\times X}(\mu,\hat{\mu})} \int_{X\times X} d_X(x,y) \Gamma (dx,dy).
	\]

	For each $t\in \mathbb{N}$, we use the abbreviation $X^{t}:=X\times \cdots \times X$ for the $t$-times Cartesian product of the set $X$. Given a sequence $(x_0,\dots,x_{t})\in X^{t+1}$ and $0\leq s \leq t$, we use the following abbreviation $x_{s:t}:=(x_s,\dots,x_{t})$. Then  we endow $X^{t+1}$ with the corresponding product topology induced by the following metric
	: for every $x_{0:t}, \tilde x_{0:t}\in X^{t+1}$, 
	\[
		d_{X^{t+1}}(x_{0:t},\tilde x_{0:t}):=
		\sum_{i=0}^td_X(x_i,\tilde x_i).
	\]
	The same convention applies to a finite Cartesian product of (possibly different) Polish spaces.

	For two Polish spaces $X$ and $Y$, the term `kernel' refers to a Borel measurable map $\lambda:X\ni x\mapsto \lambda(dy|x)\in{\cal P}(Y)$. For every $\mu\in {\cal P}(X)$ and kernel $\lambda$, we write $\mu \mathbin{\hat\otimes} \lambda \in {\cal P}(X\times Y)$ for the measure given by: for every $B\in {\cal B}_{X\times Y}$, $\mu \mathbin{\hat\otimes} \lambda (B):= \int_{X\times Y} {\bf 1}_{\{(x,y)\in B\}}\lambda(dy|x)\mu(dx)$. Moreover for every $\nu\in{\cal P}(Y)$, we write $\mu\otimes \nu\in {\cal P}(X\times Y)$ for the product measure.

	Finally, given $\mu\in {\cal P}(X)$ we use the notation $\mathscr{L}_{\mu}({\cal Z})$ for the law of a random variable~${\cal Z}$ under~${\mu}$ and use $\mathscr{L}_{\mu}({\cal Z}|{\cal Y})$ for the conditional law of ${\cal Z}$ given a random variable ${\cal Y}$ under~${\mu}$. The same convention applies to a $\sigma$-field. We denote by $\delta_x\in {\cal P}(X)$ the Dirac measure at the point $x\in X$.

	\subsection{Propagation of chaos under common noise uncertainty}\label{subsec:propagation}
	In this section, we specify what we mean by the discrete-time $N$-agent model and mean-field control (MFC)  model under common noise uncertainty. We then establish the convergence of the $N$-agent model to the MFC model as the number of agents $N$ goes to infinity.
	
	To that end, we begin by defining a canonical space for the mean-field models with infinitely many indistinguishable agents. 

	Denote by $G$ the initial information space and by $\Theta$ the randomization source space. Moreover denote by $E$ and $E^0$ idiosyncratic and common noise spaces, respectively. 
	On the space defined~by
	\begin{align*}
		\Omega:= \left\{\omega:=\big((g^i)_{i\in \mathbb{N}},(\theta_t^i)_{t\geq 0, i\in \mathbb{N}},(e_t^i)_{t\geq 1, i\in \mathbb{N}},(e^0_t)_{t\geq 1}\big) \; :
		\begin{aligned}
		&\; (g^i,\theta^i_t)\in G \times \Theta,\;\;\mbox{for $t\geq 0$},\; i\in \mathbb{N};
		\\
		&\; (e^i_t, e^0_t)\in E\times E^0,\;\;\mbox{for}\;t\geq 1,\; i\in \mathbb{N}
		\end{aligned}
		\right\},
	\end{align*}
	we denote, for every $\omega \in \Omega$,
	\begin{align}\label{eq:coordinate}
		\begin{aligned}
			\big(\gamma^i( \omega),\vartheta_0^i(\omega)\big)&:=(g^i,\theta^i_0)\in G\times \Theta \quad &&i \in \mathbb{N},\\
			\big(\vartheta_t^i(\omega),\varepsilon^i_t( \omega)\big)&:=(\theta^i_t,e_t^i)\in \Theta\times E \quad &&t\geq 1,\;\; i \in \mathbb{N},\\
			\varepsilon_t^0(\omega)&:=e_t^0\in E^0\quad &&t\geq 1,
		\end{aligned}
	\end{align}
	so that $\gamma^i$  and $(\vartheta_t^i)_{t\geq 0}$ represent the initial state information of agent $i\in \mathbb{N}$ and her 
	randomization source process, respectively. Moreover, $(\varepsilon_t^i)_{t\geq1}$ represents her  idiosyncratic noise process and $(\varepsilon_t^0)_{t\geq1}$ represents the common noise process for all agents.
	
	In what follows, we describe a set of probability measures on the space $\Omega$, which captures model {uncertainty} in the common noise process. 
	\begin{dfn}[Filtrations]\label{dfn:filtrations}
		Consider the following filtrations: for each $i\in \mathbb{N}$
		\begin{itemize}[leftmargin=3.em]
			\item [$\cdot$] $\mathbb{F}^0:=({\cal F}^0_t)_{t\geq 0}$ is 
			given by ${\cal F}^0_t:=\sigma(\varepsilon_{1:t}^0)$ for all $t\geq 1$ with ${\cal F}_0^0=\{\emptyset,\Omega\}$.
			\item [$\cdot$] $\mathbb{F}^{i}:=({\cal F}^{i}_t)_{t\geq 0}$ is given by ${\cal F}_0^{i}:= \sigma(\gamma^i)$ and ${\cal F}^{i}_t:=\sigma(\gamma^i,\vartheta_{0:t-1}^i,\varepsilon_{1:t}^i,\varepsilon_{1:t}^0)$ for all $t\geq 1$.
			\item [$\cdot$] $\mathbb{G}^{i}:=({\cal G}^{i}_t)_{t\geq 0}$ is given by 
			${\cal G}_t^{i}:={\cal F}_t^{i}\vee \sigma(\vartheta_t^i)$ for all $t\geq 0$ so that $\mathbb{F}^{i}\subseteq\mathbb{G}^{i}$. 
		\end{itemize}
		Here ${\cal F}_t^0$ represents the common noise information shared by all agents at time $t$. Both ${\cal F}_t^{i}$ and ${\cal G}_t^{i}$ represent the information of agent $i$ at time $t$, where ${\cal G}_t^{i}$ includes the {\it current} randomization source~$\vartheta_t^i$, while ${\cal F}_t^{i}$ does not. 
	\end{dfn}
	\begin{dfn}[Measures]\label{dfn:measures}
		Fix $\lambda_\gamma\in {\cal P}(G)$, $\lambda_\vartheta \in {\cal P}(\Theta),$ and $\lambda_{\varepsilon}\in {\cal P}(E)$. 
		\begin{itemize}[leftmargin=3.em]
			\item [(i)] Let $\mathfrak{P}^0\subseteq {\cal P}({E}^0)$
			be a non-empty subset of Borel probability measures on $E^0$. Then denote by $\mathcal{K}^0$ the set of all $(p_t)_{t\geq 1}$ consisting of a measure and sequence of kernels such that
			\begin{align*}
				\hspace{3.em}  p_1 \in \mathfrak{P}^0;\qquad p_t:(E^0)^{t-1}\ni e_{1:t-1}^0 \mapsto  p_t(de_t^0|e_{1:t-1}^0) \in \mathfrak{P}^0\quad \mbox{for all $t\geq 2$},
			\end{align*}
			inducing {\it model uncertainty} in the law of the common~noise process $(\varepsilon_t^0)_{t\geq 1}$.
			\item [(ii)] 
			Denote by ${\cal Q}\subseteq {\cal P}(\Omega)$ the subset of all Borel probability measures~$\mathbb{P}$ on $\Omega$ {induced by} 
			some $(p_t)_{t\geq 1}\in \mathcal{K}^0$ in the sense that for every $B_0\in \bigvee_{i \in \mathbb{N}} \mathcal{G}_0^{i}$ and $B_1 \in  \bigvee_{i \in \mathbb{N}}{\cal G}_1^{i}$
			\begin{align*}
				\hspace{4.em}
				\begin{aligned}
				&\mathbb{P}\big\{(\gamma^i,\vartheta_0^i)_{i\in \mathbb{N}}\in B_0\big\}= \hat{Q}_0
				(B_0),\quad \mathbb{P}\big\{((\gamma^i,\vartheta_{0:1}^i, \varepsilon_1^i)_{i\in \mathbb{N}},\varepsilon_1^0)\in B_1\big\}= (\hat{Q}_0 \otimes  \hat Q^{p_1})(B_1),
				\end{aligned}
			\end{align*}
			where  
			\begin{align*}
				\hspace{4.em}
				\begin{aligned}
				\hat Q_0\big((dg^i,d\theta_0^i)_{i\in \mathbb{N}}\big)&:=\mathop{\otimes}\limits_{i \in \mathbb{N}}\big\{(\lambda_\gamma\otimes\lambda_\vartheta)(dg^i,d\theta^i_0)\big\}\in {\cal P}\big((G\times\Theta)^{\mathbb{N}} \big)\\
				\hat Q^{p_1}\big((d\theta^i_1, de^i_1)_{i\in \mathbb{N}}, de^0_1\big)&:=
				\mathop{\otimes}\limits_{i \in \mathbb{N}}\big\{(\lambda_\vartheta\otimes\lambda_\varepsilon)(d\theta^i_1,de_1^i)
				\big\}p_1(de^0_1)\in {\cal P}\big((\Theta\times E)^{\mathbb{N}}\times E^0\big),
				\end{aligned}
			\end{align*}
			whereas for every $t\geq 2$ and $B_t\in  \bigvee_{i \in \mathbb{N}}{\cal G}_t^{i}$
			\begin{align*}
				\hspace{3.em} \mathbb{P}\big\{\big((\gamma^i,\vartheta^i_{0:t}, \varepsilon_{1:t}^i)_{i\in \mathbb{N}},\varepsilon_{1:t}^0\big)\in B_t\big\}= (\hat Q_0 \otimes \hat Q^{p_1}\mathbin{\hat \otimes} \hat Q^{p_2}\mathbin{\hat \otimes} \cdots \mathbin{\hat \otimes} \hat Q^{p_t})(B_t), 
			\end{align*}
			where $\hat Q^{p_t}:(E^0)^{t-1}\ni e_{1:t-1}^0\mapsto \hat Q^{p_t}((d\theta_t^{i},de_t^{i})_{i\in\mathbb{N}},de_t^0| e_{1:t-1}^0)\in {\cal P}((\Theta\times E)^{\mathbb{N}}\times E^0)$ is defined by
			\begin{align*}
				 \hspace{3.em}\hat Q^{p_t}\big((d\theta_t^{i},de_t^{i})_{i\in\mathbb{N}},de_t^0| e_{1:t-1}^0\big):=\mathop{\otimes}\limits_{i \in \mathbb{N}}\big\{(\lambda_\vartheta\otimes\lambda_\varepsilon)(d\theta^i_t,de_t^i)
				 \big\}p_t(de_t^0|e_{1:t-1}^0).
			\end{align*}
		\end{itemize}
	\end{dfn}

	\begin{rem}\label{rem:well_dfn_1}
		By Ionescu--Tulcea's theorem (see, e.g., \cite[Theorem~6.17]{kallenberg2002foundations}), the set ${\cal Q}$ given in Definition~\ref{dfn:measures} is well-defined and the following hold: for every $\mathbb{P}\in {\cal Q}$ w.r.t.\;some $(p_t)_{t\geq 1}\in \mathcal{K}^0$ 
		\begin{itemize}[leftmargin=3.em]
			\item [(i)] 
			$(\gamma^i)_{i\in \mathbb{N}}$, $(\vartheta_t^i)_{t\geq 0,i\in \mathbb{N}}$, $(\varepsilon_t^i)_{t\geq 1,i\in \mathbb{N}}$, and $(\varepsilon_t^0)_{t\geq 1}$ are mutually independent.
			\item [(ii)] $(\gamma^i)_{i\in \mathbb{N}}$ is independent and identically distributed (i.i.d.) with law $\lambda_\gamma$. Moreover, $(\vartheta_t^i)_{t\geq 0,i\in \mathbb{N}}$ is i.i.d. with law $\lambda_{\vartheta}$, and $(\varepsilon_t^i)_{t\geq 1,i\in \mathbb{N}}$ is i.i.d.\;with law $\lambda_\varepsilon$.
			\item [(iii)] $\varepsilon_1^0$ is independent of $\bigvee_{i\in \mathbb{N}}{\cal G}_0^{i}$ with law $p_1$, whereas for every~$t\geq 2$ $\varepsilon_{t}^0$ is conditionally independent of $\bigvee_{i\in \mathbb{N}}{\cal G}_{t-1}^{i}$ given ${\cal F}_{t-1}^0$ (see \cite[Lemma~6.9]{kallenberg2002foundations}), 
			satisfying
			\[
			\mathscr{L}_{\mathbb{P}}(\varepsilon_{t}^0 | {\cal F}_{t-1}^0)= p_t(\,\cdot\,|\varepsilon_{1:t-1}^0) \quad \mbox{$\mathbb{P}$-a.s.}
			\]
		\end{itemize}
		We note that when $\mathfrak{P}^0$ is a singleton (i.e., {\it without uncertainty}), the resulting probabilistic framework
		coincides with the setting in \cite[Section~2.1.2]{carmona2023model} and is also similar to the one in \cite[Section~2]{motte2022mean}. 
	\end{rem}

	We introduce a dynamical system of mean-field models with indistinguishable $N$-agents under common noise uncertainty and define the corresponding robust optimization problem. To this end, let us introduce the following elementary components:
	\begin{dfn}\label{dfn:basic_element}
		Let $S$ and $A$ be nonempty compact Polish spaces, representing the state and action spaces, respectively.
		\begin{itemize}[leftmargin=3.em]
			\item [(i)] $\operatorname{F}:S \times A \times {\cal P}(S\times A) \times E \times E^0\to S$ 
			is a Borel measurable transition function describing the dynamics of each of the $N$-agents as well as the mean-field model.
			\item [(ii)] $r:S\times A\times {\cal P}(S\times A)\to \mathbb{R}$ is a Borel measurable one-step reward function.
			\item [(iii)] $\beta \in [0,1)$ is a discount factor.  
		\end{itemize}
	\end{dfn}
	\begin{dfn}[$N$-agent model]\label{dfn:N_model}
		Recall that for each $i\in \mathbb{N}$, ${\cal F}_0^i=\sigma(\gamma^i)$ (see Definition \ref{dfn:filtrations}).  Denote for every $i\in \mathbb{N}$ by $L^0_{{\cal F}_0^i}(S)$ the set of all ${\cal F}_0^i$ measurable random variables with values in $S$.
		\begin{enumerate}[leftmargin=3.em]
			\item [(i)] Denote by $\Pi$ the set of all open-loop policies $(\pi_t)_{t\geq 0}$ in the sense that $\pi_t:G\times \Theta^{t+1}\times E^t\times (E^0)^t\to A$ is a Borel measurable function for all $t\geq 0$. Given $(\pi_t)\in \Pi$, the action process of agent $i\in \mathbb{N}$ is given by the open-loop control
			\[
			\hspace{1.em} a_{t}^{i,\pi}:= \pi_t(\gamma^i,\vartheta_{0:t}^i,\varepsilon_{1:t}^{i},\varepsilon_{1:t}^0)\quad t\geq 1,\quad \mbox{with}\;\; a_0^{i,\pi}:=\pi_0(\gamma^i,\vartheta_0^i).
			\]
			In other words, $(a_t^{i,\pi})_{t\geq 0}$ is a $\mathbb{G}^i$ adapted process (see Definition \ref{dfn:filtrations}).
			\item [(ii)] Fix the initial state $\xi^{i}\in L^0_{{\cal F}_0^i}(S)$ of agent~$i$. Given $N\in \mathbb{N}$ and $(\pi_t)_{t\geq0}\in \Pi$, the state and action processes of agent $i=1,\dots,N$ in the $N$-agent model under~$\mathbb{P}\in {\cal Q}$ are given by 
			\begin{align}\label{eq:MKV_N}
				\hspace{1.em} 
				\left\{
				\begin{aligned}
					&s_0^{i,N,\pi}:=\xi^{i},\\
					&s_{t+1}^{i,N,\pi}:=\operatorname{F}(s^{i,N,\pi}_t,a^{i,\pi}_t,\mbox{$\frac{1}{N}\sum_{j=1}^N\delta_{(s^{j,N,\pi}_t,a^{j,\pi}_t)}$},\varepsilon_{t+1}^{i},\varepsilon_{t+1}^0)\quad  t\geq 0.
				\end{aligned}
				\right.
			\end{align}
			Here we observe that both the law of the initial state and action $(s_0^{i,N,\pi},a_0^{i,\pi})$ and the law of the idiosyncratic noise process $(\varepsilon_t^i)_{t\geq 0}$ do not depend the choice of $\mathbb{P}\in {\cal Q}$ (see Definition\;\ref{dfn:measures}\;(iii)). In contrast, the law of $(s_{t}^{i,N,\pi},a_t^{i,\pi})$ for $t\geq 1$ depends on this choice, due to the model uncertainty in $(\varepsilon_{t}^0)_{t\geq 1}$.
			\item [(iii)] The contribution of agent~$i$ to the social planner’s gain over an infinite horizon under $\mathbb{P}\in {\cal Q}$ is defined by
			\[
			\quad R^{i,N,\pi}:= \sum_{t=0}^\infty \beta^t r(s_{t}^{i,N,\pi},a_{t}^{i,\pi},\mbox{$\frac{1}{N}\sum_{j=1}^N\delta_{(s^{j,N,\pi}_t,a^{j,\pi}_t)}$})\quad i=1,\dots,N.
			\]
			Then the social planner’s worst-case expected gain under the common noise uncertainty is 
			\[
			{\cal J}^{N,\pi}:=\inf_{\mathbb{P}\in {\cal Q}}\mathbb{E}^{\mathbb{P}}[ R^{N,\pi} ]\quad \mbox{where} \;\;  R^{N,\pi}:=\frac{1}{N}\sum_{i=1}^NR^{i,N,\pi},
			\]
			and the resulting $N$-agent optimization problem is given by 
			$V^{N}:= \sup_{\pi \in \Pi} {\cal J}^{N,\pi}.$ This problem is a robust analog of the classical $N$-agent optimization problem of \cite{carmona2023model,motte2022mean}. 
		\end{enumerate}
	\end{dfn}

	In light of the propagation of chaos argument, we expect and aim to show that the asymptotic version of the $N$-agent problem in Definition~\ref{dfn:N_model}, as $N \to \infty$, is given by the following:
	\begin{dfn}[MFC model]\label{dfn:MFC}
		For each $i \in \mathbb{N}$, let $\xi^i\in L_{{\cal F}_0^i}^0(S)$ be the fixed initial state of agent~$i$; see Definition~\ref{dfn:N_model}\;(ii).
		\begin{enumerate}[leftmargin=3.em]
			\item [(i)] Given $(\pi_t)_{t\geq0}\in \Pi$, the state process of agent $i\in \mathbb{N}$  in the infinite population model under $\mathbb{P}\in {\cal Q}$ is governed by the conditional McKean--Vlasov dynamics:
			\begin{align}\label{eq:MKV}
				\hspace{1.em} 
				\left\{
				\begin{aligned}
					s_0^{i,\pi,\mathbb{P}}&:=\xi^{i},\\
					s_{t+1}^{i,\pi,\mathbb{P}}&:=\operatorname{F}(s^{i,\pi,\mathbb{P}}_t,a^{i,\pi}_t,\mathbb{P}^{0}_{(s^{i,\pi,\mathbb{P}}_t,a^{i,\pi}_t)},\varepsilon_{t+1}^{i},\varepsilon_{t+1}^0)\quad t\geq 0,
				\end{aligned}\right.
			\end{align}
			where $(a^{i,\pi}_t)_{t\geq 0}$ is the open-loop control of agent $i$ as defined in Definition \ref{dfn:N_model}\;(i), and $\mathbb{P}^{0}_{(s^{i,\pi,\mathbb{P}}_t,a^{i,\pi}_t)}$ is the conditional joint law of $(s^{i,\pi,\mathbb{P}}_t,a^{i,\pi}_t)$ under $\mathbb{P}$ given the common noise trajectory $\varepsilon_{1:t}^0$, i.e.,
			\[
			\hspace{1.em}  \mathbb{P}^0_{(s^{i,\pi,\mathbb{P}}_t,a^{i,\pi}_t)}:=\mathscr{L}_{\mathbb{P}}\big((s^{i,\pi,\mathbb{P}}_t,a^{i,\pi}_t) | \varepsilon_{1:t}^0\big)\quad t\geq 1
			\]
			with the convention that $\mathbb{P}^0_{(s^{i,\pi,\mathbb{P}}_0,a^{i,\pi}_0)}:=\mathscr{L}_{\mathbb{P}}((s^{i,\pi,\mathbb{P}}_0,a^{i,\pi}_0))$. Analogously, for every $t\geq 1$ let $\mathbb{P}^0_{s^{i,\pi,\mathbb{P}}_t}$ be the conditional law of $s^{i,\pi,\mathbb{P}}_t$ under $\mathbb{P}$ given the common noise trajectory $\varepsilon_{1:t}^0$ with the convention that $\mathbb{P}^0_{s^{i,\pi,\mathbb{P}}_0}:=\mathscr{L}_{\mathbb{P}}(s^{i,\pi,\mathbb{P}}_0)$.
			\item[(ii)] The contribution of agent $i$ to the social planner’s gain under $\mathbb{P}\in {\cal Q}$ is defined~by 
			\[
			\quad R^{i,\pi,\mathbb{P}}:= \sum_{t=0}^\infty \beta^t r(s_{t}^{i,\pi,\mathbb{P}},a_{t}^{i,\pi},\mathbb{P}^0_{(s^{i,\pi,\mathbb{P}}_t,a^{i,\pi}_t)})\quad i\in \mathbb{N}.
			\]
			Then the social planner’s worst-case expected gain under the common noise uncertainty is 
			\begin{align}\label{eq:worst_gain}
			\qquad {\cal J}^{\pi}:=\inf_{\mathbb{P}\in {\cal Q}}\mathbb{E}^{\mathbb{P}}[R^{\pi,\mathbb{P}} ],\qquad \mbox{where} \;\;  R^{\pi,\mathbb{P}}:=\mathbb{E}^{\mathbb{P}^0}[R^{i,\pi,\mathbb{P}}]=\mathbb{E}^{\mathbb{P}^0}[R^{1,\pi,\mathbb{P}}]\quad i\in \mathbb{N},
			\end{align}
			where $\mathbb{E}^{\mathbb{P}^0}[\cdot]$ denotes the conditional expectation under $\mathbb{P}$ given $(\varepsilon^0_t)_{t\geq 0}$ and the quantity $ R^{\pi,\mathbb{P}}$ is independent of the choice of~$i$  due to the indistinguishability of agents. The resulting robust MFC problem is then defined as 
			$V:= \sup_{\pi \in \Pi} {\cal J}^{\pi}.$ 
		\end{enumerate}
	\end{dfn}

	The main goal of this section is to rigorously connect the $N$-agent model in Definition \ref{dfn:N_model} with the MFC model in Definition \ref{dfn:MFC}.

	We impose the following conditions on the basic components given in Definition~\ref{dfn:basic_element}. 
	\begin{as}\label{as:general}
		The following conditions hold:
		\begin{itemize}[leftmargin=3.em]
			\item [(i)] There exists some $C_{{\operatorname{F}}}>0$ such that for every $s,\tilde s \in S$, $a\in A$, $\Lambda, \tilde \Lambda \in {\cal P}(S\times A)$, and $e^0\in E^0$
			\begin{align*}
				\hspace{3.em} 
				\int_{E}d_S\big(\operatorname{F}(s,a,\Lambda,e,e^0),\operatorname{F}(\tilde s,a,\tilde\Lambda,e,e^0)\big)\lambda_{\varepsilon}(de)
				\leq C_{{\operatorname{F}}} \big(d_S(s,\tilde s)+\mathcal{W}_{{\cal P}(S\times A)}(\Lambda, \tilde \Lambda )\big),
			\end{align*}
			where $\lambda_\varepsilon$ is given in Definition \ref{dfn:N_model}\;(i).
			\item [(ii)] There exists $C_{{r}}>0$ such that for every $s,\tilde s \in S$, $a\in A$, and $\Lambda,\tilde \Lambda \in {\cal P}(S\times A)$
			\[
			\hspace{3.em} |r (s,a,\Lambda)|\leq C_{{r}},\qquad  
			|r (s,a,\Lambda)-r (\tilde s, a,\tilde \Lambda)|\leq C_{{r}}\big(d_S(s,\tilde s)+\mathcal{W}_{{\cal P}(S\times A)}(\Lambda, \tilde \Lambda )\big).
			\] 
			\item [(iii)] $\beta $ is in $[0,1\wedge (2C_{{\operatorname{F}}})^{-1})$. 
		\end{itemize} 
	\end{as}
	
	For every $N\in \mathbb{N}$, we define the following quantity 
	\begin{align}\label{eq:propag_MN}
		M_N:=\sup_{t\geq 0} \sup_{\pi\in \Pi} \sup_{\mathbb{P}\in {\cal Q}}\mathbb{E}^{\mathbb{P}}\bigg[{\cal W}_{{\cal P}(S\times A)}\Big(\frac{1}{N}\sum_{i=1}^N\delta_{(s^{i,\pi,\mathbb{P}}_t,a^{i,\pi}_t)},\, \mathbb{P}^0_{(s^{1,\pi,\mathbb{P}}_t,a^{1,\pi}_t)}\Big)\bigg],
	\end{align}
	where for each $j=1,\cdots,N$, $(s^{j,\pi,\mathbb{P}}_t,a^{j,\pi}_t)_{t\geq 0}$ are the state and action processes of agent $j$ under~$\mathbb{P}$ in the MFC model, and for each $t\geq 0$ $\mathbb{P}^0_{(s^{1,\pi,\mathbb{P}}_t,a^{1,\pi}_t)}$ is the conditional joint law of $(s^{1,\pi,\mathbb{P}}_t,a^{1,\pi}_t)$ under $\mathbb{P}$ given the common noise $\varepsilon_{1:t}^0$ (see Definition \ref{dfn:MFC}). By the indistinguishability of the $N$ agents,  $\mathbb{P}^0_{(s^{1,\pi,\mathbb{P}}_t,a^{1,\pi}_t)}$ can equivalently be replaced by $\mathbb{P}^0_{(s^{j,\pi,\mathbb{P}}_t,a^{j,\pi}_t)}$ for any $j\in \mathbb{N}$.
	
	The following estimates on the sequence $(M_N)_{N\in \mathbb{N}}$, as defined in \eqref{eq:propag_MN}, follow from standard applications of the non asymptotic bounds for the convergence rate of empirical measures in Wasserstein distance (see \cite[Theorem 1]{fournier2015rate}, \cite[Corollary~1.2]{boissard2014mean}).
	\begin{lem}\label{lem:propagation_chaos} Denote by  $\Delta_{S\times A}\in [0,\infty)$ the diameter of $S\times A$.
		Then the following hold: 
		\begin{itemize}
			\item [(i)] If $S\times A \subset \mathbb{R}^d$ for some $d\in \mathbb{N}$, then for any $q>2$ there exists some constant $C>0$ (that depends only on $d$ and $q$) such that for every $N\in \mathbb{N}$,
			\begin{align*}
			M_N \leq C \Delta_{S\times A}\cdot \alpha(N)<\infty,
			\end{align*}
			where $\alpha:\mathbb{N}\ni N\mapsto \alpha(N)\in(0,\infty)$ is given as follows: $\alpha(N):=N^{-1/2}$  for $d=1$; $\alpha(N):=N^{-1/2}\log(1+N)$ for $d=2$; $\alpha(N):=N^{-1/d}\log(1+N)$ for $d\geq 3$.
			\item [(ii)] If for every $\delta >0$ there exist some constants $k_{S\times A}>0$ and $q>2$ such that the minimal number of balls with radius $\delta$ covering $S\times A$, denoted by $\underline n(S\times A,\delta)\in \mathbb{N}$, satisfies $\underline n(S\times A,\delta) \leq k_{S\times A} \big({\Delta_{S\times A}}\cdot{\delta}^{-1}\big)^q$,  then there exists some $C>0$ (that depends only on $k_{S\times A}$ and $q$) such that for every $N\in \mathbb{N}$,
			\[
			M_N \leq C \Delta_{S\times A}\cdot N^{-\frac{1}{q}}<\infty.
			\]
		\end{itemize}
	\end{lem}
	
	
	By using Lemma \ref{lem:propagation_chaos}, we can obtain a rate of convergence when approximating the  $N$-agent model by the MFC model under model uncertainty in the common noise process. 
	\begin{thm}\label{thm:propagation_chaos}
		Suppose that Assumption \ref{as:general} holds. Moreover, we assume that $\Delta_{S\times A}$ satisfies one of the two settings imposed in Lemma~\ref{lem:propagation_chaos}. Then it holds that for every $N\in \mathbb{N}$, $i=1,\dots,N$, and $t\geq 0$  
		\begin{align}
			&\sup_{\pi \in \Pi}\sup_{\mathbb{P}\in {\cal Q}}\mathbb{E}^{\mathbb{P}}\big[d_S(s^{i,N,\pi}_t,s^{i,\pi,\mathbb{P}}_t)\big] = O(M_N),\label{eq:cvg_state}\\
			&\sup_{\pi \in \Pi} \sup_{\mathbb{P}\in {\cal Q}}\mathbb{E}^{\mathbb{P}}\bigg[{\cal W}_{{\cal P}(S\times A)}\bigg(\frac{1}{N}\sum_{j=1}^N\delta_{(s^{j,N,\pi}_t,a^{j,\pi}_t)},\, \mathbb{P}^0_{(s^{i,\pi,\mathbb{P}}_t,a^{i,\pi}_t)}\bigg)\bigg]= O(M_N), \label{eq:cvg_emprical}
		\end{align}
		where $O(\cdot)$ is the Landau symbol. Moreover, 
		there exists some constant $C>0$ (that depends only on $C_{\operatorname{F}},C_r$ and $\beta$)  such that for $N\in \mathbb{N}$ sufficiently large
		\begin{align}\label{eq:cvg_reward}
			\sup_{\pi \in \Pi} |{\cal J}^{N,\pi}-{\cal J}^{\pi}|\leq C M_N,
		\end{align} 
		which ensures that $|V^N-V|= O(M_N)$. Consequently, any $\varepsilon$-optimal policy for the robust MFC problem $V$ (see Definition \ref{dfn:MFC}) is $O(\varepsilon)$-optimal for the $N$-agent robust optimization problem $V^N$ (see Definition \ref{dfn:N_model}) if $N$ is sufficiently large such that $M_N=O(\varepsilon)$. Conversely, any $\varepsilon$-optimal policy for $V^N$ is $O(\varepsilon)$-optimal for $V$ if $N\in \mathbb{N}$ is sufficiently large such that $M_N=O(\varepsilon)$.
	\end{thm}
	The proofs of Lemma \ref{lem:propagation_chaos} and Theorem \ref{thm:propagation_chaos} can be found in Section \ref{proof:subsec:propagation}.
	
	\begin{rem}
		Theorem \ref{thm:propagation_chaos} can be viewed as a robust analog of \cite[Theorem 2.1]{motte2022mean}. The overall proof roadmap follows the arguments in the reference, where the convergence rate of the empirical measure (see Lemma \ref{lem:propagation_chaos}) plays a key role. Moreover, the Lipschitz conditions on the one-step reward and system functions in Assumption~\ref{as:general}\;(i),\;(ii)  (denoted as ${\bf \operatorname{\bf Hf}_{\operatorname{\bf lip}}}$ and ${\bf \operatorname{\bf HF}_{\operatorname{\bf lip}}}$ therein), together with a certain condition on the discount factor (similar to Assumption~\ref{as:general}\;(iii)), are imposed. While our setting is more rigid due to the uncertainty measures set ${\cal Q}$, we are able to obtain the propagation of chaos result by establishing the convergence rate of the empirical measure  uniformly over all probability measures $\mathbb{P}\in {\cal Q}$.
	\end{rem}

	\subsection{Lifted robust Markov decision processes on the space of probability measures}\label{subsec:lift_MDP}
	Theorem \ref{thm:propagation_chaos} shows that the robust MFC model in Definition \ref{dfn:MFC} serves as a macroscopic approximation of the robust $N$-agent optimization model in Definition \ref{dfn:N_model}.  
	By definition of the conditional McKean-Vlasov dynamics \eqref{eq:MKV} and the social planner's worst-case expected gain~\eqref{eq:worst_gain}, 
	we can without loss of generality consider only one representative agent. 
	
	Accordingly, we suppress the index $i\in \mathbb{N}$ representing individual agents, and denote the representative agent’s components as follows: the initial information is given by $\gamma$, the randomization source process by $(\vartheta_t)_{t\geq 0}$, the idiosyncratic noise by $(\varepsilon_t)_{t\geq 1}$, and the information processes by 
	\begin{align}\label{eq:filtrations_represent}
		\begin{aligned}
		\mathbb{F}&:=({\cal F}_t)_{t\geq 0}\quad \mbox{with ${\cal F}_0:= \sigma(\gamma)$ and ${\cal F}_t:=\sigma(\gamma,\vartheta_{0:t-1},\varepsilon_{1:t},\varepsilon_{1:t}^0)$ for all $t\geq 1$;}\\
		\mathbb{G}&:=({\cal G}_t)_{t\geq 0}\quad \mbox{with ${\cal G}_t:={\cal F}_t\vee \sigma(\vartheta_t)$ for all $t\geq 0$ so that $\mathbb{F}\subseteq\mathbb{G}$},
		\end{aligned}
	\end{align}
	see Definition~\ref{dfn:filtrations}. The initial state is then given by $\xi\in L_{{\cal F}_0}^0(S)$.
	
	Moreover, we define by
	\begin{align}\label{eq:open_loop_set}
		{\cal A}:=\left\{a:=(a_t)_{t\geq 0}:
		\begin{aligned}
			&\,\mbox{$a$ is $\mathbb{G}$ adapted and satisfies $a_t=\pi_t(\gamma,\vartheta_{0:t},\varepsilon_{1:t},\varepsilon^0_{1:t})$ for $t\geq 1$}\\
			&\,\mbox{and $a_0=\pi_0(\gamma,\vartheta_0)$ w.r.t.\;some $\pi\in \Pi$}
		\end{aligned}
		\right\},
	\end{align}
	the set of open-loop controls of the representative agent (see Definition \ref{dfn:N_model}\;(i) for the notation $\Pi$).
	
	Given $a \in {\cal A}$, the state process of the representative agent in the infinite population model under $\mathbb{P}\in {\cal Q}$ evolves according to the conditional McKean-Vlasov dynamics:
	\begin{align}\label{eq:MKV_represent}
		\hspace{1.em} 
		s_{t+1}^{\xi,a,\mathbb{P}}:=\operatorname{F}(s^{\xi,a,\mathbb{P}}_t,a_t,\mathbb{P}^{0}_{(s^{\xi,a,\mathbb{P}}_t,a_t)},\varepsilon_{t+1},\varepsilon_{t+1}^0)\quad \mbox{for $t\geq 0$},\;\; \mbox{with $\;\;s_0^{\xi,a,\mathbb{P}}:=\xi,$}
	\end{align}
	where $\mathbb{P}^0_{(s^{\xi,a,\mathbb{P}}_t,a_t)}$ is the conditional joint law of $(s^{\xi,a,\mathbb{P}}_t,a_t)$ under $\mathbb{P}$  given $\varepsilon^0_{1:t}$ for $t\geq 1$, with the convention that $\mathbb{P}^0_{(s^{\xi,a,\mathbb{P}}_0,a_0)}:=\mathscr{L}_{\mathbb{P}}((s^{\xi,a,\mathbb{P}}_0,a_0))$. Here we note that $(s_t^{\xi,a,\mathbb{P}})_{t\geq 0}$ is $\mathbb{F}$ adapted and $(\mathbb{P}^{0}_{(s^{\xi,a,\mathbb{P}}_t,a_t)})_{t\geq 0}$ is $\mathbb{F}^0$ adapted (see Lemma \ref{lem:measurability}\;(ii)). 
	
	Then the social planner’s worst-case expected gain under the common noise uncertainty is 
	\begin{align}\label{eq:worst_represent}
	\;\;{\cal J}^{a}(\xi):=\inf_{\mathbb{P}\in {\cal Q}}\mathbb{E}^{\mathbb{P}}[R^{a,\mathbb{P}}(\xi) ],\quad \mbox{where} \;\;  R^{a,\mathbb{P}}(\xi):=\mathbb{E}^{\mathbb{P}^0}\bigg[\sum_{t=0}^\infty \beta^t r(s_{t}^{\xi,a,\mathbb{P}},a_{t},\mathbb{P}^0_{(s^{\xi,a,\mathbb{P}}_t,a_t)})\bigg].
	\end{align}
	Accordingly, the robust MFC problem of the social planner is defined by
	\begin{align}\label{eq:MFC_represent}
		V(\xi):=\sup_{a\in {\cal A}} {\cal J}^{a}(\xi),\quad\;\xi \in L_{{\cal F}_0}^0(S).
	\end{align}
	This formulation coincides with Definition \ref{dfn:MFC} (by suppressing the agent index $i$).
	
	We now show how the robust MFC problem given in 
	\eqref{eq:MFC_represent} can be {\it lifted} to a robust Markov decision process (MDP) under model uncertainty in the space of probability measures. Given $\xi\in  L_{{\cal F}_0}^0(S)$, $a\in {\cal A}$, and $\mathbb{P}\in {\cal Q}$, we define the following $\mathbb{F}^0$ adapted processes:
	\begin{align}
		&(\mu_t^{\xi,a,\mathbb{P}})_{t\geq 0}:=(\mathbb{P}^0_{s_t^{\xi,a,\mathbb{P}}})_{t\geq 0}\subseteq {\cal P}(S)\label{eq:MDP_state},\\
		&(\Lambda_t^{\xi,a,\mathbb{P}})_{t\geq 0}:=(\mathbb{P}^0_{(s_t^{\xi,a,\mathbb{P}},a_t)})_{t\geq 0}\subseteq {\cal P}(S\times A).\label{eq:MDP_action}
	\end{align}
	We refer to \eqref{eq:MDP_state} and \eqref{eq:MDP_action} as the {lifted} state and lifted action processes, respectively. Note that the lifted processes satisfy the following marginal constraint: $\mathbb{P}$-a.s.,
	\begin{align}\label{eq:marginal_const}
		\operatorname{pj}_S(\Lambda_t^{\xi,a,\mathbb{P}})= \mu_t^{\xi,a,\mathbb{P}}\quad \mbox{for all $t\geq 0$,}
	\end{align}
	where $\operatorname{pj}_S:{\cal P}(S\times A)\ni\Lambda\mapsto \operatorname{pj}_S(\Lambda):=\Lambda(\cdot \times A)\in {\cal P}(S)$ denotes the projection function that maps $\Lambda$ onto its marginal on $S$. 
	
	Based on this observation, we first characterize the dynamics of the {lifted} state processes.
	To that end, let us introduce some notation and functions defined on the spaces of probability measure, ${\cal P}(S)$ and ${\cal P}(S\times A)$ (we refer to them as the `lifted' spaces), which is convenient to characterize the dynamics and then to obtain the lifted dynamic programming principle.
	\begin{dfn}\label{dfn:lift_maps}
		Let $\lambda_\varepsilon\in {\cal P}(E)$ be given in Definition \ref{dfn:measures}. Moreover, let $\operatorname{F}$
		and $r$ be the transition function and one-step reward function, respectively, as defined in Definition \ref{dfn:basic_element}\;(i).
		\begin{enumerate}[leftmargin=3.em]
			\item [(i)] Denote by
			\begin{align*}
				\hspace{2.em}\mathfrak{U}: {\cal P}(S)\ni \mu \twoheadrightarrow\mathfrak{U}(\mu):= \{\Lambda\in {\cal P}(S\times A):\operatorname{pj}_S(\Lambda)=\mu\} \subseteq {\cal P}(S\times A)
			\end{align*}
			the {correspondence} (i.e., a set-valued map) inducing the marginal constraint on $S$. Moreover, denote by $\operatorname{gr}(\mathfrak{U})$ the graph of $\mathfrak{U}$, i.e., $\operatorname{gr}(\mathfrak{U}):= \{(\mu,\Lambda)\in {\cal P}(S)\times {\cal P}(S\times A):\Lambda \in \mathfrak{U}(\mu) \}.$  
			\item [(ii)] Denote by $\overline{\mathrm{F}}:\operatorname{gr}(\mathfrak{U})\times E^0\ni (\mu,\Lambda,e^0)\mapsto \overline{\mathrm{F}}(\mu,\Lambda,e^0)\in{\cal P}(S)$ the {lifted} transition function given by 
			\begin{align*}
				\hspace{1.em} \overline{\mathrm{F}}(\mu,\Lambda,e^0)(ds'):=\big((\Lambda \otimes \lambda_{\varepsilon})\circ \operatorname{F}(\cdot,\cdot,\Lambda,\cdot,e^0)^{-1}\big)(ds'),
			\end{align*}
			i.e., the push-forward of $\Lambda \otimes \lambda_{\varepsilon}\in {\cal P}(S\times A\times E)$ by $\operatorname{F}(\cdot,\cdot,\Lambda,\cdot,e^0):S\times A\times E \to S$.
			\item [(iii)] Let $\overline p:\operatorname{gr}(\mathfrak{U})\times {\cal P}(E^0)\ni (\mu,\Lambda,p)\mapsto \overline{p}(d\mu'|\mu,\Lambda,p)\in{\cal P}({\cal P}(S))$ be a kernel defined~by 
			\begin{align*}
				\overline p(d\mu'|\mu,\Lambda,p):= \big(p\circ \overline{\mathrm{F}}(\mu,\Lambda,\cdot )^{-1}\big)(d\mu'),
			\end{align*}
			i.e., the push-forward of ${p} \in {\cal P}(E^0)$ by $\overline{\mathrm{F}}(\mu,\Lambda,\cdot ):E^0 \to {\cal P}(S)$.
			\item [(iv)] Denote by $\overline{r}: \operatorname{gr}(\mathfrak{U})\ni (\mu,\Lambda) \mapsto \overline{r}(\mu,\Lambda)\in \mathbb{R}$ the {lifted} reward function defined~by
			\[
			\overline{r}(\mu,\Lambda):= \int_{S\times A} r(s,a,\Lambda)\Lambda(ds,da).
			\]
		\end{enumerate}
	\end{dfn}

	The following lemma shows that indeed $(\mu^{\xi,a,\mathbb{P}}_t)_{t\geq 0}$ given in \eqref{eq:MDP_state} can be seen as an MDP on the space of probability measures.

	\begin{pro}\label{pro:lift_dynamics}
		Let $\overline{\operatorname{F}}$ and $\overline{p}$ be given in Definition \ref{dfn:lift_maps}. 
		Let $\xi\in  L_{{\cal F}_0}^0(S)$, $a\in {\cal A},$ and $\mathbb{P}\in {\cal Q}$ be given where $\mathbb{P}$ is {induced by} some couple $(p_t)_{t\geq 1}\in \mathcal{K}^0$ (see Definition \ref{dfn:measures}). Then the lifted~state and action processes $(\mu_t^{\xi,a,\mathbb{P}})_{t\geq 0}$ and $(\Lambda_t^{\xi,a,\mathbb{P}})_{t\geq 0}$ (see \eqref{eq:MDP_state}, \eqref{eq:MDP_action}) satisfy for every $t\geq 0$, $\mathbb{P}$-a.s.
		\begin{align}
			\mu_{t+1}^{\xi,a,\mathbb{P}} = \overline{\mathrm{F}}(\operatorname{pj}_S(\Lambda_t^{\xi,a,\mathbb{P}}),\,\Lambda_t^{\xi,a,\mathbb{P}},\,\varepsilon_{t+1}^0) \label{eq:lift_1},
		\end{align}
		which implies that $\mathbb{P}$-a.s.
		\begin{align}
			\begin{aligned}
			\mathscr{L}_{\mathbb{P}}(\mu_{1}^{\xi,a,\mathbb{P}})&=\overline p(\,\cdot\,|\operatorname{pj}_S(\Lambda_0^{\xi,a,\mathbb{P}}),\Lambda_0^{\xi,a,\mathbb{P}},p_1(\cdot)),\\
			\mathscr{L}_{\mathbb{P}}(\mu_{t+1}^{\xi,a,\mathbb{P}})&=\overline p(\,\cdot\,|\operatorname{pj}_S(\Lambda_t^{\xi,a,\mathbb{P}}),\Lambda_t^{\xi,a,\mathbb{P}},p_{t+1}(\,\cdot\,|\varepsilon^0_{1:t}))\quad \mbox{for all $t\geq 1$}. \label{eq:lift_2}
			\end{aligned}
		\end{align}
	\end{pro}
	The proof of Proposition \ref{pro:lift_dynamics} can be found in Section \ref{proof:subsec:lift_MDP}.
	
	\begin{rem}\label{rem:lift_gains}
		Let $\xi\in  L_{{\cal F}_0}^0(S)$, $a\in {\cal A}$, and $\mathbb{P}\in {\cal Q}$  be given. 
		Note that for every $t\geq 0$,
		\begin{align}\label{eq:lift_reward}
			\begin{aligned}
			\mathbb{E}^{\mathbb{P}}\big[{r}(s_{t}^{\xi,a,{\mathbb{P}}},a_t,{\Lambda}_{t}^{\xi,a,\mathbb{P}})\big]&= \mathbb{E}^{\mathbb{P}}\bigg[\mathbb{E}^{\mathbb{P}}\bigg[\int_{S\times A}{r}(\tilde s,\tilde a,{\Lambda}_{t}^{\xi,a,\mathbb{P}}){\Lambda}_{t}^{\xi,a,\mathbb{P}}(d\tilde s,d\tilde a)\,\bigg|\,{\cal F}_t^0\bigg]\bigg]\\
			&=\mathbb{E}^{\mathbb{P}}\big[\overline{r}(\operatorname{pj}_S({\Lambda}_t^{\xi,a,\mathbb{P}}),{\Lambda}_t^{\xi,a,\mathbb{P}})\big]=\mathbb{E}^{\mathbb{P}}\big[\overline{r}(\mu_t^{\xi,a,\mathbb{P}},{\Lambda}_t^{\xi,a,\mathbb{P}})\big],
			\end{aligned}
		\end{align}
		where the first equality holds by ${\cal F}_t^0$-measurability of ${\Lambda}_{t}^{\xi,a,\mathbb{P}}$ (see Lemma \ref{lem:measurability}\;(ii)), the second equality follows from the definition of $\overline r$ (see Definition~\ref{dfn:lift_maps}\;(iv)), and the third equality follows from the marginal constraint \eqref{eq:marginal_const}.\footnote{Since $(\mu_t^{\xi,a,\mathbb{P}}\,,\,{\Lambda}_t^{\xi,a,\mathbb{P}})\in \operatorname{gr}(\mathfrak{U})$, $\mathbb{P}$-a.s., for all $t\geq 0$, the term $\overline{r}(\mu_t^{\xi,a,\mathbb{P}}\,,\,{\Lambda}_t^{\xi,a,\mathbb{P}})$ is well-defined in the $\mathbb{P}$-a.s.\;sense.}
		
		Moreover, since $r$ is bounded and $\beta<1$ (see Assumption \ref{as:general}), by the dominated convergence theorem we can rewrite ${\cal J}^a(\xi)$ (given in \eqref{eq:worst_represent}) by 
		\begin{align}\label{eq:lift_worst}
			{\cal J}^a(\xi)=\inf_{\mathbb{P}\in {\cal Q}}\mathbb{E}^{\mathbb{P}}\bigg[\sum_{t=0}^\infty \beta^t \overline{r}(\mu_t^{\xi,a,\mathbb{P}}\,,\,{\Lambda}_t^{\xi,a,\mathbb{P}})\bigg].
		\end{align}
		\end{rem}

        Using Proposition \ref{pro:lift_dynamics}--particularly the MDP given in~\eqref{eq:lift_2} and the representations \eqref{eq:lift_reward} and \eqref{eq:lift_worst} in Remark \ref{rem:lift_gains}--we can view the robust MFC problem \eqref{eq:MFC_represent} as a robust MDP with state and action processes $(\mu_t^{\xi,a,\mathbb{P}},\Lambda_t^{\xi,a,\mathbb{P}})_{t\geq 0}$ given in \eqref{eq:MDP_state} and \eqref{eq:MDP_action}. This leads us to consider the following Bellman-Isaacs operator $\cal T$ defined on $C_{b}({\cal P}({S});\mathbb{R})$: for every $\overline {V}\in C_{b}({\cal P}({S});\mathbb{R})$ 
	\begin{align}\label{dfn:oper_T}
		{\cal T}\overline{V}(\mu):= \sup_{\Lambda \in \mathfrak{U}(\mu)} \bigg\{\overline{r}(\mu,\Lambda ) + \beta \inf_{p\in \mathfrak{P}^0} \int_{{\cal P}(S)} \overline V(\mu') \overline p(d\mu'|\mu,\Lambda,p)\bigg\}\quad\;\mu \in {\cal P}(S),
	\end{align}
	where $\mathfrak{P}^0$ is given in Definition \ref{dfn:measures}\;(i), and $\mathfrak{U}$, $\overline r$ and $\overline p$ are given in Definition \ref{dfn:lift_maps}.
	
	Following the framework of the `{local to global paradigm}' for robust MDP problems (see, e.g., \cite{neufeld2023markov,neufeld2024non,langner2024markov}), we first aim to characterize the local (i.e., one time-step) optimizers of the Bellman–Isaacs operator~${\cal T}$, and subsequently establish the fixed point theorem. This will then enable us to construct the global optimizers of the robust MFC problem \eqref{eq:MFC_represent}.
	
	To that end, we impose the following conditions on the basic components given in Definition~\ref{dfn:basic_element}. These conditions are (slightly) stronger than those in Assumption \ref{as:general}, as they contain certain regularity on the arguments in $A$ and $E^0$ along with others on the arguments in $S$ and ${\cal P}(S\times A)$.  However, they allow us to have some useful properties on the lifted functions and mappings given in Definition~\ref{dfn:lift_maps}, which are similar to and appear in a framework for robust MDP problems under model uncertainty (see, e.g., \cite{neufeld2023markov,neufeld2024non,langner2024markov}). 
	\begin{as}\label{as:general2} 
		The following conditions hold:
		\begin{itemize}[leftmargin=3.em]
			\item [(i)] The subset $\mathfrak{P}^0$ (see Definition \ref{dfn:measures}\;(i)) is compact.
			\item [(ii)] There is some $\overline C_{\operatorname{F}}>0$ such that\footnote{As noted in Section \ref{subsec:notation}, the product space $S \times A \times {\cal P}(S\times A)  \times E^0$ is endowed with the corresponding product topology induced by the following metric: for every $(s,a,\Lambda,e^0),(\tilde s,\tilde a,\tilde\Lambda,\tilde e^0)\in S \times A \times {\cal P}(S\times A) \times E^0$,
				\begin{align*}
					d_{S \times A \times {\cal P}(S\times A)  \times E^0}((s,a,\Lambda,e^0),(\tilde s,\tilde a,\tilde\Lambda,\tilde e^0)) :=d_S(s,\tilde s)+d_A(a,\tilde a)+ {\cal W}_{{\cal P}(S\times A)}(\Lambda,\tilde \Lambda)+d_{E^0}(e^0,\tilde e^0 ) .
				\end{align*}
				The same convention applies to $S\times A\times {\cal P}(S\times A)$ appearing in (iii).
			} for every $(s,a,\Lambda,e^0),(\tilde s,\tilde a,\tilde\Lambda,\tilde e^0)\in S \times A \times {\cal P}(S\times A)  \times E^0$ 
			\begin{align*}
				\hspace{1.5em} 
				\int_{E}d_S\big(\operatorname{F}(s,a,\Lambda,e,e^0),\operatorname{F}(\tilde s,\tilde a,\tilde\Lambda,e,\tilde e^0)\big)\lambda_{\varepsilon}(de)
				\leq \overline C_{{\operatorname{F}}} d_{S \times A \times {\cal P}(S\times A)  \times E^0}\big((s,a,\Lambda,e^0),(\tilde s,\tilde a,\tilde\Lambda,\tilde e^0)\big).
			\end{align*}
			\item [(iii)] The reward function $r$ is Lipschitz continuous, in the sense that there is some $\overline C_{{r}}>0$ such that for every $(s,a,\Lambda),(\tilde s,\tilde a,\tilde \Lambda)\in S\times A \times {\cal P}({S\times A})$
			\[
			\hspace{3.em} |r (s,a,\Lambda)-r (\tilde s,\tilde a,\tilde \Lambda)|\leq \overline C_{{r}} d_{S\times A \times {\cal P}(S\times A)}\big((s,a,\Lambda),(\tilde s,\tilde a,\tilde \Lambda)\big).
			\] 
			\item [(iv)] $\beta $ is in $[0,1\wedge (2\overline C_{{\operatorname{F}}})^{-1})$.
		\end{itemize} 
	\end{as}


	In the following proposition, we characterize the local optimizers of the Bellman-Isaacs operator~${\cal T}$ given in \eqref{dfn:oper_T}. To that end, we recall that given $L\geq 0$, $\operatorname{Lip}_{b,L}({\cal P}(S);\mathbb{R})\subset C_{b}({\cal P}(S);\mathbb{R})$ is the set of all $L$-Lipschitz continuous~functions defined on ${\cal P}(S)$. 
	\begin{pro}\label{pro:dpp}
		Suppose that Assumption \ref{as:general2}\;(i)--(iii) are satisfied. Then the following holds: For every $L\geq 0$ and every $\overline{V}\in\operatorname{Lip}_{b,L}({\cal P}(S);\mathbb{R})$,
		\begin{itemize}
			\item [(i)] (Local minimizer) There exists a measurable selector $\overline{p}^*:{\cal P}(S\times A)\ni \Lambda\mapsto \overline{p}^*(\Lambda)\in \mathfrak{P}^0$ such that for every $\Lambda\in   {\cal P}(S\times A)$ 
			\begin{align}\label{eq:local_min}
				\hspace{1.em}\int_{{\cal P}(S)} \overline V(\mu')\overline{p}(d\mu'|\operatorname{pj}_S(\Lambda),\Lambda,\overline{p}^*(\Lambda))= \inf_{p \in \mathfrak{P}^0}\int_{{\cal P}(S)} \overline V(\mu') \overline p(d\mu'|\operatorname{pj}_S(\Lambda),\Lambda,p).
			\end{align}
			\item [(ii)] (Local maximizer) There exists a measurable selector $\overline{\pi}^*:{\cal P}(S)\ni \mu \mapsto \overline{\pi}^*(\mu)\in \mathfrak{U}(\mu)$ satisfying that for every $\mu\in {\cal P}(S)$
			\begin{align}\label{eq:local_max}
				\hspace{1.em} \overline{r}(\mu,\overline \pi^*(\mu)) + \beta \inf_{p \in \mathfrak{P}^0} \int_{{\cal P}(S)} \overline V(\mu') \overline p(d\mu'|\mu,\overline \pi^*(\mu),p)={\cal T}\overline{V}(\mu).
			\end{align}
		\end{itemize}
	\end{pro}

	We now apply the Banach fixed-point theorem (see, e.g., \cite[Theorem A~3.5]{bauerle2011markov}) for the Bellman-Isaacs operator~${\cal T}$ given in \eqref{dfn:oper_T}.  
	\begin{pro}\label{pro:fixed_point}
		Suppose that Assumption \ref{as:general2} is satisfied, and let $\overline{L} \geq 2\overline C_r/(1-2\beta \overline C_{{\operatorname{F}}})$. Then it holds that ${\cal T}(\operatorname{Lip}_{b,\overline L}({\cal P}(S);\mathbb{R}) )\subseteq \operatorname{Lip}_{b,\overline L}({\cal P}(S);\mathbb{R})$,  and for every $\overline{V}^1,\overline{V}^2\in \operatorname{Lip}_{b,\overline L}({\cal P}(S);\mathbb{R})$
		\begin{align}\label{eq:contraction}
			\|{\cal T}\overline{V}^1- {\cal T}\overline{V}^2\|_{\infty}\leq \beta  \|\overline{V}^1- \overline{V}^2\|_{\infty}.
		\end{align}
		In particular, there exists a unique $\overline{V}^* \in \operatorname{Lip}_{b,\overline L}({\cal P}(S);\mathbb{R})$ satisfying that ${\cal T}\overline{V}^*=\overline{V}^*$. Moreover, it holds for every $\overline{V}\in \operatorname{Lip}_{b,\overline L}({\cal P}(S);\mathbb{R})$ that $\overline{V}^*= \lim_{n\to \infty} {\cal T}^n \overline{V}$.
	\end{pro}

	The proofs of Propositions~\ref{pro:dpp} and \ref{pro:fixed_point} can be found in Section \ref{proof:subsec:lift_MDP}.
	
	\subsection{Verification theorem}\label{subsec:main_thm}
	This section aims to establish that the fixed point $\overline{V}^*$ of the Bellman-Isaacs operator ${\cal T}$ (see Proposition \ref{pro:fixed_point}) coincides with the robust MFC problem $V$ of the representative agent (see \eqref{eq:MFC_represent}) in the sense that\footnote{\label{footnote:invar_xi}By construction of the set ${\cal Q}$ (see Definition \ref{dfn:measures}\;(ii)), the law of $\xi\in L_{{\cal F}_0}^0(S)$ is invariant w.r.t.\;the choice of supporting probability measure $\mathbb{P}\in {\cal Q}$. Therefore, we can and do write $\mathscr{L}(\xi):=\mathscr{L}_{\mathbb{P}}(\xi)\in {\cal P}(S)$ for any $\mathbb{P}\in{\cal Q}$.} $V(\xi)=\overline{V}(\mathscr{L}(\xi))$ for all $\xi\in L_{{\cal F}_0}^0(S)$. 
	
	To that end, we first construct a measure in ${\cal Q}$ for each open-loop control in ${\cal A}$ (see~\eqref{eq:open_loop_set}), using the local minimizer from Proposition \ref{pro:dpp}\;(i). This will later be used in the verification theorem to derive a worst-case measure in ${\cal Q}$ by suitably choosing an optimal control in ${\cal A}$.
    
	\begin{lem}\label{lem:worst_lift_dynamics}
		Suppose that Assumption \ref{as:general2} is satisfied. 
		Let $\xi \in L_{{\cal F}_0}^0(S)$ be the initial state of the representative agent. Then for every $a\in {\cal A}$, there exists $\underline {\mathbb{P}}^{\xi,a}\in {\cal Q}$ induced by some $(\underline{p}_t^{\xi,a})_{t\geq 1}\in \mathcal{K}^0$ (see Definition~\ref{dfn:measures})
		such that $\underline{\mathbb{P}}^{\xi,a}$-a.s.
		\begin{align}\label{eq:worst_mdp}
			\begin{aligned}
			\hspace{3.em} 
			&\mathscr{L}_{\underline {\mathbb{P}}^{\xi,a}}(\varepsilon_1^0)=
			\underline{p}_1^{\xi,a} = \overline{p}^*(\underline{ \Lambda}_{0}^{\xi,a}),\\
			&\mathscr{L}_{\underline {\mathbb{P}}^{\xi,a}}(\varepsilon_{t+1}^0\,| {\cal F}_{t}^0)= 
			\underline{p}_{t+1}^{\xi,a}(\,\cdot\,| \varepsilon_{1:t}^0)=\overline{p}^*(\underline{ \Lambda}_{t}^{\xi,a})\quad \mbox{for all $t\geq 1$},
			\end{aligned}
		\end{align}
		where $\overline p^*$ is the local minimizer given in Proposition \ref{pro:dpp}\;(i), $\underline{ \Lambda}_{0}^{\xi,a}$ is the joint law of $(s_{0}^{\xi,a,\underline{\mathbb{P}}^{\xi,a}},a_0)$ under $\underline {\mathbb{P}}^{\xi,a}$, and for $t\geq 1$ $\underline{ \Lambda}_{t}^{\xi,a}$ is the conditional joint law of $(s_{t}^{\xi,a,\underline{\mathbb{P}}^{\xi,a}},a_t)$ under $\underline {\mathbb{P}}^{\xi,a}$ given $\varepsilon^0_{1:t}$. 
		 Consequently, we have 
		\begin{align}\label{eq:worst_mdp2}
				\mathscr{L}_{\underline{\mathbb{P}}^{\xi,a}}(\underline{\mu}_{t+1}^{\xi,a})=\overline p(\,\cdot\,\big| \operatorname{pj}_S(\underline{\Lambda}_{t}^{\xi,a}), \underline{\Lambda}_{t}^{\xi,a}, \overline{p}^*(\underline{\Lambda}_{t}^{\xi,a})),\quad  \mbox{$\underline{\mathbb{P}}^{\xi,a}$-a.s.,}\quad \mbox{for all $t\geq 0$},
		\end{align}
			where $\overline p$ is given in Definition \ref{dfn:lift_maps}, and $\underline{\mu}_{t+1}^{\xi,a}$ is the conditional law of $s_{t+1}^{\xi,a,\underline{\mathbb{P}}^{\xi,a}}$ under $\underline {\mathbb{P}}^{\xi,a}$ given~$\varepsilon^0_{1:t+1}$.
	\end{lem}

	We now construct an open-loop control in ${\cal A}$, using the local maximizer from Proposition\;\ref{pro:dpp}\;(ii). Then we will verify that this open-loop control is indeed a maximizer of the robust MFC problem given in \eqref{eq:MFC_represent}.
    
    We impose the following condition.
	\begin{as}\label{as:randomize}
		$\lambda_\vartheta\in {\cal P}(\Theta)$ given in Definition \ref{dfn:measures} is atomless.
	\end{as}

	\begin{rem}\label{rem:randomize} Assumption \ref{as:randomize} also appears in \cite{carmona2023model} (see Section~2.1.2). Moreover, \cite{motte2022mean} incorporates this assumption by assuming the existence of a uniform random variable that is independent of the given initial state (see Section~3 therein). This assumption is crucial for constructing an optimal control/policy from the lifted dynamic programming results presented in both references---and consequently in this article as well. In particular, we often use the following properties.
		
	Since $\mathscr{L}_{\mathbb{P}}(\vartheta)=\lambda_\vartheta$ for all $\mathbb{P}\in {\cal Q}$ (see Remark \ref{rem:well_dfn_1}\;(ii)), Assumption \ref{as:randomize} implies the existence of a sequence $(h_t)_{t\geq 0}$ of Borel measurable functions $h_t:\Theta \to [0,1]$ such that under any $\mathbb{P}\in {\cal Q}$, 
		\[
		(h_t(\vartheta_t))_{t\geq 0}\;\;\mbox{is i.i.d.\;with law ${\cal U}_{[0,1]}$, }
		\]
		i.e., uniform distribution on $[0,1]$; see \cite[Theorem~9.2.2]{bogachev2007constructions}.   Since all the agents are indistinguishable, such a sequence exists for each agent $i\in \mathbb{N}$, and we denote it by~$(h^i_t)_{t\geq 0}$. 
	\end{rem}

	\begin{lem}\label{lem:robust_cont}
		Suppose that Assumptions \ref{as:general2} and \ref{as:randomize} are satisfied. Let $\xi \in L_{{\cal F}_0}^0(S)$ be the initial state of the representative agent. Then there exists $a^*\in {\cal A}$ such that for every $\mathbb{P}\in {\cal Q}$,
		\begin{align}\label{eq:law_control}
			{ \Lambda}_{t}^{\xi,a^*,\mathbb{P}}= \overline \pi^*({\mu}_t^{\xi,a^*,\mathbb{P}}),\quad \mbox{ $\mathbb{P}$-a.s., for all $t\geq 0$,}
		\end{align}
		where $\overline \pi^*$ is the local maximizer given in Proposition \ref{pro:dpp}\;(ii), and $({ \mu}_{t}^{\xi,a^*,\mathbb{P}})_{t\geq 0}$ and $({\Lambda}_{t}^{\xi,a^*,\mathbb{P}})_{t\geq 0}$ are given in \eqref{eq:MDP_state} and \eqref{eq:MDP_action}, respectively, under $(a^*,\mathbb{P})$.
	\end{lem}
	
	We are now ready to state the verification theorem for the constructed open-loop control and probability measure in the preceding two lemmas.
	The proofs of the theorem and preceding lemmas are provided in Section~\ref{proof:subsec:main_thm}.
	\begin{thm}\label{thm:origin_MDP}
		Suppose that Assumptions \ref{as:general2} and \ref{as:randomize} are satisfied. Let $\overline{L} \geq 2C_r/(1-2\beta C_{{\operatorname{F}}})$ be given, and  let $\overline{V}^* \in \operatorname{Lip}_{b,\overline{L}}({\cal P}(S);\mathbb{R})$ be such that ${\cal T}\overline V^*=\overline V^*$ (see Proposition \ref{pro:fixed_point}). Moreover, let $a^*\in {\cal A}$ be such that 
		\eqref{eq:law_control} holds for every $\mathbb{P}\in {\cal Q}$ (see Lemma \ref{lem:robust_cont}). 
		Moreover, let ${\cal J}^{a^*}$  and $V$ be given in \eqref{eq:worst_represent} and \eqref{eq:MFC_represent}, respectively. Then, 
        for every $\xi \in L_{{\cal F}_0}^0(S)$ the following hold:
        \begin{itemize}
            \item [(i)] $\overline{V}^*(\mathscr{L}(\xi))= V(\xi)$, where $\mathscr{L}(\xi)\in {\cal P}(S)$ is the law of $\xi$ (see Footnote \ref{footnote:invar_xi}). 
            \item [(ii)] $a^*\in {\cal A}$ and $\underline{\mathbb{P}}^{\xi,a^*}\in {\cal Q}$ induced by $(\underline{p}_t^{\xi,a^*})_{t\geq 1} \in \mathcal{K}^0$ satisfying \eqref{eq:worst_mdp},\;\eqref{eq:worst_mdp2} (see Lemma\;\ref{lem:worst_lift_dynamics}) are optimal in the sense that 
            \begin{align}\label{eq:verify}
			{V}(\xi)
			={\cal J}^{a^*}(\xi)
			= \mathbb{E}^{\underline{\mathbb{P}}^{\xi,a^*}}\big[R^{a^*,\underline{\mathbb{P}}^{\xi,a^*}}(\xi) \big].
		\end{align}
        \end{itemize}
	\end{thm}
	
	\begin{rem}
		As a consequence of Theorems \ref{thm:propagation_chaos} and \ref{thm:origin_MDP}, under Assumptions \ref{as:general2} and \ref{as:randomize} the optimal open-loop policy $\pi^*\in \Pi$ of the robust MFC problem $V$ (see Definition~\ref{dfn:MFC})---which can be obtained from the optimal open-loop control $a^*\in {\cal A}$ in Theorem \ref{thm:origin_MDP} of the representative robust MFC problem $V(\xi)$ in \eqref{eq:MFC_represent}---serves as an approximate of the $N$-agent optimization problem $V^N$ (see Definition \ref{dfn:N_model}) when $N\in \mathbb{N}$ is sufficiently large.
		
		Lastly, we note that computing the local optimizers from the lifted dynamic programming principle (given in Proposition \ref{pro:dpp}) is crucial for deriving the optimal open-loop control of the robust MFC problem. In particular, this step involves implementation of $Q$-learning (or policy iteration) algorithms for the lifted dynamic programming principle and analyzing their convergence, together with the discretization error arising from of the lifted state and action spaces. While we defer these aspects to future research, in Section \ref{sec:numerics} we present some numerical examples based on a value iteration type scheme to implement the lifted dynamic programming principle. 
	\end{rem}

	\subsection{Connection with a closed-loop Markov policy framework}\label{subsec:closed_Markov} In this section, we introduce the notion of a closed-loop Markov policy for the robust MFC problem. In particular, 
	following \cite[Definition 10]{carmona2023model},  we consider a relaxed version of the robust MFC problem in Definition~\ref{dfn:MFC}, in which individual agents are allowed to sample their actions randomly according to a policy specified by the social planner.
	
	As in Sections \ref{subsec:lift_MDP} and \ref{subsec:main_thm}, we suppress the index $i\in \mathbb{N}$ representing individual agents and consider the following representation agent's robust MFC problem with closed-loop Markov policies.
	\begin{dfn}\label{dfn:MFC_closed} Let ${\cal Q}$ be the uncertainty measures set given in Definition~\ref{dfn:measures}. Moreover, let $\mathbb{F},\mathbb{G}$ be the filtrations given in~\eqref{eq:filtrations_represent}, and let $\mathbb{F}^0$ be the filtration generated by the common noise.
		\begin{itemize}[leftmargin=3.em]
			\item [(i)] Denote by $\Pi^{c}$ the set of all closed-loop Markov policies $\pi^c:=(\pi_t^{{c}})_{t\geq 0}$ such that for every $t\geq 0$ the~kernel 
			\[
			\pi^{c}_t:S\times {\cal P}(S)\ni (s,\mu)\mapsto \pi^{c}_t(da|s,\mu) \in {\cal P}(A)
			\]
			induces a randomized action given a couple of a state and a probability measure~on~$S$. 
			\item [(ii)] Let $\xi \in L^0_{{\cal F}_0}(S)$ be the fixed initial state. Assume that for any $(\pi^c,\mathbb{P})\in \Pi^c\times  {\cal Q}$, the~state and action processes $(s_{t}^{\xi,\pi^c,\mathbb{P}},a_{t}^{\pi^c,\mathbb{P}})_{t\geq 0}$ for the representative agent in the inifinite population model satisfy that\footnote{We refer to Remark \ref{rem:well-defined-closed-loop} for the well-posedness of $(s_{t}^{\xi,\pi^c,\mathbb{P}},a_{t}^{\pi^c,\mathbb{P}})_{t\geq 0}$ defined as in Definition \ref{dfn:MFC_closed}\;(ii).} $(s_{t}^{\xi,\pi^c,\mathbb{P}})_{t\geq 0}$ is $\mathbb{F}$-adapted, $(a^{\pi^c,\mathbb{P}}_t)_{t\geq 0}$ is $\mathbb{G}$-adapted, and they satisfy
			\begin{align}\label{eq:MKV_represent_closed}
				\begin{aligned}
					\hspace{1.em} 
					&s_{t+1}^{\xi,\pi^c,\mathbb{P}}:=\operatorname{F}(s^{\xi,\pi^c,\mathbb{P}}_t,a_t^{\pi^c,\mathbb{P}},\mathbb{P}^{0}_{(s^{\xi,\pi^c,\mathbb{P}}_t,a_t^{\pi^c,\mathbb{P}})},\varepsilon_{t+1},\varepsilon_{t+1}^0)\quad \mbox{for $t\geq 0$},\;\; \mbox{with $\;\;s_0^{\xi,\pi^c,\mathbb{P}}:=\xi,$}\\
					&\mathscr{L}_{\mathbb{P}}(a^{\pi^c,\mathbb{P}}_t|{\cal F}_t)=\pi^c_t(\,\cdot\,|s^{\xi,\pi^c,\mathbb{P}}_t,\mathbb{P}_{s^{\xi,\pi^c,\mathbb{P}}_t}^0)\quad \mbox{$\mathbb{P}$-a.s.\quad for $t\geq 0$},
				\end{aligned}
			\end{align}
			where $\mathbb{P}^{0}_{(s^{\xi,\pi^c,\mathbb{P}}_t,a_t^{\pi^c,\mathbb{P}})}$ is the conditional joint law of $(s^{\xi,\pi^c,\mathbb{P}}_t,a_t^{\pi^c,\mathbb{P}})$ under $\mathbb{P}$ given $\varepsilon^0_{1:t}$ for~$t\geq 1$, with the convention that $\mathbb{P}^{0}_{(s^{\xi,\pi^c,\mathbb{P}}_0,a_0^{\pi^c,\mathbb{P}})}:=\mathscr{L}_{\mathbb{P}}((s^{\xi,\pi^c,\mathbb{P}}_0,a_0^{\pi^c,\mathbb{P}}))$. In analogy, $\mathbb{P}^{0}_{s^{\xi,\pi^c,\mathbb{P}}_t}$ is the conditional law of $s^{\xi,\pi^c,\mathbb{P}}_t$ under $\mathbb{P}$ given $\varepsilon^0_{1:t}$ for $t\geq 1$ with $\mathbb{P}^{0}_{s^{\xi,\pi^c,\mathbb{P}}_0}:=\mathscr{L}_{\mathbb{P}}(s^{\xi,\pi^c,\mathbb{P}}_0)$.
			
			\item [(iii)] Accordingly, the robust MFC problem under closed-loop Markov policies is 
			\begin{align}\label{eq:MFC_represent_closed}
				V^c(\xi):=\sup_{\pi^c\in \Pi^c} {\cal J}^{\pi^c}(\xi),\quad\;\xi \in L_{{\cal F}_0}^0(S),
			\end{align}
			where ${\cal J}^{\pi^c}(\xi)$ is defined as ${\cal J}^{\pi^c}(\xi):=\inf_{\mathbb{P}\in {\cal Q}}\mathbb{E}^{\mathbb{P}}[R^{\pi^c,\mathbb{P}}(\xi) ]$ with
			\begin{align*}
				R^{\pi^c,\mathbb{P}}(\xi):=\mathbb{E}^{\mathbb{P}^0}\bigg[\sum_{t=0}^\infty \beta^t r(s_{t}^{\xi,\pi^c,\mathbb{P}},a^{\pi^c,\mathbb{P}}_{t},\mathbb{P}^0_{(s^{\xi,\pi^c,\mathbb{P}}_t,a^{\pi^c,\mathbb{P}}_t)})\bigg].
			\end{align*}
		\end{itemize}
	\end{dfn}
	\begin{rem}\label{rem:well-defined-closed-loop}
		Under Assumption \ref{as:randomize}, the conditional McKean-Vlasov dynamics with closed-loop Markov policies, as given in Definition \ref{dfn:MFC_closed}\;(ii), are well-defined. Indeed, by using the random variable $h_t(\vartheta_t)\sim {\cal U}_{[0,1]}$ (see Remark~\ref{rem:randomize}) and the Blackwell--Dubins function $\rho_A:{\cal P}(A)\times[0,1]\to A$ (see Lemma \ref{lem:BlackDubins}), we can define, for any $\pi^c\in \Pi^c$ and $\mathbb{P}\in {\cal Q}$, 
		\[
		a_t^{\pi^c,\mathbb{P}}:=\rho_A\big(\pi_t^c(\,\cdot\,|\,s_t^{\xi,\pi^c,\mathbb{P}},\mathbb{P}^0_{s_t^{\xi,\pi^c,\mathbb{P}}}),h_t(\vartheta_t)\big)\quad t\geq 0.
		\]
		By the same arguments presented for the proof of Lemma \ref{lem:measurability}\;(ii), we note that $s_t^{\xi,\pi^c,\mathbb{P}}$ is ${\cal F}_t$~measurable and $\mathbb{P}^0_{s_t^{\xi,\pi^c,\mathbb{P}}}$ is ${\cal F}_t^0$ measurable. Consequently, $a_t^{\pi^c,\mathbb{P}}$ is ${\cal G}_t$ measurable by the construction above.  Furthermore, since ${\cal F}_t$ is independent of $\vartheta_t$, the property of $\rho_A$ ensures that $a_t^{\pi^c,\mathbb{P}}$ satisfies the distributional constraint given in \eqref{eq:MKV_represent_closed}.  
	\end{rem}

	We aim to show that the robust MFC problem $V^c$ given in \eqref{eq:MFC_represent_closed} coincides with the open-loop robust MFC problem $V$ given in \eqref{eq:MFC_represent}.  This equivalence will be established by demonstrating that $V^c(\xi)=\overline{V}^*(\mathscr{L}(\xi))$ 
	for all $\xi \in L_{{\cal F}_0}^0(S)$, where $\overline{V}^*$ is the fixed point of the Bellman–Isaacs operator $\mathcal{T}$ given in Proposition \ref{pro:fixed_point}, and $\mathscr{L}(\xi)\in {\cal P}(S)$ is the law of $\xi$ (see Footnote \ref{footnote:invar_xi}).
	
	To this end, and following the approach in Section \ref{subsec:lift_MDP},  we begin by examining the dynamics of the lifted state and action processes, defined as follows:  for every $\pi^c\in \Pi^c$ and $\mathbb{P}\in {\cal Q}$, 
	\begin{align}\label{eq:MDP_state_closed}
		\begin{aligned}
			(\mu_t^{\xi,\pi^c,\mathbb{P}})_{t\geq 0}&:=(\mathbb{P}^0_{s_t^{\xi,\pi^c,\mathbb{P}}})_{t\geq 0}\subseteq {\cal P}(S),\\
			(\Lambda_t^{\xi,\pi^c,\mathbb{P}})_{t\geq 0}&:=(\mathbb{P}^0_{(s_t^{\xi,\pi^c,\mathbb{P}},a_t^{\pi^c,\mathbb{P}})})_{t\geq 0}\subseteq {\cal P}(S\times A).
		\end{aligned}
	\end{align}
	Here we note that both processes are $\mathbb{F}^0$ adapted (see Lemma \ref{lem:measurability}).
	\begin{lem}\label{lem:MDP_closed_loop}
		Suppose that Assumptions \ref{as:general2} and \ref{as:randomize} are satisfied. Let $\pi^c\in \Pi^c$ be given and let $\mathbb{P}\in {\cal Q}$ be {induced by} some $(p_t)_{t\geq 1}\in \mathcal{K}^0$ (see Definition \ref{dfn:measures}). Then, 
		\begin{align}\label{eq:kernel_prod}
		\Lambda_t^{\xi,\pi^c,\mathbb{P}} =\mu_t^{\xi,\pi^c,\mathbb{P}}\mathbin{\hat \otimes} \pi_t^c(\,\cdot\,|\,\cdot,\mu_t^{\xi,\pi^c,\mathbb{P}})\quad\mbox{$\mathbb{P}$-a.s. for all $t\geq 0$}.
		\end{align}
		Consequently, it holds that $\mathbb{P}$-a.s.
		\begin{align}
			\begin{aligned}
				\mathscr{L}_{\mathbb{P}}(\mu_{1}^{\xi,\pi^c,\mathbb{P}})&=\overline p(\,\cdot\,|\,\mu_0^{\xi,\pi^c,\mathbb{P}},\,\,\mu_0^{\xi,\pi^c,\mathbb{P}}\mathbin{\hat \otimes} \pi_0^c(\,\cdot\,|\,\cdot,\mu_0^{\xi,\pi^c,\mathbb{P}}),\,p_1(\cdot)),\\
				\mathscr{L}_{\mathbb{P}}(\mu_{t+1}^{\xi,\pi^c,\mathbb{P}})&=\overline p(\,\cdot\,|\,\mu_t^{\xi,\pi^c,\mathbb{P}},\,\,\mu_t^{\xi,\pi^c,\mathbb{P}}\mathbin{\hat \otimes} \pi_t^c(\,\cdot\,|\,\cdot,\mu_t^{\xi,\pi^c,\mathbb{P}}),\,p_{t+1}(\,\cdot\,|\varepsilon^0_{1:t}))\quad \mbox{for all $t\geq 1$}. \label{eq:lift_2_closed}
			\end{aligned}
		\end{align}
	\end{lem}

	Then, as in Lemma \ref{lem:worst_lift_dynamics}, we construct a measure in ${\cal Q}$ for each closed-loop policy in $\Pi^c$ (see Definition \ref{dfn:MFC_closed}), using the local minimizer from Proposition \ref{pro:dpp}\;(i). 
	\begin{lem}\label{lem:worst_lift_dynamics_closed}
		Suppose that Assumptions \ref{as:general2} and \ref{as:randomize} are satisfied.
		For every $\pi^c\in \Pi^c$, there exists $\underline {\mathbb{P}}^{\xi,\pi^c}\in {\cal Q}$ induced by some $(\underline{p}_t^{\xi,\pi^c})_{t\geq 1}\in \mathcal{K}^0$ (see Definition~\ref{dfn:measures})
		such that $\underline{\mathbb{P}}^{\xi,\pi^c}$-a.s.
		\begin{align}\label{eq:worst_mdp_closed}
			\begin{aligned}
				\hspace{3.em} 
				&\mathscr{L}_{\underline {\mathbb{P}}^{\xi,\pi^c}}(\varepsilon_1^0)=
				\underline{p}_1^{\xi,\pi^c} = \overline{p}^*(\underline{ \Lambda}_{0}^{\xi,\pi^c}),\\
				&\mathscr{L}_{\underline {\mathbb{P}}^{\xi,\pi^c}}(\varepsilon_{t+1}^0\,|\,{\cal F}_{t}^0)= 
				\underline{p}_{t+1}^{\xi,\pi^c}(\cdot\,|\,\varepsilon_{1:t}^0)=\overline{p}^*(\underline{ \Lambda}_{t}^{\xi,\pi^c})\quad \mbox{for all $t\geq 1$},
			\end{aligned}
		\end{align}
		where $\overline p^*$ is the local minimizer in Proposition \ref{pro:dpp}\;(i), 
		$\underline{ \Lambda}_{0}^{\xi,\pi^c}$ is the joint law of $(s_{0}^{\xi,\pi^c,\underline{\mathbb{P}}^{\xi,\pi^c}},a^{\pi^c,\underline {\mathbb{P}}^{\xi,\pi^c}}_0)$ under $\underline{\mathbb{P}}^{\xi,\pi^c}$, and for $t\geq 1$ $\underline{ \Lambda}_{t}^{\xi,\pi^c}$ is the conditional joint law of $(s_{t}^{\xi,\pi^c,\underline{\mathbb{P}}^{\xi,\pi^c}},a^{\pi^c,\underline {\mathbb{P}}^{\xi,\pi^c}}_t)$  under $\underline{\mathbb{P}}^{\xi,\pi^c}$ given $\varepsilon^0_{1:t}$. 
		Consequently, we have
		\begin{align}\label{eq:worst_mdp2_closed}
			\mathscr{L}_{\underline{\mathbb{P}}^{\xi,\pi^c}}(\underline{\mu}_{t+1}^{\xi,\pi^c})=\overline p(\,\cdot\,| \operatorname{pj}_S(\underline{\Lambda}_{t}^{\xi,\pi^c}), \underline{\Lambda}_{t}^{\xi,\pi^c}, \overline{p}^*(\underline{\Lambda}_{t}^{\xi,\pi^c})),\quad \mbox{ $\underline{\mathbb{P}}^{\xi,\pi^c}$-a.s., for all $t\geq 0$,}
		\end{align}
		where $\overline p$ is given in Definition \ref{dfn:lift_maps}, and $\underline{\mu}_{t+1}^{\xi,\pi^c}$ is the conditional law of $s_{t+1}^{\xi,\pi^c,\underline{\mathbb{P}}^{\xi,\pi^c}}$ given $\varepsilon^0_{1:t+1}$.
	\end{lem}
	The proofs of Lemma \ref{lem:MDP_closed_loop} and \ref{lem:worst_lift_dynamics_closed} are presented in Section \ref{proof:subsec:closed_Markov}.
	\begin{rem}
		While the construction of $(\underline{p}_t^{\xi,\pi^c})_{t\geq 1}\in \mathcal{K}^0$ given in Lemma \ref{lem:worst_lift_dynamics_closed} proceeds inductively (as in the proof of Lemma \ref{lem:worst_lift_dynamics}), the arguments differ from those used therein. This is due to the fact that a closed-loop Markov policy $\pi^c\in \Pi^c$ does not determine a fixed action process, but a randomly sampled one. For this, we rely on the Blackwell-Dubins function given in Lemma\;\ref{lem:BlackDubins} together with Remark \ref{rem:randomize} and some measure-theoretic arguments.
	\end{rem}

	Finally, we conclude that the robust MFC problem under the closed-loop Markov policy framework coincides with the fixed point $\overline V$, and hence with the robust MFC problem under the open-loop policy framework.
	\begin{cor}\label{cor:origin_MDP_closed}
		Suppose that Assumptions \ref{as:general2} and \ref{as:randomize} are satisfied. Let $\overline{L} \geq 2C_r/(1-2\beta C_{{\operatorname{F}}})$ be given, and let $\overline{V}^* \in \operatorname{Lip}_{b,\overline{L}}({\cal P}(S);\mathbb{R})$ be such that ${\cal T}\overline V^*=\overline V^*$ (see Proposition \ref{pro:fixed_point}).  Define  
		\begin{align}\label{eq:local_max_ext}
			{\pi}_{\operatorname{loc}}^{c,*}:S\times{\cal P}(S)\ni (s,\mu)\mapsto {\pi}_{\operatorname{loc}}^{c,*}(\,\cdot\,|s,\mu):={\cal K}_{S\times A}(\,\cdot\,|s, \overline\pi^*(\mu),\mu)\in{\cal P}(A), 
		\end{align}
		i.e., the universal disintegration kernel of $\overline{\pi}^*(\mu)$ w.r.t.\;$\operatorname{pj}_S(\overline{\pi}^*(\mu))=\mu$ (see Lemma \ref{lem:univ_disint}) so that 
		\begin{align}\label{eq:local_max_ext_integ}
			\overline\pi^*(\mu)=\mu  \mathbin{\hat \otimes} {\pi}_{\operatorname{loc}}^{c,*}(\,\cdot\,|\,\cdot,\mu).
		\end{align}
		Define $\pi^{c,*}:=(\pi_t^{c,*})_{t\geq 0}\in \Pi^c $ by $\pi_t^{c,*}:= \pi_{\operatorname{loc}}^{c,*}$ for every $t\geq 0$ (i.e., stationary closed-loop Markov policy). 
		Moreover, let $V^c$ and ${\cal J}^{\pi^{c,*}}$ be given in \eqref{eq:MFC_represent_closed}, and let $V$ be given in \eqref{eq:MFC_represent}. Then, for every $\xi \in L_{{\cal F}_0}^0(S)$ the following hold:
        \begin{itemize}
            \item [(i)] $\overline{V}^*(\mathscr{L}(\xi))= V^{c}(\xi)=V(\xi)$, where $\mathscr{L}(\xi)\in {\cal P}(S)$ is the law of $\xi$ (see Footnote \ref{footnote:invar_xi}).
            \item [(ii)] $\pi^{c,*}\in \Pi^c$ and $\underline{\mathbb{P}}^{\xi,\pi^{c,*}}$ induced by $(\underline{p}_t^{\xi,\pi^{c,*}})_{t\geq 1}\in \mathcal{K}^0$ satisfying \eqref{eq:worst_mdp_closed},\;\eqref{eq:worst_mdp2_closed} (see Lemma\;\ref{lem:worst_lift_dynamics_closed}) are optimal in the sense that
            \begin{align}\label{eq:verify_closed}
			{V}^c(\xi)
			={\cal J}^{\pi^{c,*}}(\xi)
			= \mathbb{E}^{\underline{\mathbb{P}}^{\xi,\pi^{c,*}}}\big[R^{\pi^{c,*},\underline{\mathbb{P}}^{\xi,\pi^{c,*}}}(\xi) \big].
		\end{align}
        \end{itemize}
	\end{cor}
	
	
	\section{Numerical examples}
	\label{sec:numerics}
	In this section, we apply our robust MFC framework under common noise uncertainty to illustrative examples in distribution matching and financial systemic risk, thereby emphasizing the critical role of incorporating common noise uncertainty into the analysis.
	In both examples, the algorithm implementing the lifted dynamic programming principle in Proposition \ref{pro:dpp} together with the verification theorem in Theorem \ref{thm:origin_MDP} (or Corollary \ref{cor:origin_MDP_closed}) builds upon the value iteration algorithm for the robust MDP framework of \cite[Section 4.4.1]{neufeld2023markov}.
	
	\subsection{Example 1: Distribution matching}\label{sec:ex1}
	We first consider an example inspired by Example~1 in~\cite{carmona2023model}, in which the goal for the central planner is to make the population distribution match a given target distribution. Common noise makes the task harder because it may randomly shift the distribution.  
	
	To be specific, consider the following basic elements (recall Definition~\ref{dfn:basic_element}):\footnote{The code is provided for the sake of completeness at \url{https://github.com/mlauriere/RobustMFMDP}.}
	\begin{itemize}[leftmargin=3.em]
		\item $S = \{1,2,\dots, |S|\}$ representing a one-dimensional grid world with $|S|$ states; in the experiments, we use $|S| =7$ states.
		\item $A = \{-1,0,1\}$, where the actions are interpreted respectively as moving to the left, staying or moving to the right.
		\item $E = \{0\}$, which means that there is no idiosyncratic noise.
		\item $E^0 = \{-1,0,1\}$, where the common noise values are interpreted as the actions but they affect the whole population.
		\item $\operatorname{F}:S \times A \times {\cal P}(S\times A) \times E \times E^0\to S$ is given by 
		\[
			\operatorname{F}(s,a,\Lambda,e,e^0) = \max(1, \min(|S|,s+a+e^0)),
		\]
		which represents the fact that the agent's movement is determined by her action and the common noise, and the agent remains at $1$ (resp.\;$7$) if she tries to move to the left (resp.\;right) of this state.
		\item $r:S\times A\times {\cal P}(S\times A)\to \mathbb{R}$ is given by 
		\[
			r(s,a,\Lambda) = \|\operatorname{pj}_{S}(\Lambda) - \mu^*\|_2^2 = \sum_{s \in S} |\operatorname{pj}_{S}(\Lambda)(s) - \mu^*(s)|^2,
		\]
		where $\mu^* \in {\cal P}(S)$ is a fixed target distribution which is part of the model's definition.
		\item {$\beta = 0.4$ is the discount factor so that Assumptions~\ref{as:general}\;(iii) and~\ref{as:general2}\;(iv) are satisfied.}  
	\end{itemize}

	For the common noise probability measure, we consider the following situation:
	\begin{itemize}
		\item The true common noise distribution $p_{\textrm{true}}\in {\cal P}(E^0)$ is given by 
		\begin{align}\label{eq:true}
			p_{\textrm{true}}:=v_{{\textrm{true}},1}\delta_{\{\varepsilon^0=-1\}}+v_{{\textrm{true}},2}\delta_{\{\varepsilon^0=0\}}+v_{{\textrm{true}},3}\delta_{\{\varepsilon^0=1\}},
		\end{align}
		with some probability vector $v_{{\textrm{true}}}:=(v_{{\textrm{true}},1},v_{{\textrm{true}},2},v_{{\textrm{true}},3})\in [0,1]^3$, i.e., a simplex.
		\item However, we consider that the central planner does not  know this true distribution; she has estimated the common noise distribution to be approximately equal to a reference probability measure $p_{\textrm{ref}}\in{\cal P}(E^0)$ with the corresponding probability vector $v_{\textrm{ref}}\in [0,1]^3$.
	\end{itemize}

	As a baseline, the central planner can learn a policy $\pi_{\textrm{ref}}$ which is optimal for the MFC model with common noise distribution $p_{\textrm{ref}}$. Alternatively, she can solve the robust MFC problem and learn a policy $\pi_{\textrm{robust}}$ which may be suboptimal for the model with $p_{\textrm{ref}}$ but which performs better than $\pi_{\textrm{ref}}$ in the true model with common noise distribution $p_{\textrm{true}}$.

	We consider the 
	uncertainty set $\mathfrak{P}^0$ which consists of all perturbed measures $p\in {\cal P}(E^0)$ of the reference measure $p_{\textrm{ref}}$,  whose corresponding probability vector $v\in [0,1]^3$ is 
	\begin{align}\label{eq:perturb}
		v:=\mathrm{renorm}(\max(0, v_{\textrm{ref}} + v_{\textrm{perturb}})),
	\end{align}
	where $v_{\textrm{perturb}}\in \mathbb{R}^3$ is a perturbation vector constructed as follows: each coordinate is sampled uniformly from $[-\delta_{\textrm{perturb}}, \delta_{\textrm{perturb}}]$, with a small $\delta_{\textrm{perturb}}>0$ representing the uncertainty level. The average of the $3$ coordinates is then subtracted to each coordinate to ensure that the average of $v_{\textrm{perturb}}$ over coordinates is $0$. Under this construction, Assumption \ref{as:general2}\;(i) is satisfied.  
	
	\begin{figure}[t]
		\centering
		\includegraphics[width=0.5\linewidth]{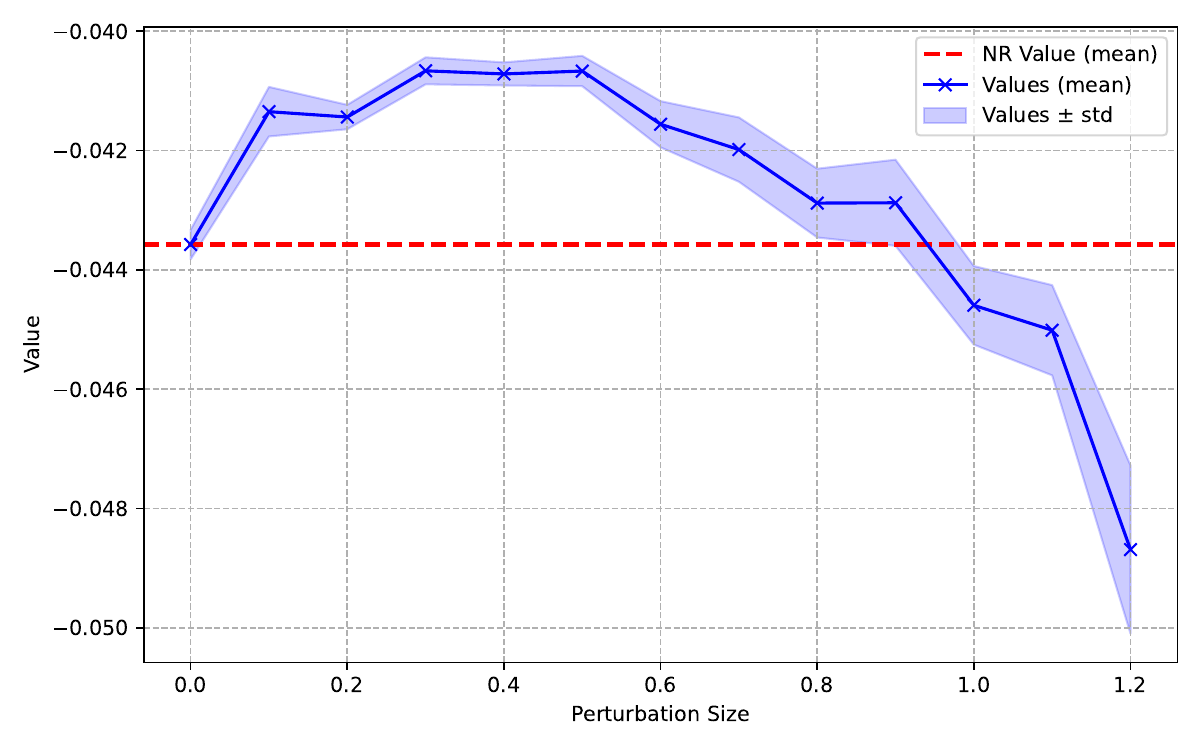}
		\caption{Values achieved under $p_{\textrm{true}}$ when using the optimal policy for the MFC under $p_{\textrm{ref}}$ (red dashed line) or the robust MFC under the uncertainty level $\delta_{\textrm{perturb}} \in \{0.0, 0.1, 0.2, 0.3, 0.4, 0.5, 0.6, 0.7, 0.8, 0.9, 1.0, 1.1, 1.2\}$ (blue curve) in Example 1. Shaded areas represents $\pm$ standard deviation over 8 independent runs.}
		\label{fig:distrib-matching-value-cmp}
	\end{figure}

	We implement the above model with: $v_{\textrm{true}} = (0.2, 0.7, 0.1)$, $v_{\textrm{ref}} = (0,1,0)$ and $\delta_{\textrm{perturb}}$ varying between $0.0$ and $0.8$. 
	\textit{Figure~\ref{fig:distrib-matching-value-cmp} shows that for moderately small $\delta_{\textrm{perturb}}$, the robust policy performs better than the non-robust policy.} For large values of $\delta_{\textrm{perturb}}$ however, the robust policy yields a smaller value: being robust against a large set of possible common noise distributions prevents the policy from performing well on the true distribution. The results are averaged over 8 different runs and the plots shows the mean value and its standard deviation. 
	
\begin{figure}[t]
		\centering
		\includegraphics[width=0.675\linewidth]{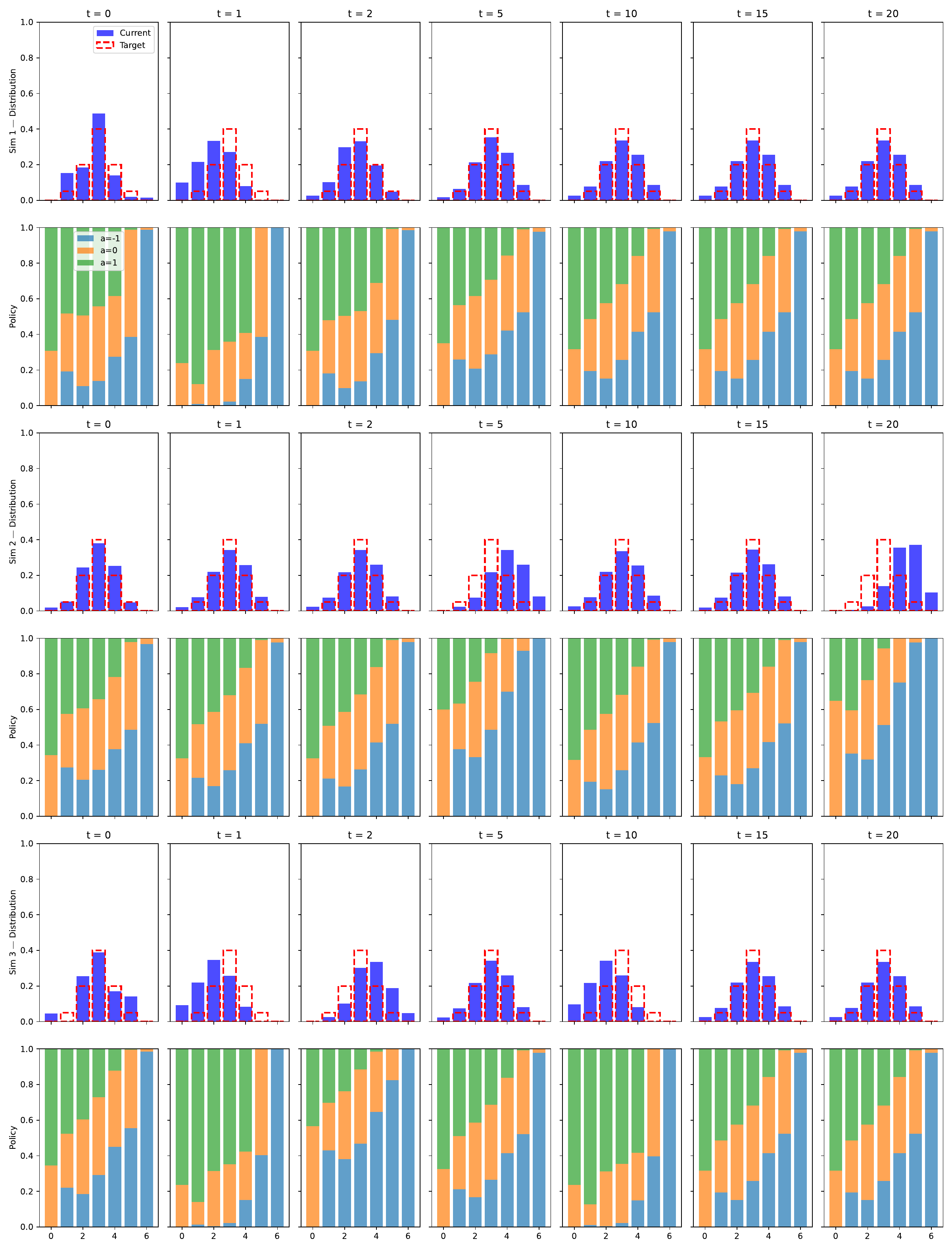}
		\caption{Three sample trajectories of the population distribution and
			corresponding action distribution for each state in Example 1. The target
			distribution to be matched is shown by dashed red lines.}
		\label{fig:distrib-matching-traj}
	\end{figure}

	Figure~\ref{fig:distrib-matching-traj} shows three realizations of trajectories, starting from random initial distributions. We display a few time steps between $0$ and $20$. We observe that the learnt policy uses the actions with varying proportions depending on the individual state and also depending on the current population distribution. Overall, it uses mostly action $1$ (resp.\;$-1$) when the state is below (resp.\;above) the middle state because the target distribution is centered around the middle state. 
	
	\begin{figure}
		\centering
		\includegraphics[width=0.65\linewidth]{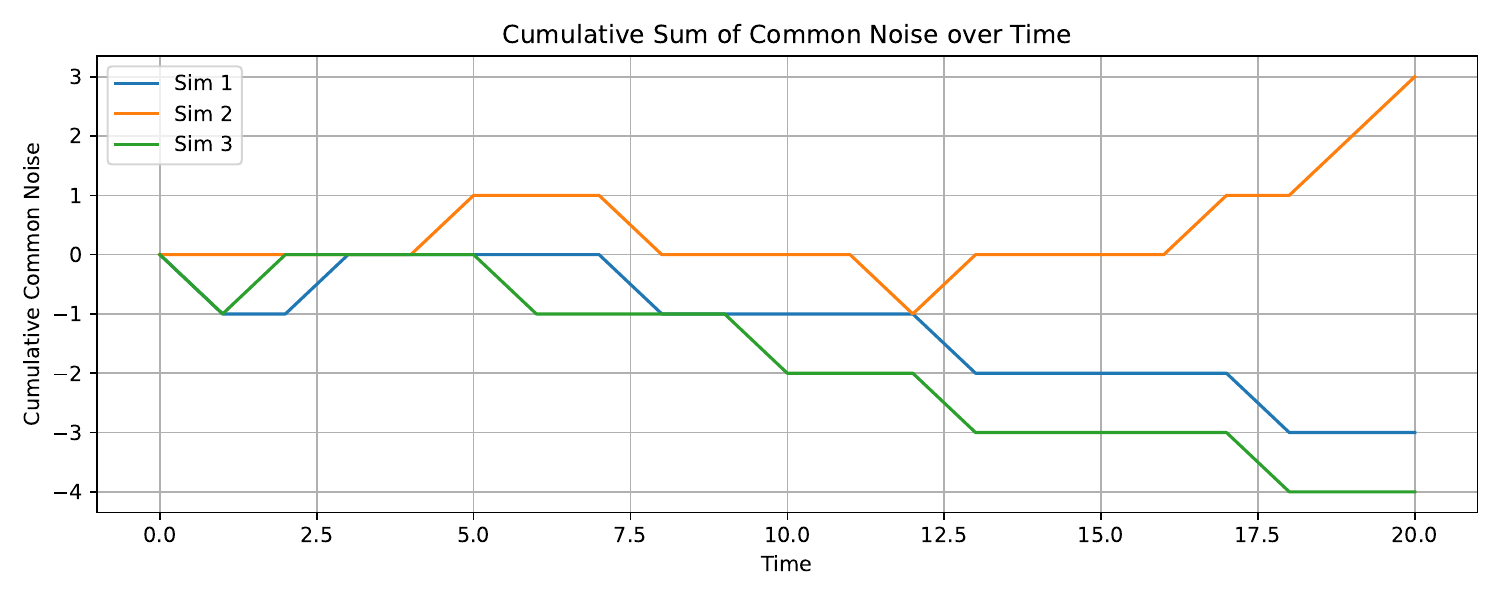}
		\caption{The three trajectories of common noise associated with Figure \ref{fig:distrib-matching-traj}}
		\label{fig:distrib-matching-traj-CN}
\end{figure}

	The fact that the target distribution is not perfectly matched is due to the impact of the common noise, whose trajectories are displayed in Figure~\ref{fig:distrib-matching-traj-CN}. Notice that for the second simulation, the common noise takes several positive values on time steps 17, 18 and 19, leaving little time for the population distribution to adapt and shift back to the target distribution (recall that the possible actions are $\{-1,0,1\}$, just as the possible common noise values).

	\subsection{Financial systemic risk}
	We now consider an example inspired by the systemic risk model proposed by~\cite{MR3325083}. In this model, the agents are financial institutions, represented by a state which is their log-reserve. They interact by borrowing and lending to each other, or to a central bank. Their evolution is impacted by a common noise which can be interpreted as macroscopic events affecting the whole economy. If a financial institution touches a given threshold, it defaults. 
	There are two main differences between the model we present below and the original model one: first, the model of~\cite{MR3325083} was a mean field game (corresponding to non-cooperative players) while we consider a mean field control problem (corresponding to cooperative players); furthermore, the original model was written in continuous space and time whereas we consider a discrete space and time model for the sake of numerical experiments. However, the main ideas underpinning the model are similar. The central planner is to make the population distribution match a given target distribution.
	\begin{figure}[t]
		\centering
		\includegraphics[width=0.6\linewidth]{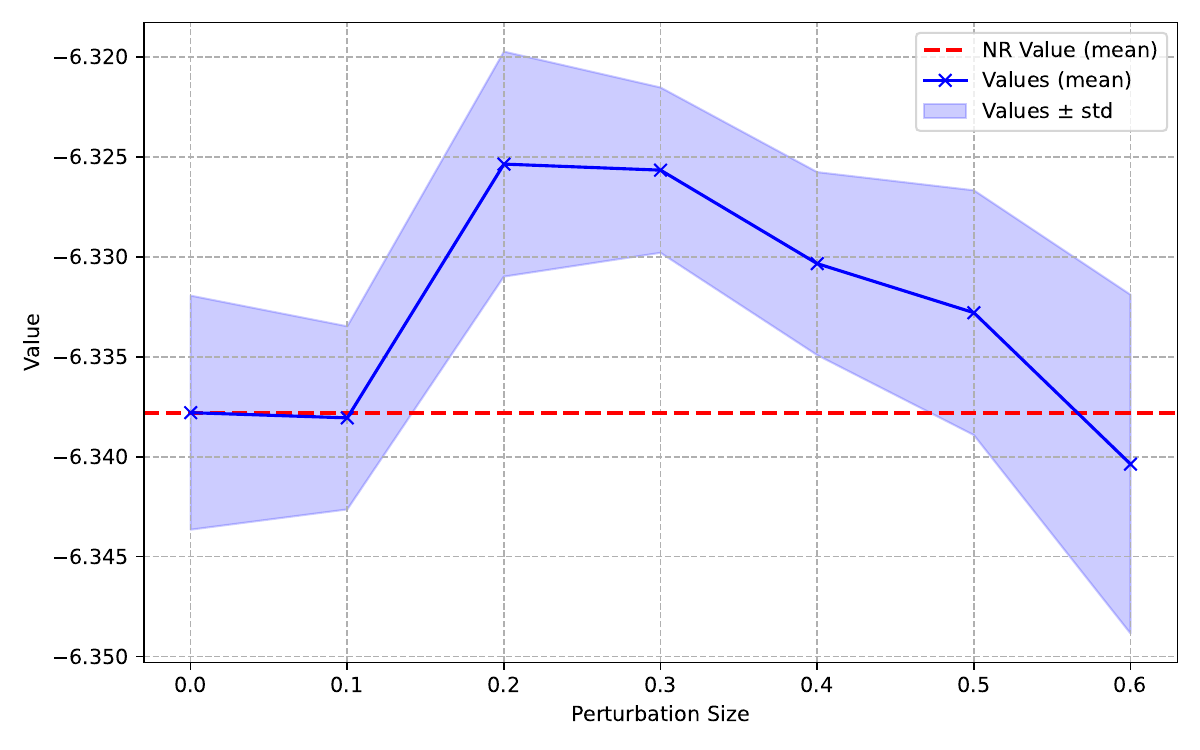}
		\caption{Value achieved under $p_{\textrm{true}}$ when using the optimal policy for the MFC with $p_{\textrm{ref}}$ (red dashed line) or the optimal policy for the robust MFC with $\delta_{\textrm{perturb}} \in \{0.0, 0.1, 0.2, 0.3, 0.4, 0.5, 0.6\}$ (blue curve) in Example 2. Shaded areas represents $\pm$ standard deviation over 8 independent runs.}
		\label{fig:systemic-risk-value-cmp}
	\end{figure} 
	
	To be specific, consider the following basic elements (recall Definition~\ref{dfn:basic_element}):
	\begin{itemize}[leftmargin=3.em]
		\item $S = \{s_{\operatorname{min}},s_{\operatorname{min}}+1,\dots, s_{\operatorname{max}}\}$, which represents a one-dimensional grid world with $|S| = s_{\operatorname{max}} - s_{\operatorname{min}} + 1$ states; in the experiments, we use $s_{\operatorname{min}} = -1$, $s_{\operatorname{max}} = 4$, $|S| =5$ states.
		\item $A = \{-1,0,1\}$, which corresponds to lending (if negative) or borrowing (if positive) units.
		\item $E = \{-1,0,1\}$, which corresponds to idiosyncratic noise. Moreover, the probaility vector of its law $\lambda_\varepsilon\in {\cal P}(E)$ is given by $(0.05,0.9,0.05)$.
		\item $E^0 = \{-2,-1,0,1,2\}$, which corresponds to common noise affecting the whole population.
		\item $\operatorname{F}:S \times A \times {\cal P}(S\times A) \times E \times E^0\to S$ 
		is given by 
		\[
		\operatorname{F}(s,a,\Lambda,e,e^0) = \max(s_{\operatorname{min}}, \min(s_{\operatorname{max}},s+a+e+e^0))\quad \mbox{if $s > s_{\operatorname{min}}$},
		\]
		and $\operatorname{F}(s_{\operatorname{min}},a,\Lambda,e,e^0) = s_{\operatorname{min}}$, which represents the fact that the agent's log-reserve evolution is determined by her action, the individual noise and the common noise, the agent remains at $1$ (resp.\;$7$) if she tries to move to the left (resp.\;right) of this state, and the agent remains stuck at $s=1$ if she ever reaches this state.
		\item $r:S\times A\times {\cal P}(S\times A)\to \mathbb{R}$ is given by  
		\[
			r(s,a,\Lambda) = -a^2 + q a (m({\Lambda}) - s)^2 - 0.5 \epsilon (m({\Lambda})- s)^2 + (m({\Lambda}) - s_{\mathrm{target}})^2,
		\]
		where $m({\Lambda})$ is given by $m({\Lambda}):=\int_{S}s'\operatorname{pj}_S(\Lambda)(ds')$ (i.e., the first moment of the state), the constants $q,\epsilon$ are non-positive and satisfy $q^2 \le \epsilon$, and $s_{\mathrm{target}}$ is a target state taken equal to $2$ in the experiments. The first term is a cost of borrowing\,/\,lending, the second and third terms have a mean-reverting effect, and the last term means that the regulator has a target level for the mean of the log-reserves. Here, $q$ represents the incentive to borrowing or lending. We refer to~\cite{MR3325083} for more details.
		\item {$\beta = 0.15$ is the discount factor so that Assumptions~\ref{as:general}\;(iii) and~\ref{as:general2}\;(iv) are satisfied.} 
	\end{itemize}

	\begin{figure}[t]
		\centering
		\includegraphics[width=0.7\linewidth]{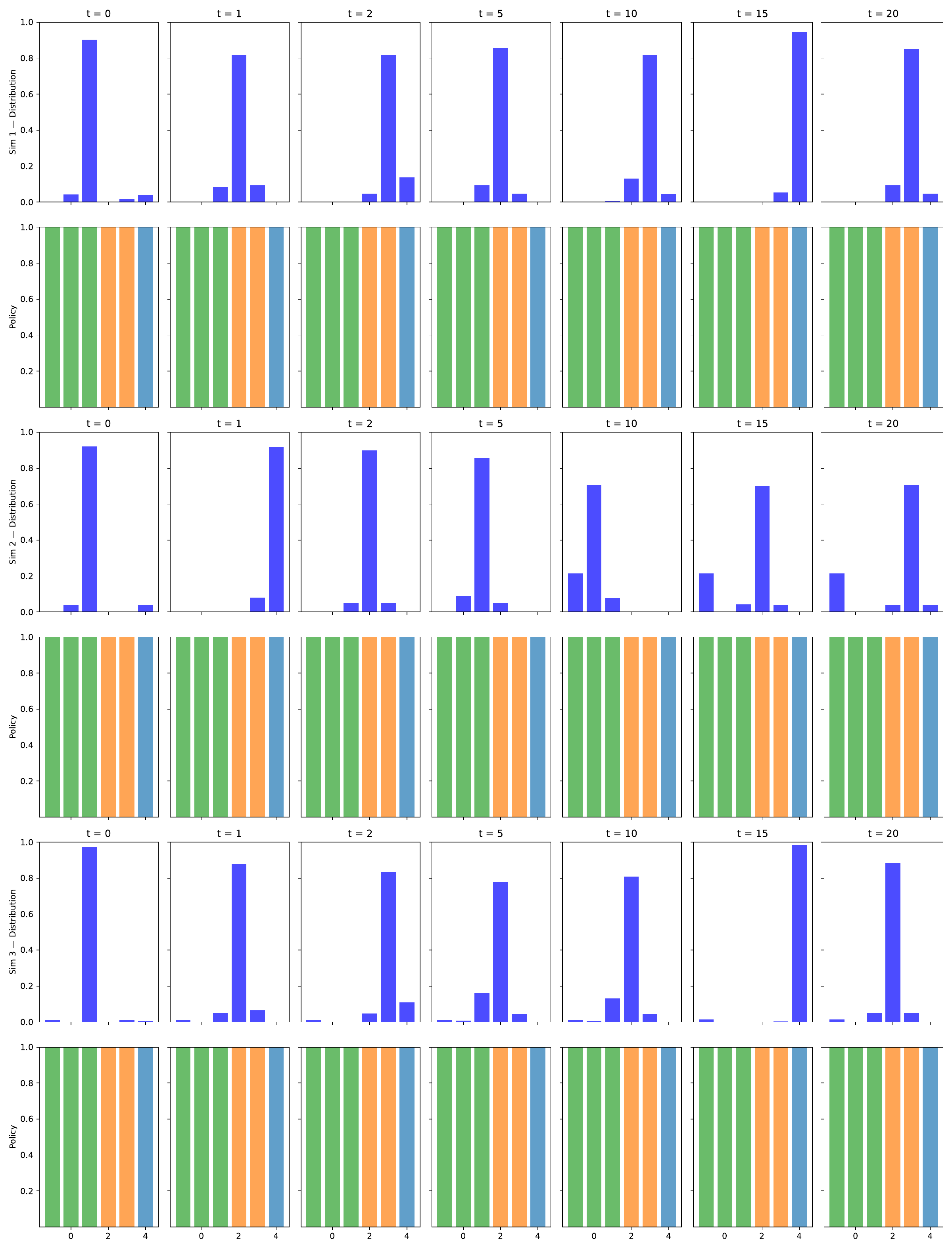}
		\caption{Three sample trajectories of the population distribution and corresponding action distribution for each state in Example 2. 
		}
		\label{fig:systemic-risk-traj}
	\end{figure}
	\begin{figure}
		\centering
		\includegraphics[width=0.6\linewidth]{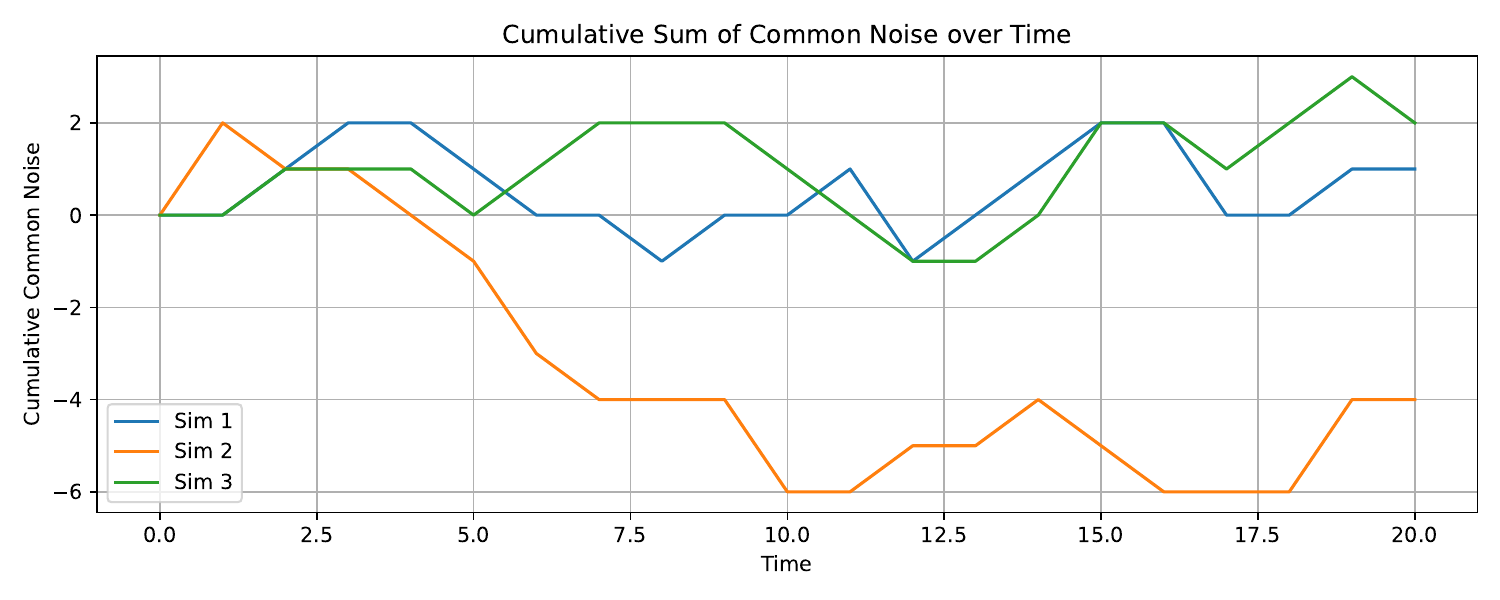}
		\caption{Three sample trajectories of common noise, associated to the three distribution trajectories presented in Figure~\ref{fig:systemic-risk-traj}.}
		\label{fig:systemic-risk-traj-CN}
	\end{figure}
	
	For the common noise probability measure, we proceed as in the previous example of Section~\ref{sec:ex1}. The true common noise measure is denoted by $p_{\textrm{true}}\in{\cal P}(E^0)$ (as in \eqref{eq:true}, but now represented by a 5-dimensional probability vetor $v_{\textrm{true}}\in[0,1]^5$). The central planner does not know this true measure and instead relies on a reference probability measure $p_{\textrm{ref}}\in{\cal P}(E^0)$ with corresponding probability vector $v_{\textrm{ref}}\in[0,1]^5$. We then compare, in the true model with $p_{\textrm{true}}$, the performance of $\pi_{\textrm{ref}}$ (an optimal policy for the model with common noise distribution $p_{\textrm{ref}}$) and the performance of $\pi_{\textrm{robust}}$ (a robust policy for $p_{\textrm{ref}}$). The uncertainty set $\mathfrak{P}^0$ is defined as in \eqref{eq:perturb}, but adapted to the 5-dimensional setting so that Assumption \ref{as:general2}\;(i) also holds.
	
	We implement the above model with: $v_{\textrm{true}} = (0.1, 0.2, 0.4, 0.2, 0.1)$, $v_{\textrm{ref}} = (0,0, 1, 0, 0)$ and $\delta_{\textrm{perturb}}$ varying between $0.0$ and $0.6$. 
	\textit{Figure~\ref{fig:systemic-risk-value-cmp} shows that for moderately small $\delta_{\textrm{perturb}}$, the robust policy performs better than the non-robust policy.} For large values of $\delta_{\textrm{perturb}}$ however, the robust policy yields a smaller value: being robust against a large set of possible common noise distributions prevents the policy from performing well on the true distribution. The results are averaged over 15 different runs and the plots shows the mean value and its standard deviation. Figure~\ref{fig:systemic-risk-traj} shows three realizations of trajectories, starting from random initial distributions. We display a few time steps between $0$ and $20$. We observe that the learnt policy is pure at the agent level, meaning that in each state, the agent uses one action with probability $1$. In fact, the agent uses actions that tend to make the state move towards state $2$ or $3$. The distribution concentrates (but not completely due to the idiosyncratic noise which tends to make the agent spread). Moreover, the peak is not always at state $2$ or $3$ due to the impact of the common noise, whose trajectories are displayed in Figure~\ref{fig:systemic-risk-traj-CN}.


	\section{Proof of results in Section \ref{subsec:propagation}}\label{proof:subsec:propagation}
	
	We begin by verifying the measurability of the state dynamics appearing in both models. 
	We recall the filtrations given in Definition \ref{dfn:filtrations}. 
	\begin{lem}\label{lem:measurability}
		For any $\pi \in \Pi$ and $\mathbb{P}\in {\cal Q}$, the following statements hold: 
		\begin{itemize}
			\item [(i)] For every $N\in \mathbb{N}$, $i=1,\dots,N$, and $t\geq 0$, $s_t^{i,N,\pi}$ given in \eqref{eq:MKV_N} is $\big(\bigvee_{j=1}^N{\cal F}_t^{j}\big)$ measurable. 
			\item [(ii)] For every $i\in \mathbb{N}$ and $t\geq 0$, $s_t^{i,\pi,\mathbb{P}}$ in \eqref{eq:MKV} is ${\cal F}^i_t$ measurable, and both $\mathbb{P}^0_{(s_t^{i,\pi,\mathbb{P}},a_t^{i,\pi})}$ and $\mathbb{P}_{s_t^{i,\pi,\mathbb{P}}}^0$ are ${\cal F}_t^0$ measurable.
		\end{itemize}
	\end{lem}
	\begin{proof} We start proving (i). Let $N\in \mathbb{N}$ and $i=1,\dots,N$ be given. The statement is shown via an induction over $t\geq 0$: Since $s_0^{i,N,\pi}=\xi^i\in L^0_{{\cal F}_0^i}(S)$ (see Definition \ref{dfn:N_model}), the claim for $t=0$ holds.  
	
	Now assume that the induction claim holds for some $t\geq 0$. Note that $s_{t+1}^{i,N,\pi}$ satisfies 
	\[
	s_{t+1}^{i,N,\pi}=\operatorname{F}(s^{i,N,\pi}_t,a^{i,\pi}_t,\mbox{$\frac{1}{N}\sum_{j=1}^N\delta_{(s^{j,N,\pi}_t,a^{j,\pi}_t)}$},\varepsilon_{t+1}^{i},\varepsilon_{t+1}^0)
	\]
	where the first three terms are $(\bigvee_{j=1}^N{\cal G}_t^j)$ measurable because of the induction assumption and 
	the definition of the open-loop control $\alpha_t^{i,\pi}$ in Definition \ref{dfn:N_model}\;(i), and the fact that $\bigvee_{j=1}^N{\cal F}_t^j\subset\bigvee_{j=1}^N{\cal G}_t^j$. Hence by the Borel measurability of $\operatorname{F}$, $s_{t+1}^{i,N,\pi}$ is $(\bigvee_{j=1}^N{\cal F}_{t+1}^j)$ measurable (see Definition \ref{dfn:filtrations}). 
	
	By the induction hypothesis, the statement in (i) holds for all $t\geq 0$.
	
	\vspace{0.5em}
	\noindent The part (ii) is also shown via an induction over $t$ given any $i\in \mathbb{N}$. Since $s_0^{i,\pi,\mathbb{P}}=\xi^i\in L_{{\cal F}_0^i}^0(S)$ (see Definition~\ref{dfn:MFC}), $s_0^{i,\pi,\mathbb{P}}$ is ${\cal F}_0^i$ measurable. Moreover, since ${\cal F}_0^0$ is trivial, both $\mathbb{P}^0_{(s_0^{i,\pi,\mathbb{P}},a_0^{i,\pi})}$ and $\mathbb{P}_{s_0^{i,\pi,\mathbb{P}}}^0$ are ${\cal F}_t^0$ measurable obviously.
	
	We assume that the claim holds for some $t\geq 0$. Note that $s_{t+1}^{i,\pi,\mathbb{P}}$ satisfies 
	\[
	s_{t+1}^{i,\pi,\mathbb{P}}=\operatorname{F}(s^{i,\pi,\mathbb{P}}_t,a^{i,\pi}_t,\mathbb{P}^{0}_{(s^{i,\pi,\mathbb{P}}_t,a^{i,\pi}_t)},\varepsilon_{t+1}^{i},\varepsilon_{t+1}^0),
	\]
	where the first three terms are ${\cal G}_t^i$ measurable because of the induction assumption and the fact that ${\cal F}_t^0\subset {\cal G}_t^i$. Hence by the Borel measurability of $\operatorname{F}$, $s_{t+1}^{i,\pi,\mathbb{P}}$ is ${\cal F}_{t+1}^i$ measurable (see Definition~\ref{dfn:filtrations}). 
	
	Moreover, since $a^{i,\pi}_{t+1}$ is ${\cal G}_{t+1}^i$ measurable and  
	$(\gamma^i,\vartheta_{0:t+1}^i,\varepsilon^i_{1:t+1})$ is independent of $\varepsilon^0_{1:t+1}$ (see Remark~\ref{rem:well_dfn_1}\;(i)), we apply Lemma \ref{lem:regular_cond}\;(ii) to have that both $\mathbb{P}^0_{(s_{t+1}^{i,\pi,\mathbb{P}},a_{t+1}^{i,\pi})}$ and $\mathbb{P}_{s_{t+1}^{i,\pi,\mathbb{P}}}^0$ are ${\cal F}_{t+1}^0$ measurable. By the induction hypothesis, the statement in (ii) holds. 
	\end{proof}

		\subsection{Proof of Lemma \ref{lem:propagation_chaos}}
		 We start proving (i). Let $q>2$ be given. Note that by Lemma~\ref{lem:measurability}\;(ii), the definition of open-loop controls (see Definition \ref{dfn:N_model}\;(i)), and recalling that $\mathbb{F}^i\subset \mathbb{G}^i$ for any~$i\in \mathbb{N}$ $(s_{t}^{i,\pi,\mathbb{P}},a_{t}^{i,\pi})$ is ${\cal G}_t^i$ measurable.  
		 
		 Moreover, since the private components $(\gamma^i)_{i\in \mathbb{N}}$, $(\vartheta_t^i)_{t\geq 0,i\in \mathbb{N}}$, and $(\varepsilon_t^i)_{t\geq 1, i\in \mathbb{N}}$ are mutually independent (see Remark~\ref{rem:well_dfn_1}\;(i)) and all agents are indistinguishable, it holds for every $t\geq 0$, $\pi \in \Pi$, and $\mathbb{P}\in {\cal Q}$ that $(s_{t}^{i,\pi,\mathbb{P}},a_{t}^{i,\pi})_{i\in \mathbb{N}}$ is (conditionally) i.i.d.\;given the common noise information~${\cal F}^0_t$ with law $\mathbb{P}^0_{(s_{t}^{1,\pi,\mathbb{P}},a_{t}^{1,\pi})}$.  Therefore, it follows from \cite[Theorem 1]{fournier2015rate} that 
			\[
			\mathbb{E}^{\mathbb{P}^0} \Big[{\cal W}_{{\cal P}(S\times A)}\big(\mbox{$\frac{1}{N}\sum_{i=1}^N\delta_{(s^{i,\pi,\mathbb{P}}_t,a^{i,\pi}_t)}$},\, \mathbb{P}^0_{(s^{1,\pi,\mathbb{P}}_t,a^{1,\pi}_t)}\big) \Big]\leq C \big(K_q(\mathbb{P}^0_{(s^{1,\pi,\mathbb{P}}_t,a^{1,\pi}_t)})\big)^{1/q}  \alpha(N),
			\]
			where $C>0$ does not depends on $\mathbb{P}^0$ and $N$ but on $d$ and $q$, $\alpha(\cdot)$ is defined as in the statment, and $K_q(\mathbb{P}^0_{(s^{1,\pi,\mathbb{P}}_t,a^{1,\pi}_t)})$ is given by 
			\[
			K_q(\mathbb{P}^0_{(s^{1,\pi,\mathbb{P}}_t,a^{1,\pi}_t)}):=\int_{S\times A} |(s,a)|^q\;\mathbb{P}^0_{(s^{1,\pi,\mathbb{P}}_t,a^{1,\pi}_t)}(ds,da).
			\]
			Since $S\times A$ is a compact subset of $\mathbb{R}^d$, the above quantitiy is uniformly bounded by $(\Delta_{S\times A})^{q}$ for every $t\geq 0$, $\pi \in \Pi$, and $\mathbb{P}\in {\cal Q}$. Hence the estimate in part (i)~holds.
			
			\vspace{0.5em}
			\noindent Last, we prove (ii). Let $q>2$ be given. In part (i), we have verified that for every $t\geq 0$, $\pi \in \Pi$, and $\mathbb{P}\in {\cal Q}$, $(s_{t}^{i,\pi,\mathbb{P}},a_{t}^{i,\pi})_{i\in \mathbb{N}}$ is (conditionally) i.i.d.\;given ${\cal F}^0_t$ with law $\mathbb{P}^0_{(s_{t}^{1,\pi,\mathbb{P}},a_{t}^{1,\pi})}$.

			Hence, we can apply \cite[Corollary~1.2]{boissard2014mean} to obtain that for every $t\geq 0$, $\pi \in \Pi$, and $\mathbb{P}\in {\cal Q}$
			\[
			\mathbb{E}^{\mathbb{P}^0} \Big[{\cal W}_{{\cal P}(S\times A)}\big(\mbox{$\frac{1}{N}\sum_{i=1}^N\delta_{(s^{i,\pi,\mathbb{P}}_t,a^{i,\pi}_t)}$},\, \mathbb{P}^0_{(s^{1,\pi,\mathbb{P}}_t,a^{1,\pi}_t)}\big) \Big]\leq c \Big(\frac{2}{q-2}\Big)^{\frac{2}{q}}  (k_{S\times A})^{\frac{1}{q}} \Delta_{S\times A}  N^{-\frac{1}{q}},
			\]
			with some $c\leq 64/3$. Therefore, we can obtain the estimate in part (ii), as claimed.  \qed

			\subsection{Proof of Theorem \ref{thm:propagation_chaos}}  For notational simplicity, throuhgout this proof, denote for every $N\in \mathbb{N}$, $i=1,\dots,N$, $t\geq 0$, $\pi \in \Pi$, and $\mathbb{P}\in {\cal Q}$ by 
				\begin{align*}
					\begin{aligned}
						&\Lambda_t^{N,\pi}:=\mbox{$\frac{1}{N}\sum_{j=1}^N\delta_{(s^{j,N,\pi}_t,a^{j,\pi}_t)}$},\qquad \Lambda_t^{N,\infty,\pi,\mathbb{P}}:=\mbox{$\frac{1}{N}\sum_{j=1}^N\delta_{(s^{j,\pi,\mathbb{P}}_t,a^{j,\pi}_t)}$},\\
						&\tilde \Lambda_t^{i,\pi,\mathbb{P}}:=\mathbb{P}^0_{(s^{i,\pi,\mathbb{P}}_t,a^{i,\pi}_t)}.
					\end{aligned}
				\end{align*}
				
				Let $N\in \mathbb{N}$ and $i=1,\dots,N$ be given. We first prove \eqref{eq:cvg_state} and \eqref{eq:cvg_emprical}. The proof uses an induction over $t\geq 0$: Since $s^{i,N,\pi}_0=s_0^{i,\pi,\mathbb{P}}$ for every $\pi\in \Pi$,  and $\mathbb{P}\in {\cal Q}$ (see Definitions \ref{dfn:N_model} and~\ref{dfn:MFC}), the claim for $t=0$~holds.  
				
				Now assume that the induction claim holds true for some $t\geq 1$. Let $\pi \in \Pi$ and $\mathbb{P}\in {\cal Q}$ be given. 
				Since $\bigvee_{j=1}^N{\cal F}_t^j\subset \bigvee_{j=1}^N{\cal G}_t^j$ (see Definition \ref{dfn:filtrations}), both $s_t^{i,N,\pi}$ given in \eqref{eq:MKV_N} and $s_t^{i,\pi,\mathbb{P}}$ given in \eqref{eq:MKV} are $(\bigvee_{j=1}^N{\cal G}_t^j)$ measurable 
				(see Lemma~\ref{lem:measurability}).
				Moreover, $a_t^{i,\pi}$ is ${\cal G}_t^i$ measurable (see Definition~\ref{dfn:N_model}~(i)).

				Since $\varepsilon_{t+1}^i$ is independent of $\bigvee_{j=1}^N{\cal G}_{t}^j$ and $\varepsilon_{t+1}^0$ (see Remark \ref{rem:well_dfn_1}\;(i),\;(ii)), we can have the following conditioning
				\begin{align}\label{eq:state_induct_prpc1}
					\mathbb{E}^{\mathbb{P}}[d_S(s^{i,N,\pi}_{t+1},s^{i,\pi,\mathbb{P}}_{t+1})]=\mathbb{E}^{\mathbb{P}}[\operatorname{D}^{i,\mathbb{P}}(s^{i,N,\pi}_{t},s^{i,\pi,\mathbb{P}}_{t}, a_{t}^{i,\pi,\mathbb{P}}, \Lambda_t^{N,\pi},\tilde \Lambda_t^{i,\pi,\mathbb{P}}, e^0)],
				\end{align}
				where for every $(s,\tilde s)\in S$, $a\in A$, $\Lambda,\tilde \Lambda \in {\cal P}(S\times A)$, and $e^0\in E^0$
				\begin{align}\label{eq:state_induct_prpc2}
					\begin{aligned}
						\operatorname{D}^{i,\mathbb{P}}(s,\tilde s, a, \Lambda, \tilde \Lambda, e^0):=& 
						\int_E d_S\big(\operatorname{F}(s,a,\Lambda,e,e^0),\operatorname{F}(\tilde s,a,\tilde \Lambda,e,e^0)\big)\lambda_\varepsilon(de)\\
						\leq& C_{\operatorname{F}}\big(d_S(s,\tilde s)+\mathcal{W}_{{\cal P}(S\times A)}(\Lambda, \tilde \Lambda ) \big),
					\end{aligned}
				\end{align}
				where the inequality follows from  Assumption \ref{as:general}\;(i).
				
				On the other hand, it holds that
				\begin{align}
					\mathbb{E}^{\mathbb{P}}[\mathcal{W}_{{\cal P}(S\times A)}(\Lambda_t^{N,\pi},\tilde \Lambda_t^{i,\pi,\mathbb{P}})]
					&\leq \mathbb{E}^{\mathbb{P}}[\mathcal{W}_{{\cal P}(S\times A)}(\Lambda_t^{N,\pi},\Lambda_t^{N,\infty,\pi,\mathbb{P}})]+\mathbb{E}^{\mathbb{P}}[\mathcal{W}_{{\cal P}(S\times A)}(\Lambda_t^{N,\infty,\pi,\mathbb{P}},\tilde \Lambda_t^{i,\pi,\mathbb{P}})] \nonumber\\
					&\leq \mathbb{E}^{\mathbb{P}}[d_S(s^{i,N,\pi}_{t},s^{i,\pi,\mathbb{P}}_{t})] +  M_N,\label{eq:state_induct_prpc3}
				\end{align}
				where the second inequality follows from the definition of $M_N$ given in \eqref{eq:propag_MN} and the fact that $\mathcal{W}_{{\cal P}(S\times A)}(\Lambda_t^{N,\pi},\Lambda_t^{N,\infty,\pi,\mathbb{P}})\leq \frac{1}{N}\sum_{j=1}^Nd_S(s^{j,N,\pi}_{t},s^{j,\pi,\mathbb{P}}_{t})$ together with the indistinguishability.
				
				Combining \eqref{eq:state_induct_prpc1} with \eqref{eq:state_induct_prpc2} and \eqref{eq:state_induct_prpc3}, we have that
				\begin{align}\label{eq:state_induct_prpc4}
					\mathbb{E}^{\mathbb{P}}[d_S(s^{i,N,\pi}_{t+1},s^{i,\pi,\mathbb{P}}_{t+1})]\leq C_{\operatorname{F}} \big(2\mathbb{E}^{\mathbb{P}}[d_S(s^{i,N,\pi}_{t},s^{i,\pi,\mathbb{P}}_{t})] +  M_N \big).
				\end{align}
				Since the estimate \eqref{eq:state_induct_prpc4} holds for any $\pi \in \Pi$ and $\mathbb{P}\in {\cal Q}$, by the induction hypothesis we have that the estimate \eqref{eq:cvg_state} holds for all $t\geq 0$, as claimed.
				
				Moreover, since the estimate \eqref{eq:state_induct_prpc3} holds for any $\pi \in \Pi$ and $\mathbb{P}\in {\cal Q}$, by using \eqref{eq:cvg_state} we have that the other estimate \eqref{eq:cvg_emprical} holds for all $t\geq 0$, as claimed. As $N\in \mathbb{N}$ and $i=1,\dots,N$ are given arbitrary, we can conclude that \eqref{eq:cvg_state} and \eqref{eq:cvg_emprical} hold for all $N\in \mathbb{N}$, $i=1,\dots,N$, and $t\geq 0$.
				
				\vspace{0.5em}
				\noindent We now prove \eqref{eq:cvg_reward}. Note that for every $N\in \mathbb{N}$ and $\pi \in \Pi$
				\begin{align}\label{eq:reward_1}
					\begin{aligned}
						&\big|{\cal J}^{N,\pi}-{\cal J}^{\pi}\big|= \bigg|\inf_{\mathbb{P}\in {\cal Q}}\mathbb{E}^{\mathbb{P}}\bigg[ \frac{1}{N}\sum_{i=1}^NR^{i,N,\pi} \bigg]- \inf_{\mathbb{P}\in {\cal Q}}\mathbb{E}^{\mathbb{P}}\bigg[ \frac{1}{N}\sum_{i=1}^NR^{i,\pi,\mathbb{P}} \bigg]\bigg|\\ 
						&\qquad \leq \sup_{\mathbb{P}\in {\cal Q}} \frac{1}{N}\sum_{i=1}^N \mathbb{E}^{\mathbb{P}}\Big[ |R^{i,N,\pi}- R^{i,\pi,\mathbb{P}}| \Big]= \sup_{\mathbb{P}\in {\cal Q}} \mathbb{E}^{\mathbb{P}}\Big[ |R^{1,N,\pi}- R^{1,\pi,\mathbb{P}}| \Big] \\
						&\qquad \leq \sum_{t=0}^\infty \beta^t  \sup_{\mathbb{P}\in {\cal Q}} \mathbb{E}^{\mathbb{P}}\Big[ \big|r(s_{t}^{1,N,\pi},a_{t}^{1,\pi,},\Lambda_t^{N,\pi})- r(s_{t}^{1,\pi,\mathbb{P}},a_{t}^{1,\pi},\tilde \Lambda_t^{1,\pi,\mathbb{P}})\big| \Big]=:\operatorname{I}^{N,\pi},
					\end{aligned}
				\end{align}
				where the equalities follow from the indistinguishability and the last inequality holds because $r$ is bounded (see Assumption \ref{as:general}\;(ii)).
				
				Moreover, by the Lipschitz continuity of $r(\cdot,a,\cdot):S\times {\cal P}(S\times A)\to \mathbb{R}$ for any $a\in A$ (see Assumption \ref{as:general}\;(ii)), for every $N\in \mathbb{N}$ and $\pi \in \Pi$
				\begin{align}\label{eq:reward_2}
					\begin{aligned}
						\operatorname{I}^{N,\pi}&\leq C_{r} \sum_{t=0}^\infty \beta^t\sup_{\mathbb{P}\in {\cal Q}} \mathbb{E}^{\mathbb{P}}\Big[d_S(s^{1,N,\pi}_{t},s^{1,\pi,\mathbb{P}}_{t})+{\cal W}_{{\cal P}(S\times A)}(\Lambda_t^{N,\pi},\tilde \Lambda_t^{1,\pi,\mathbb{P}} )  \Big]\\
						&\leq  C_{r} \bigg( 2 \sum_{t=0}^\infty \beta^t \delta_t^N+ \frac{M_N}{1-\beta} \bigg),
					\end{aligned}
				\end{align}
				where $\delta_t^N:= \sup_{\pi \in \Pi}\sup_{\mathbb{P}\in {\cal Q}}\mathbb{E}^{\mathbb{P}}[d_S\big(s^{1,N,\pi}_t,s^{1,\pi,\mathbb{P}}_t\big)]$ for $t\geq 0$. 
				
				Since the estimate given in \eqref{eq:reward_1} coincides with that of \cite[Theorem~2.1]{motte2022mean}---specifically Eq.\;(2.17) therein---and Assumption \ref{as:general}\;(iii) ensures that $2\beta C_r <1 $, we can follow the same calculations as in the proof of \cite[Theorem~2.1]{motte2022mean} (replacing $K_F$ with $C_r$). This yields that $\sum_{t=0}^\infty \beta^t \delta_t^N \leq C M_N$ for some constant $C>0$ (that do not depend on $N$ and $\pi$); see also \cite[Remark~2.4]{motte2022mean}. 
				
				Combining this with \eqref{eq:reward_1} and \eqref{eq:reward_2} establishes the estimate in \eqref{eq:cvg_reward}. 
				\qed

				\section{Proof of results in Section \ref{subsec:lift_MDP}}\label{proof:subsec:lift_MDP}
				\subsection{Proof of Proposition \ref{pro:lift_dynamics}} We first prove \eqref{eq:lift_1}. For simplicity, denote for every $t\geq 0$ by
					\begin{align}\label{eq:abbrev}
						\mu_t:=\mu_t^{\xi,a,\mathbb{P}},\qquad \Lambda_t:=\Lambda_t^{\xi,a,\mathbb{P}},\qquad \nu_{t+1}:=\mathscr{L}_{\mathbb{P}}(\varepsilon_{t+1}^0|{\cal F}_t^0).
					\end{align}
					
					Since $\mu_{t+1} $ is ${\cal F}_{t+1}^0$ measurable, it is sufficient to show that for any bounded Borel measurable functions $\hat g:(E^0)^{t+1}\to \mathbb{R}$ and $\hat f:S\to \mathbb{R}$,
					\begin{align}\label{eq:claim2}
						\mathbb{E}^{\mathbb{P}}[\hat g(\varepsilon_{1:t+1}^0) \hat f(s_{t+1}^{\xi,a,\mathbb{P}}) ] = \mathbb{E}^{\mathbb{P}}\bigg[\hat g(\varepsilon^0_{1:t+1})\int_S\hat f(s')\overline{\mathrm{F}}(\operatorname{pj}_S(\Lambda_t),\Lambda_t,\varepsilon_{t+1}^0)(ds')\bigg],
					\end{align}
					where we note that $(\operatorname{pj}_S(\Lambda_t),\Lambda_t)\in \operatorname{gr}(\mathfrak{U})$ (see Definition \ref{dfn:lift_maps}\;(i)).
					
					Note that by Remark \ref{rem:well_dfn_1}\;(i) and (ii), $\varepsilon_{t+1}$ is independent of $\varepsilon_{1:t+1}^0,s_t,a_t,\mathbb{P}_{(s_t,a_t)}^0$ (since they are all ${\cal G}_t\vee \sigma(\varepsilon_{t+1}^0)$ measurable) with $\mathscr{L}_{\mathbb{P}}(\varepsilon_{t+1})=\lambda_\varepsilon$. Moreover, by \eqref{eq:MKV_represent} and Fubini's theorem (noting that $\hat g$ and $\hat f$ are both bounded)
					\begin{align*}
						\mathbb{E}^{\mathbb{P}}[\hat g(\varepsilon_{1:t+1}^0) \hat f(s_{t+1}^{\xi,a,\mathbb{P}}) ] 
						&=\mathbb{E}^{\mathbb{P}}\bigg[\mathbb{E}^{\mathbb{P}}\Big[\hat g(\varepsilon_{1:t+1}^0) \hat f(\operatorname{F}(s^{\xi,a,\mathbb{P}}_t,a_t,\Lambda_t,\varepsilon_{t+1},\varepsilon_{t+1}^0)) \Big|\,e=\varepsilon_{t+1} \Big] \bigg]\\
						&=\int_{E} \mathbb{E}^{\mathbb{P}}\bigg[\hat g(\varepsilon_{1:t+1}^0) \hat f(\operatorname{F}(s^{\xi,a,\mathbb{P}}_t,a_t,\Lambda_t,e,\varepsilon_{t+1}^0)) \bigg]  \lambda_{\varepsilon}(de)=:\operatorname{I}.
					\end{align*}
					
					Note that $\varepsilon_{1:t}^0,s^{\xi,a,\mathbb{P}}_t,a_t,$ and $\Lambda_t$ are all ${\cal G}_t$ measurable. Since $\varepsilon_{t+1}^0$ is conditionally independent of ${\cal G}_t$ given ${\cal F}_t^0$ (see Remark~\ref{rem:well_dfn_1}\;(iii)), by definition of  $\nu_{t+1}$ (see \eqref{eq:abbrev}) 
					\begin{align*}
						\begin{aligned}
							\operatorname{I}&=\int_{E} \mathbb{E}^{\mathbb{P}}\bigg[ \mathbb{E}^{\mathbb{P}}\bigg[ \int_{E^0}  \hat g(\varepsilon_{1:t}^0,e^0) \hat f\big(\operatorname{F}(s^{\xi,a,\mathbb{P}}_t,a_t,\Lambda_t,e,e^0)\big)\nu_{t+1}(de^0) \bigg| {\cal F}_t^0 \bigg]\bigg]  \lambda_{\varepsilon}(de)\\
							&= \int_{E} \mathbb{E}^{\mathbb{P}}\bigg[ \int_{E^0}  \hat g(\varepsilon_{1:t}^0,e^0) \mathbb{E}^{\mathbb{P}}\Big[\hat f(\operatorname{F}(s^{\xi,a,\mathbb{P}}_t,a_t,\Lambda_t,e,e^0))\Big| {\cal F}_{t}^0\Big]  \nu_{t+1}(de^0) \bigg] \lambda_\varepsilon(de)
							=:\operatorname{II}.
						\end{aligned}
					\end{align*}
					
					Moreover by definition of $\Lambda_t$ (see \eqref{eq:abbrev}) and Fubini's theorem
					\begin{align*}
						\operatorname{II}
						&= \int_{E} \mathbb{E}^{\mathbb{P}}\bigg[ \int_{E^0}\hat g(\varepsilon_{1:t}^0,e^0) \int_{S\times A}\hat f({\operatorname{F}}(s,a,\Lambda_t,e,e^0))\Lambda_t(ds,da) \nu_{t+1}(de) \bigg]\lambda_\varepsilon(de)\\
						&=\mathbb{E}^{\mathbb{P}}\bigg[ \hat g(\varepsilon_{1:t+1}^0) \int_{S\times A\times E}\hat f({\operatorname{F}}(s,a,\Lambda_t,e,\varepsilon_{t+1}^0))\Lambda_t(ds,da) \lambda_{\varepsilon}(de) \bigg].
					\end{align*}
				
					By definition of $\overline{\mathrm{F}}$ (see Definition \ref{dfn:lift_maps}\;(ii)), the last term above is equal to the second term given in \eqref{eq:claim2}, as claimed.
					
					\vspace{0.5em}
					\noindent We now prove \eqref{eq:lift_2}. Note that by Remark \ref{rem:well_dfn_1}\;(iii) $(\nu_t)_{t\geq 0}$ given in \eqref{eq:abbrev} satisfies $\mathbb{P}$-a.s.
					\[
					\nu_1= p_1, \quad \nu_t=p_t(\cdot|\varepsilon_{1:t-1}^0)\quad \mbox{for all $t\geq 2$},
					\]
					where $(p_t)_{t\geq 1}\in \mathcal{K}^0$ induces the measure $\mathbb{P}\in {\cal Q}$.
					
					Let $t\geq 1$. It is sufficient to show that for any bounded Borel~measurable function $\tilde f: {\cal P}(S)\to \mathbb{R}$
					\begin{align}\label{eq:claim3}
						\mathbb{E}^{\mathbb{P}}[\tilde f (\mu_{t+1}) ] = \mathbb{E}^{\mathbb{P}}\bigg[\int_{{\cal P}(S)} \tilde f (\mu' ) \overline p\big(d\mu'|\operatorname{pj}_S(\Lambda_t),\Lambda_t,\nu_{t+1}\big) \bigg].
					\end{align}
					
					By \eqref{eq:lift_1}, we have $\mu_{t+1}=\overline{\mathrm{F}}(\operatorname{pj}_S(\Lambda_t),\Lambda_t,\varepsilon_{t+1}^0) $ $\mathbb{P}$-a.s.. Moreover, since $\varepsilon_{t+1}^0$ is conditionally independent of $(\operatorname{pj}_S(\Lambda_t),\Lambda_t)$ given ${\cal F}_t^0$ (as $\Lambda_t$ is ${\cal G}_t$ measurable) with $\mathscr{L}_{\mathbb{P}}(\varepsilon^0_{t+1}|{\cal F}_{t}^0)=\nu_{t+1}$, it follows that
					\begin{align*}
						\mathbb{E}^{\mathbb{P}}[\tilde f (\mu_{t+1}) ] =\mathbb{E}^{\mathbb{P}}\Big[\mathbb{E}^{\mathbb{P}}\big[f\big(\overline{\mathrm{F}}\big(\operatorname{pj}_S(\Lambda_t),\Lambda_t,\varepsilon_{t+1}^0\big) \big)\big|{\cal F}_t^0\big] \Big]= \mathbb{E}^{\mathbb{P}}\bigg[\int_{E^0}f\big(\overline{\mathrm{F}}\big(\operatorname{pj}_S(\Lambda_t),\Lambda_t,e^0\big) \big) \nu_{t+1}(de^0)\bigg].
					\end{align*}
					By definition of $\overline p$ (see Definition \ref{dfn:lift_maps}\;(iii)),
					the claim~\eqref{eq:claim3} holds. 
					
					For the case $t=0$, note that $\mathscr{L}_{\mathbb{P}}(\varepsilon^0_{1})=p_1$ and $\Lambda_0\in {\cal P}(S\times A)$ is deterministic. Thus, it is straightforward to verify that \eqref{eq:lift_2} holds also for $t=0$.
					
					This completes the proof. \qed

			\subsection{Proof of Proposition \ref{pro:dpp}} 
			In what follows, we often make use of the following coupling result along with the continuity of the projection map $\operatorname{pj}_S:{\cal P}(S\times A)\to {\cal P}(S)$.
			\begin{lem}\label{lem:element_1Wass}
				The following properties hold:
				\begin{enumerate}
					\item [(i)] For every $(\mu,\zeta),(\tilde \mu, \tilde \zeta) \in {\cal P}( S)\times  {\cal P}( A)$ and every $\Lambda \in \operatorname{Cpl}_{S\times A}(\mu,\zeta)$, there exists a coupling $\tilde \Lambda^* \in \operatorname{Cpl}_{S\times A}(\tilde \mu,\tilde \zeta)$ such~that 
					\[
					{\cal W}_{{\cal P}({S\times A})}(\Lambda,\tilde \Lambda^*)\leq {\cal W}_{{\cal P}(S)}(\mu,\tilde \mu)+{\cal W}_{{\cal P}(A)}(\zeta,\tilde \zeta).
					\]
					\item [(ii)] For every $\Lambda,\tilde \Lambda\in {\cal P}(S\times A) $, it holds that 
					\[
					{\cal W}_{{\cal P}(S)}(\operatorname{pj}_S(\Lambda),\operatorname{pj}_S(\tilde \Lambda))\leq {\cal W}_{{\cal P}(S\times A)}(\Lambda,\tilde \Lambda).
					\]
					Thus $\operatorname{pj}_S:{\cal P}(S\times A)\to {\cal P}(S)$ is continuous.
				\end{enumerate}
			\end{lem}
			\begin{proof} We start by proving (i). 
				Let $(\mu,\zeta),(\tilde \mu, \tilde \zeta) \in {\cal P}( S)\times  {\cal P}( A)$ and $\Lambda \in \operatorname{Cpl}_{S\times A}(\mu,\zeta)$.~Denote by 
				\begin{align}\label{eq:const1}
					\Gamma\in \operatorname{Cpl}_{S\times S}(\mu,\tilde \mu),\qquad \Upsilon \in \operatorname{Cpl}_{A\times A}(\zeta,\tilde \zeta) 
				\end{align}
				the optimal couplings for ${\cal W}_{{\cal P}({S})}(\mu,\tilde \mu)$ and ${\cal W}_{{\cal P}({A})}(\zeta,\tilde \zeta)$, respectively (whose existence is ensured by \cite[Theorem 4.1]{villani2008optimal}). 
				Then we define $\Xi \in {\cal P}((S\times A)^2)$ by
				\begin{align*}
					\Xi (ds,da,d\tilde s,d\tilde a):=\Upsilon_\zeta(d\tilde a|a)\Lambda_\mu(da|s)\Gamma (d s,d \tilde s),
				\end{align*}
				where $\Lambda_\mu:S\ni s \mapsto \Lambda_\mu(da|s)\in {\cal P}(A)$ denotes a disintegrating kernel of $\Lambda$ with respect to its marginal $\mu=\operatorname{pj}_S(\Lambda)$, i.e.,
				\begin{align}\label{eq:dis_integ}
					\Lambda(ds,da)=\Lambda_\mu(da|s)\mu(ds).
				\end{align} 
				In a similar manner, $\Upsilon_\zeta:A\ni a \mapsto \Upsilon_\zeta(d\tilde a |a)\in {\cal P}(A)$ denotes a disintegrating kernel of~$\Upsilon$ with respect to its marginal $\zeta=\operatorname{pj}_A(\Upsilon)$.
				
				Then, by \eqref{eq:const1} and \eqref{eq:dis_integ}, it holds that $\int_{(\tilde s, \tilde a)\in S\times A}\Xi(ds,da,d\tilde s,d\tilde a)= \Lambda(ds,da).$ Moreover by setting $\tilde \Lambda^{\diamond}(d\tilde s, d\tilde a):= \int_{(s, a)\in S\times A}
				\Xi(ds,da,d\tilde s,d\tilde a)$, we have that 
				\begin{align*}
					\tilde \Lambda^\diamond\in \operatorname{Cpl}_{S\times A}(\tilde \mu,\tilde \zeta),\qquad \Xi\in \operatorname{Cpl}_{(S\times A)^2}(\Lambda,\tilde \Lambda^\diamond).
				\end{align*}
				This implies  that
				\begin{align*}
					\begin{aligned}
						\inf_{\tilde \Lambda \in \operatorname{Cpl}_{S\times A}(\tilde \mu,\tilde \zeta)}{\cal W}_{{\cal P}({S\times A})}(\Lambda,\tilde \Lambda)\leq {\cal W}_{{\cal P}({S\times A})}(\Lambda, \tilde \Lambda^{\diamond})
						&\leq  \int_{(S\times A)^2} d_{S\times A}((s,a),(\tilde s,\tilde a))\Xi (ds,da,d\tilde s, d\tilde a)\\
						=&\int_{S\times S}d_S(s,\tilde s)\Gamma(ds,d\tilde s)+\int_{A\times A}d_A(a,\tilde a)\Upsilon(da,d\tilde a)\\
						=& {\cal W}_{{\cal P}(S)}(\mu,\tilde \mu)+ {\cal W}_{{\cal P}(A)}(\zeta,\tilde \zeta),
					\end{aligned}
				\end{align*}
				where the last equality follows from the optimality of $\Gamma$ and $\Upsilon$ (see \eqref{eq:const1}).
				
				Combining this with the compactness of $\operatorname{Cpl}_{S\times A}(\tilde \mu,\tilde \zeta)$ (see \cite[Theorem 4.1 \& Lemma 4.4]{villani2008optimal}), one can choose $\tilde \Lambda^* \in \operatorname{Cpl}_{S\times A}(\tilde \mu,\tilde \zeta)$ so that 
				\[
				{\cal W}_{{\cal P}({S\times A})}(\Lambda,\tilde \Lambda^*)=\inf_{\tilde \Lambda \in \operatorname{Cpl}_{S\times A}(\tilde \mu,\tilde \zeta)}{\cal W}_{{\cal P}({S\times A})}(\Lambda,\tilde \Lambda)\leq {\cal W}_{{\cal P}({S})}(\mu,\tilde \mu)+ {\cal W}_{{\cal P}({A})}(\zeta,\tilde \zeta),
				\]
				as claimed.
				
				\vspace{0.5em}
				\noindent Next we prove the part (ii). Let $\Lambda,\tilde \Lambda \in {\cal P}(S\times A)$. Denote by $\Xi^*\in \operatorname{Cpl}_{(S\times A)^2}(\Lambda,\tilde \Lambda )$ the optimal coupling for ${\cal W}_{{\cal P}(S\times A)}(\Lambda,\tilde \Lambda)$. By setting $h(s,a):=s$ for every $(s,a)\in S\times A$ (i.e., a projection map onto $S$), denote by 
				\[
				\Xi^\diamond:= \big(\Xi^*\circ (h\times h)^{-1} \big)\in {\cal P}(S\times S)
				\]
				the push-forward of $\Xi^*$ by the map $(h\times h):(S\times A)^2\to S^2$. 
				
				Clearly $\Xi^\diamond$ is in $\operatorname{Cpl}_{S\times S}(\operatorname{pj}_S(\Lambda),\operatorname{pj}_S(\tilde \Lambda))$. Thus, 
				\begin{align*}
					{\cal W}_{{\cal P}(S)}(\operatorname{pj}_S(\Lambda),\operatorname{pj}_S(\tilde \Lambda))\leq \int_{S\times S}d_S(s,\tilde s) \Xi^\diamond (ds,d\tilde s) 
					= \int_{(S\times A)^2}d_S(h(s,a),h(\tilde s,\tilde a)) \Xi^*(ds,da,d\tilde s,d\tilde a).
				\end{align*}
				Moreover, since $d_S(h(s,a),h(\tilde s,\tilde a))= d_S(s,\tilde s)\leq d_{S\times A}((s,a),(\tilde s, \tilde a))$ for every $(s,a),(\tilde s, \tilde a)\in S\times A$, by the optimality of $\Xi^*\in \operatorname{Cpl}_{(S\times A)^2}(\Lambda,\tilde \Lambda )$, the assertion for the part (ii) holds, as claimed.
			\end{proof}
			The following lemma provides useful properties of the lifted functions defined in Definition~\ref{dfn:lift_maps}.
			\begin{lem}\label{lem:corr_S}
				Suppose that Assumption \ref{as:general2}\;(ii),\;(iii) are satisfied.~Let $\mathfrak{U}$, $\overline{\operatorname{F}}$, $\overline{r}$ be given in~Definition~\ref{dfn:lift_maps}. Then the following hold:
				\begin{itemize}
					\item [(i)] $\mathfrak{U}$ is non-empty, compact-valued and continuous.\footnote{A correspondence between topological spaces is continuous if it is both lower- and upper-hemicontinuous (see,
						e.g., \cite[Definition 17.2, p.\;558]{CharalambosKim2006infinite}).}
					\item [(ii)] $\overline{\operatorname{F}}$ satisfies that for every $(\mu,\Lambda,e^0), (\tilde \mu,\tilde \Lambda,\tilde e^0)\in \operatorname{gr}(\mathfrak{U})\times E^0$,
					\[
					\hspace{3.em}{\cal W}_{{\cal P}(S)}\big(\overline{\operatorname{F}}(\mu,\Lambda,e^0),\overline{\operatorname{F}}(\tilde \mu,\tilde \Lambda,\tilde e^0)\big)\leq \overline C_{{\operatorname{F}}}\big(2 {\cal W}_{{\cal P}(S\times A)}(\Lambda,\tilde \Lambda)+d_{E^0}(e^0,\tilde e^0)\big).
					\]
					\item [(iii)] $\overline{r}$ is bounded. Furthermore, for every $(\mu,\Lambda), (\tilde \mu,\tilde \Lambda)\in \operatorname{gr}(\mathfrak{U})$
					\[
					|\overline r(\mu,\Lambda)-\overline r(\tilde \mu,\tilde \Lambda)| \leq 2\overline C_r {\cal W}_{{\cal P}(S\times A)}(\Lambda,\tilde \Lambda).
					\]
				\end{itemize}
			\end{lem}
			\begin{proof}
				We start by proving (i). Both the non-emptyness and the compact-valuedness of $\mathfrak{U}$ are clear. Indeed, for every $\mu \in {\cal P}(S)$ one can consider the Dirac measure~$\delta_{\hat a}(da)\in {\cal P}(A)$ at some $\tilde a \in A$ to obtain that $\delta_{\tilde a}(da)\mu(ds)\in\mathfrak{U}(\mu).$ Therefore $\mathfrak{U}(\mu)$ is non-empty. 
				
				Moreover, since $\operatorname{pj}_S: {\cal P}(S\times A)\to {\cal P}(S)$ is continuous (see Lemma \ref{lem:element_1Wass}\;(ii)) and ${\cal P}(S\times A)$ is compact (as $S\times A$ is compact), $\mathfrak{U}(\mu)\subseteq {\cal P}(S\times A)$ is compact for every $\mu\in {\cal P}(S)$, as claimed. 
				
				\vspace{0.5em}
				We now claim that $\mathfrak{U}$ is both upper and lower hemicontinuous. Let $\mu\in {\cal P}(S)$ be given. 
				
				Recalling that $\operatorname{gr}(\mathfrak{U})=\{(\mu,\Lambda)\in {\cal P}(S)\times {\cal P}(S\times A)\;|\;\Lambda \in \mathfrak{U}(\mu) \}$, let us consider a sequence $(\mu^{(n)},\Lambda^{(n)})_{n\in \mathbb{N}}\in\operatorname{gr}(\mathfrak{U})$ such that $\mu^{(n)}\rightharpoonup \mu$ as $n\to \infty$. Since the subset $\operatorname{gr}(\mathfrak{U})\subseteq {\cal P}(S)\times {\cal P}(S\times A)$ is compact (by the continuity of $\operatorname{pj}_S:{\cal P}(S\times A)\to {\cal P}(S)$ and the compactness of ${\cal P}(S)\times {\cal P}(S\times A)$), there exists a subsequence 
				\[
				(\mu^{(n_k)},\Lambda^{(n_k)})_{k\in \mathbb{N}} \subseteq (\mu^{(n)},\Lambda^{(n)})_{n\in \mathbb{N}}\quad \mbox{s.t. $(\mu^{(n_k)},\Lambda^{(n_k)})\rightharpoonup (\mu^\star,\Lambda^\star)$ as $k\to \infty$}
				\]
				with some $(\mu^\star,\Lambda^\star)\in \operatorname{Gr}(\mathfrak{U})$. Combined with the limit $\mu^{(n)}\rightharpoonup \mu=\mu^\star$, this ensures that $(\Lambda^{(n)})_{n\in\mathbb{N}}$ has a limit point $\Lambda^{\star}\in\mathfrak{U}(\mu)= \mathfrak{U}(\mu^\star)$. Thus, by \cite[Theorem~17.20]{CharalambosKim2006infinite}, $\mathfrak{U}$ is upper hemicontinuous.
				
				It remains to show the lower hemicontinuity of $\mathfrak{U}$. First note that for every~$\mu\in {\cal P}(S)$ the set $\mathfrak{U}(\mu)\subseteq {\cal P}(S\times A)$ can be represented by
				\begin{align}\label{eq:union_lhc}
					\mathfrak{U}(\mu)=\bigcup_{\zeta \in {\cal P}(A)}\operatorname{Cpl}_{S\times A}(\mu,\zeta).
				\end{align}
				
				Then we claim that $\operatorname{Cpl}_{S\times A}:{\cal P}(S)\times {\cal P}(A)\ni(\mu,\zeta)\twoheadrightarrow \operatorname{Cpl}_{S\times A}(\mu,\zeta)\subseteq {\cal P}(S\times A)$ is lower-hemicontinuous. To that end, let $(\mu,\zeta)\in {\cal P}(S)\times {\cal P}(A)$ and $\Lambda \in \operatorname{Cpl}_{S\times A}(\mu,\zeta)$ be given, and consider a sequence $(\mu^{(n)},\zeta^{(n)})_{n\in\mathbb{N}}\subseteq {\cal P}(S)\times {\cal P}(A)$ such that $(\mu^{(n)},\zeta^{(n)})\rightharpoonup (\mu,\zeta)$ as $n\to \infty$. 
				
				By Lemma \ref{lem:element_1Wass}, for every $n\in \mathbb{N}$ there exists $\Lambda^{(n),*}\in \operatorname{Cpl}_{S\times A}(\mu^{(n)}, \zeta^{(n)})$ such that
				\[
				{\cal W}_{{\cal P}(S\times A)}(\Lambda,\Lambda^{(n),*})\leq {\cal W}_{{\cal P}(S)}(\mu,\mu^{(n)})+{\cal W}_{{\cal P}(A)}(\zeta,\zeta^{(n)}).
				\]
				Combined with the limit $(\mu^{(n)},\zeta^{(n)})\rightharpoonup (\mu,\zeta)$, this ensures that $\Lambda^{(n),*}\rightharpoonup \Lambda$ as $n\rightarrow\infty$. Thus, by \cite[Theorem 17.21]{CharalambosKim2006infinite}, $\operatorname{Cpl}_{S\times A}$ is lower hemicontinuous. 
				
				Moreover, by the lower hemicontinuity of $\operatorname{Cpl}_{S\times A}$ and the representation given in \eqref{eq:union_lhc}, \cite[Theorem 17.27]{CharalambosKim2006infinite} asserts that $\mathfrak{U}$ is lower hemicontinuous.  Therefore, $\mathfrak{U}$ is continuous, as claimed. 
				
				\vspace{0.5em}
				\noindent Now we prove the part (ii). Let $(\mu,\Lambda,e^0), (\tilde \mu,\tilde \Lambda,\tilde e^0)\in \operatorname{gr}(\mathfrak{U})\times {\cal P}(E)\times E^0$. For simplicity,~let 
				\begin{align}\label{eq:abbre_msr_mu}
					\begin{aligned}
						\mu':=\overline{\operatorname{F}}(\mu,\Lambda,e^0)
						, \qquad  \tilde \mu':=\overline{\operatorname{F}}(\tilde \mu,\tilde \Lambda,\tilde e^0).
					\end{aligned}
				\end{align}
				Then, set $\operatorname{id}_E:E\ni e\mapsto \operatorname{id}_E(e):=(e,e)\in E^2$. Then we denote the diagonal coupling of $\lambda_\varepsilon$ by
				\begin{align}\label{eq:coupl_diag}
					\Xi_1:= \lambda_\varepsilon \circ (\operatorname{id}_E(\cdot))^{-1} \in \operatorname{Cpl}_{E\times E}(\lambda_\varepsilon,\lambda_\varepsilon)
				\end{align}
				so that ${\cal W}_{{\cal P}(E)}(\lambda_\varepsilon,\lambda_\varepsilon)= \int_{E\times E}d_E(e,\tilde e)\Xi_1(de,d \tilde e)=0$. 
				
				Furthermore, we denote the optimal coupling for ${\cal W}_{{\cal P}(S\times A)}(\Lambda,\tilde\Lambda)$ (see \cite[Theorem 4.1]{villani2008optimal}) by 
				\begin{align}\label{eq:coupl_opt}
					\Xi_2\in \operatorname{Cpl}_{(S\times A)^2}(\Lambda,\tilde \Lambda).
				\end{align}
				Using the couplings $\Xi_1$ and $\Xi_2$, we define a coupling $\Xi_3\in \operatorname{Cpl}_{(S\times A\times E)^2}(\Lambda\otimes \lambda_\varepsilon,\tilde \Lambda\otimes \lambda_\varepsilon)$ by
				\begin{align}\label{eq:coupl_prod}
					\Xi_3(ds,da,de,d\tilde s, d\tilde a, d\tilde e):= \Xi_1(de,d\tilde e)\Xi_2(ds,da,d\tilde s,d\tilde a).
				\end{align}
				
				By the definition of~$\overline{\operatorname{F}}$ (see Definition \ref{dfn:lift_maps}\;(ii)) and the setting \eqref{eq:abbre_msr_mu}, it holds~that
				\[
				\Xi_3\circ \big(\operatorname{F}(\cdot,\cdot,\Lambda,\cdot,e^0)\times \operatorname{F}(\cdot,\cdot,\tilde \Lambda,\cdot,\tilde e^0)\big)^{-1}\in \operatorname{Cpl}_{S\times S}(\mu',\tilde \mu'),
				\]
				i.e., the push-forward of $\Xi_3$ by $\operatorname{F}(\cdot,\cdot,\Lambda,\cdot,e^0)\times \operatorname{F}(\cdot,\cdot,\tilde \Lambda,\cdot,\tilde e^0):(S,A,E)^2\to S^2$.  
				
				Then it holds that
				\begin{align}
					{\cal W}_{{\cal P}(S)}(\mu',\tilde \mu')
					&\leq  \int_{S\times S }  d_S(s,s') \big(\Xi_3\circ (\operatorname{F}(\cdot,\cdot,\Lambda,\cdot,e^0)\times \operatorname{F}(\cdot,\cdot,\tilde \Lambda,\cdot,\tilde e^0)\big)^{-1}\big)(ds,ds')\nonumber \\
					&= \int_{(S\times A\times E)^2}d_S(\operatorname{F}(s,a,\Lambda,e,e^0),\operatorname{F}(\tilde s,\tilde a,\tilde \Lambda,\tilde e,\tilde e^0))\Xi_3(ds,da,de,d\tilde s, d\tilde a, d\tilde e)\label{eq:est_1} \\
					&=\int_{(S\times A)^2} \int_{E} d_S(\operatorname{F}(s,a,\Lambda,e,e^0),\operatorname{F}(\tilde s,\tilde a,\tilde \Lambda,e,\tilde e^0)) \lambda_\varepsilon(de) \Xi_2(ds,da,d\tilde s,d\tilde a)
					=:\operatorname{I}, \nonumber
				\end{align}
				where the last line follows from the definition of $\Xi_1$ and $\Xi_3$ (see \eqref{eq:coupl_diag}, \eqref{eq:coupl_prod}) and by applying Fubini's theorem (noting that $\operatorname{F}$ maps into the compact space $S$).
				
				By Assumption \ref{as:general2}\;(ii) and the triangle inequality, 
				\begin{align*}
					\begin{aligned}
						\operatorname{I}\leq & \overline C_{\operatorname{F}} \bigg( \int_{(S\times A)^2}d_{S\times A}\big((s,a),(\tilde s,\tilde a)\big)\Xi_2(ds,da,d\tilde s, d\tilde a)+{\cal W}_{{\cal P}(S\times A)}(\Lambda,\tilde \Lambda)+d_{E^0}(e^0,\tilde e^0 ) \bigg)\\
						=&\overline C_{\operatorname{F}} \big(2{\cal W}_{{\cal P}(S\times A)}(\Lambda, \tilde \Lambda)+d_{E^0}(e^0,\tilde e^0)\big),
					\end{aligned}
				\end{align*}
				where the last equality follows from the optimality of $\Xi_2$ (see \eqref{eq:coupl_opt}).

				Combined with \eqref{eq:est_1}, this ensures the estimates for $\overline{\operatorname{F}}$ to hold.

				\vspace{0.5em}
				\noindent We next prove the part (iii). Since $S$, $A$, and ${\cal P}(S\times A)$ are all compact and $r$ is continuous (by Assumption \ref{as:general}\;(i) and Assumption \ref{as:general2}\;(iii)), $\overline r$ is bounded. We prove its $2\overline C_r$-Lipschitz continuity. Let $(\mu,\Lambda),(\tilde \mu, \tilde \Lambda)\in \operatorname{gr}(\mathfrak{U})$ be given. Then it follows from Assumption~\ref{as:general2}\;(iii) and the triangle inequality that for every $\Xi \in \operatorname{Cpl}_{S\times A}(\Lambda,\tilde \Lambda)$
				\begin{align*}
					\begin{aligned}
						|\overline r(\mu,\Lambda)-\overline r(\tilde \mu,\tilde \Lambda)| &=\bigg|\int_{(S\times A)^2}\big(r(s,a,\Lambda)-r(\tilde s,\tilde a,\tilde \Lambda) \big)\Xi(ds,da,d\tilde s, d\tilde a) \bigg|\\
						&\leq \overline C_r \bigg(\int_{S\times A} d_{S\times A}\big((s,a),(\tilde s,\tilde a)\big)\Xi(ds,da,d\tilde s, d\tilde a) + {\cal W}_{{\cal P}(S\times A)}(\Lambda, \tilde \Lambda) \bigg).
					\end{aligned}
				\end{align*}
				By taking inifimum over all $\Xi \in\operatorname{Cpl}_{S\times A}(\Lambda,\tilde \Lambda)$ into the above, we can obtain the estimate for~$\overline{r}$.	
				
				This completes the proof. 
			\end{proof}

			Using the two preceding lemmas, we now proceed to prove Proposition \ref{pro:dpp}. 
			\begin{proof}[Proof of Proposition \ref{pro:dpp}]We start by proving (i). Let $L\geq0$ and $\overline{V}\in \operatorname{Lip}_{b,L}({\cal P}(S);\mathbb{R})$ be given. Set
				\(
				{\cal S}:={\cal P}(S\times A)\times \mathfrak{P}^0.
				\) 
				Recalling the definition of $\overline p$ (see Definition \ref{dfn:lift_maps}\;(iii)), define $G: {\cal S}\ni(\Lambda,p)\mapsto  G(\Lambda,p)\in \mathbb{R}$~by 
				\begin{align}\label{eq:map_G}
					G(\Lambda,p):= \int_{{\cal P}(S)} \overline{V}(\mu') \overline p(d\mu'|\operatorname{pj}_S(\Lambda),\Lambda,p)
					= \int_{E^0} \overline{V}(\overline{\operatorname{F}}(\operatorname{pj}_S(\Lambda),\Lambda,e^0)) p(de^0).
				\end{align}
				
				We claim that $G$ is continuous. Consider a sequence $(\Lambda^{(n)},p^{(n)})_{n\in \mathbb{N}}\subseteq {\cal S}$ such that $(\Lambda^{(n)},p^{(n)})\rightharpoonup (\Lambda^{\star},p^{\star})$ as $n\to \infty$, with some $(\Lambda^{\star},p^{\star})\in {\cal S}$.
				
				By the triangle inequality, for every $n\in \mathbb{N}$,
				\begin{align*}
					\big|G(\Lambda^{(n)},p^{(n)})-G(\Lambda^{\star},p^{\star})\big|&\leq \big|G(\Lambda^{\star},p^{(n)})-G(\Lambda^{\star},p^{\star})\big|+\big|G(\Lambda^{(n)},p^{(n)})-G(\Lambda^{\star},p^{(n)})\big|\\
					&=:\operatorname{I}^{(n)}+\operatorname{II}^{(n)}.
				\end{align*}
				We will show that $\operatorname{I}^{(n)}$and $\operatorname{II}^{(n)}$ vanish as $n\rightarrow \infty$. 
				
				Since $\overline{V}\in \operatorname{Lip}_{b,L}({\cal P}(S);\mathbb{R})$ and $\overline{\operatorname{F}}$ is continuous (see Lemma \ref{lem:corr_S}\;(ii)), it holds that
				$g^\star(\cdot):=\overline{V}(\overline{\operatorname{F}}(\operatorname{pj}_S(\Lambda^\star),\Lambda^{\star},\cdot))\in C_b(E_0;\mathbb{R})$. Combined with the limit $p^{(n)}\rightharpoonup p^{\star}$, this ensures that
				\begin{align*}
					\lim_{n\to \infty}\operatorname{I}^{(n)}=\lim_{n\to \infty} \bigg|\int_{E^0}g^\star(e^0)p^{(n)}(d e^0) -\int_{E^0} g^\star(\tilde e^0)p^{\star}(d\tilde e^0)\bigg|= 0.
				\end{align*}
				
				It remains to show the limit of $\operatorname{II}^{(n)}$. We use the $L$-Lipschitz continuity of $\overline{V}$, the estimate of~$\overline{\operatorname{F}}$ given in Lemma \ref{lem:corr_S}\;(ii), and the limits $\Lambda^{(n)}\rightharpoonup\Lambda^\star$ and $p^{(n)}\rightharpoonup p^{\star}$ to obtain
				\begin{align*}
					\begin{aligned}
						\lim_{n\rightarrow\infty }\operatorname{II}^{(n)} &\leq  \lim_{n\rightarrow\infty}\int_{E^0}\Big|\overline{V}\big(\overline{\operatorname{F}}(\operatorname{pj}_S(\Lambda^{(n)}),\Lambda^{(n)},e^0)\big)-\overline{V}\big(\overline{\operatorname{F}}(\operatorname{pj}_S(\Lambda^{\star}),\Lambda^{\star},e^0)\big) \Big| p^{(n)}(de^0)\\
						&\leq  2L \overline C_{{ \operatorname{F}}}  \lim_{n\rightarrow\infty}{\cal W}_{{\cal P}(S\times A)}(\Lambda^{(n)},\Lambda^\star) =0.
					\end{aligned}
				\end{align*}
				Therefore $G$ given in \eqref{eq:map_G} is continuous, as claimed. 
				
				Since $\mathfrak{P}^0$ is compact (see Assumption \ref{as:general}\;(i)) and $G$ is continuous, an application of Berge's maximum theorem (see, e.g., \cite[Theorem 17.31]{CharalambosKim2006infinite}) ensures the continuity of the map $\overline{J}:{\cal P}(S\times A)\ni \Lambda \mapsto \overline{J}(\Lambda)\in \mathbb{R}$ given by 
				\begin{align}\label{eq:aux_local_J}
					\overline{J}(\Lambda):= \inf_{p \in \mathfrak{P}^0}\int_{{\cal P}(S)} \overline V(\mu') \overline p(d\mu'|\operatorname{pj}_S(\Lambda),\Lambda,p),
				\end{align}
				and the existence of the measurable selector $\overline{p}^*: {\cal P}(S\times A)\ni \Lambda\mapsto \overline{p}^*(\Lambda)\in \mathfrak{P}^0$ satisfying~\eqref{eq:local_min}. 
				
				\vspace{0.5em}
				\noindent We now prove the part (ii). In analogy to the part (i), the key idea is to apply Berge's maximum theorem. To that end, we first show that a map $H: \operatorname{gr}(\mathfrak{U})\in (\mu,\Lambda)\mapsto H(\mu,\Lambda)\in\mathbb{R}$ defined by
				\begin{align}\label{eq:map_H}
					H(\mu,\Lambda):=\overline{r}(\mu,\Lambda)+\beta\cdot \overline{J}(\Lambda),
				\end{align}
				with $\overline{J}:{\cal P}(S\times A) \to  \mathbb{R}$ defined in \eqref{eq:aux_local_J} is continuous. That will be achieved in two steps. 
				
				\vspace{0.5em}
				Consider a sequence $(\mu^{(n)},\Lambda^{(n)})_{n\in \mathbb{N}}\subseteq \operatorname{gr}(\mathfrak{U})$ such that $(\mu^{(n)},\Lambda^{(n)})\rightharpoonup (\mu^\star,\Lambda^\star)$ as $n\to \infty$, with some $(\mu^\star,\Lambda^\star)\in \operatorname{gr}(\mathfrak{U})$. By the triangle inequality, it holds that for every $n\in \mathbb{N}$,
				\begin{align*}
					|H(\mu^{(n)},\Lambda^{(n)})-H(\mu^{\star},\Lambda^{\star})| &\leq |\overline{r}(\mu^{(n)},\Lambda^{(n)})-\overline{r}(\mu^{\star},\Lambda^{\star})|+\beta\cdot |\overline{J}(\Lambda^{(n)})-\overline{J}(\Lambda^\star)|\\
					&=: \operatorname{III}^{(n)}+\beta \cdot|\operatorname{IV}^{(n)}|.
				\end{align*}
				
				The limit of $\operatorname{III}^{(n)}$ is straightforward. Indeed, by Lemma \ref{lem:corr_S}\;(iii) and the limit~$\Lambda^{(n)} \rightharpoonup \Lambda^\star$,
				\begin{align*}
					\lim_{n\to\infty}\operatorname{III}^{(n)} \leq 2\overline C_r \lim_{n\to \infty} {\cal W}_{{\cal P}(S\times A)}(\Lambda^{(n)},\Lambda^\star)=0.
				\end{align*}
				
				It remains to show the limit of $|\operatorname{IV}^{(n)}|$. Recalling the measuarable selector $\overline{p}^*$ defined as in the part\;(i), denote by $p^\star:=\overline{p}^*(\Lambda^\star)\in \mathfrak{P}^0$. Then it holds that
				\begin{align}\label{eq:local_nu1}
					\overline{J}(\Lambda^\star)
					=\int_{{\cal P}(S)} \overline V(\mu') \overline p(d\mu'|\operatorname{pj}_S(\Lambda^\star),\Lambda^\star,p^{\star})=\int_{E^0}\overline{V}(\overline{\operatorname{F}}(\mu^\star,\Lambda^\star,e^0))p^{\star}(de^0),
				\end{align}
				noting that $\operatorname{pj}_S(\Lambda^\star)=\mu^\star$ as $(\mu^\star,\Lambda^\star)\in \operatorname{gr}(\mathfrak{U})$.
				
				On the other hand, as $p^\star\in \mathfrak{P}^0$ does not necessarily optimize $\overline J(\Lambda^{(n)})$, it holds that
				\begin{align}\label{eq:local_nu2}
					\overline{J}(\Lambda^{(n)})
					\leq \int_{{\cal P}(S)} \overline V(\mu') \overline p(d\mu'|\operatorname{pj}_S(\Lambda^{(n)}),\Lambda^{(n)},p^\star)=\int_{E^0}\overline{V}(\overline{\operatorname{F}}(\mu^{(n)},\Lambda^{(n)},e^0))p^{\star}(de^0),
				\end{align}
				with $\operatorname{pj}_S(\Lambda^{(n)})=\mu^{(n)}$.
				
				By \eqref{eq:local_nu1} and \eqref{eq:local_nu2}, it holds that for every $n\in\mathbb{N}$ and every $\Gamma \in \operatorname{Cpl}_{E^0\times E^0}(p^{\star},p^{\star} )$,
				\begin{align}\label{eq:local_nu2_1}
					\begin{aligned}
						\operatorname{IV}^{(n)} &\leq \int_{E^0}\overline{V}(\overline{\operatorname{F}}(\mu^{(n)},\Lambda^{(n)},e^0))p^{\star}(de^0)-\int_{E^0}\overline{V}(\overline{\operatorname{F}}(\mu^\star,\Lambda^\star,e^0))p^{\star}(de^0)\\
						&=\int_{E^0\times E^0} \Big(\overline{V}(\overline{\operatorname{F}}(\mu^{(n)},\Lambda^{(n)},e^0))-\overline{V}(\overline{\operatorname{F}}(\mu^\star,\Lambda^\star,\tilde e^0)) \Big) \Gamma (de^0,d\tilde e^0)\\
						&\leq 2L\overline C_{{\operatorname{F}}} \cdot \bigg( {\cal W}_{{\cal P}(S\times A)}(\Lambda^{(n)},\Lambda^\star)+ \int_{E^0\times E^0} d_{E^0}(e^0,\tilde e^0) \Gamma (de^0,d\tilde e^0) \bigg),
					\end{aligned}
				\end{align}
				where the last inequality follows from the $L$-Lipschitz continuity of $\overline V$ and the estimate of $\overline{\operatorname{F}}$ given in Lemma \ref{lem:corr_S}\;(ii).
				
				By taking infimum over $\Gamma\in\operatorname{Cpl}_{E^0\times E^0}(p^{\star},p^{\star} )$ in the last equation of \eqref{eq:local_nu2_1}, we have 
				\begin{align}\label{eq:local_nu3}
					\operatorname{IV}^{(n)} \leq 2L\overline C_{{\operatorname{F}}}  {\cal W}_{{\cal P}(S\times A)}(\Lambda^{(n)},\Lambda^\star).
				\end{align}
				
				Using the same arguments as presented for \eqref{eq:local_nu3}, one can have the lower bound with the same constant, i.e.,  $\operatorname{IV}^{(n)}\geq -  2L\overline C_{{\operatorname{F}}}  {\cal W}_{{\cal P}(S\times A)}(\Lambda^{(n)},\Lambda^\star).$
				
				Combined with the limit $\Lambda^{(n)}\rightharpoonup \Lambda^\star$, this ensures that $|\operatorname{IV}^{(n)}|$ vanishes as $n\to \infty$. Therefore $H$ given in \eqref{eq:map_H} is continuous as claimed. 
				
				\vspace{0.5em}
				Since $\mathfrak{U}$ is is non-empty, compact-valued, and continuous (see Lemma \ref{lem:corr_S}\;(ii)) and $H$ is continuous, an application of Berge's maximum theorem ensures the continuity of ${\cal T}\overline {V}$ (see \eqref{dfn:oper_T}) and the existence of the measurable selector $\overline{\pi}^*:{\cal P}(S)\ni \mu \mapsto \overline{\pi}^*(\mu)\in \mathfrak{U}(\mu)$ satisfying \eqref{eq:local_max}. This completes the proof.
			\end{proof}

			\subsection{Proof of Proposition \ref{pro:fixed_point}} Let $\overline{V}\in \operatorname{Lip}_{b,\overline{L}}({\cal P}(S);\mathbb{R})$. We claim that ${\cal T}\overline{V}\in \operatorname{Lip}_{b,\overline{L}}({\cal P}(S);\mathbb{R})$. 
				From Lemma \ref{lem:corr_S}\;(iii) and the fact that $\overline{V}\in\operatorname{Lip}_{b,\overline {L}}({\cal P}(S);\mathbb{R})$, the boundedness of ${\cal T}\overline{V}$ is straightforward. To verify the $\overline{L}$-Lipschitz continuity of ${\cal T}\overline{V}$, let $\mu,\tilde \mu\in {\cal P}(S)$ and denote~by
				\begin{align}\label{eq:lips_0}
					\operatorname{D}(\mu,\tilde \mu):={\cal T}\overline{V}(\mu)-{\cal T}\overline{V}(\tilde \mu).
				\end{align}
				
				Then let $\overline{\pi}^*(\mu)\in \mathfrak{U}(\mu)$ be the local maximizer of ${\cal T}\overline{V}(\mu)$ (see Proposition \ref{pro:dpp}\;(ii)). Then, denote by $\zeta^\diamond:=\operatorname{pj}_A(\overline{\pi}^*(\mu)) \in {\cal P}(A)$ the marginal of $\overline{\pi}^*(\mu)\in \mathfrak{U}(\mu)\subset {\cal P}(S\times A)$ on $A$. Since $\overline{\pi}^*(\mu)\in \operatorname{Cpl}_{S\times A}(\mu, \zeta^\diamond)$, by Lemma \ref{lem:element_1Wass}\;(i) there exists a coupling $\tilde \Lambda^\diamond \in \operatorname{Cpl}_{S\times A}(\tilde \mu, \zeta^\diamond)$ such that 
				\begin{align}\label{eq:lips_1}
					{\cal W}_{{\cal P}(S\times A)}(\overline \pi^*(\mu), \tilde \Lambda^\diamond) \leq {\cal W}_{{\cal P}(S)}(\mu,\tilde \mu).
				\end{align}
				
				Then since $\tilde \Lambda^\diamond \in  \mathfrak{U}(\tilde \mu)$ (which does not necessarily maximize ${\cal T}\overline{V}(\tilde \mu)$), it holds that
				\begin{align*}
					\operatorname{D}(\mu,\tilde \mu) \leq 
					\overline{r}(\mu,\overline \pi^*(\mu) ) - \overline{r}(\tilde \mu, \tilde \Lambda^\diamond)+\beta \cdot \overline {J}(\overline \pi^* (\mu))- \beta \cdot\overline {J}(\tilde \Lambda^\diamond)=:\operatorname{D}^1(\mu,\tilde \mu),
				\end{align*}
				recalling $\overline J:{\cal P}(S\times A)\to \mathbb{R}$ defined in \eqref{eq:aux_local_J} (with noting that $\operatorname{pj}_S(\overline \pi^* (\mu))=\mu$ and $\operatorname{pj}_S(\tilde \Lambda^\diamond)=\tilde \mu$).
				
				Let $\overline p^* (\tilde \Lambda^\diamond)\in \mathfrak{P}^0$ be the local minimizers of~$\overline J (\tilde \Lambda^\diamond)$ (see Proposition \ref{pro:dpp}\;(i)). Since they do not necessarily minimize $\overline  {J}(\overline \pi^*(\mu))$, it holds that  
				\begin{align}\label{eq:D_estimate}
					\begin{aligned}
						\operatorname{D}^1(\mu,\tilde \mu)&\leq \overline{r}(\mu,\overline \pi^*(\mu) ) - \overline{r}(\tilde \mu, \tilde \Lambda^\diamond)+\beta\int_{E^0} \overline V (\overline{\operatorname{F}}(\mu,\overline \pi^* (\mu),e^0)) \overline{p}^*(\tilde \Lambda^\diamond)(de^0)\\
						&\quad -\beta \int_{{\cal P}(S)} \overline V(\overline{\operatorname{F}}(\tilde \mu, \tilde \Lambda^\diamond,\tilde e^0)) \overline{p}^*(\tilde \Lambda^\diamond)(d\tilde e^0)=: \operatorname{D}^2(\mu,\tilde \mu),
					\end{aligned}
				\end{align}
				recalling the definition of $\overline{p}$ given in Definition \ref{dfn:lift_maps}\;(iii).

				Let $\Gamma \in \operatorname{Cpl}_{E^0\times E^0}(\overline{p}^*(\tilde \Lambda^\diamond),\overline{p}^*(\tilde \Lambda^\diamond))$ be some arbitrary. Then, by the estimates for $\overline r$ and $\overline{\operatorname{F}}$ (given in Lemma \ref{lem:corr_S}\;(ii),\;(iii)) and $\overline{V}\in \operatorname{Lip}_{b,\overline{L}}({\cal P}(S);\mathbb{R})$, it holds that
				\begin{align}
					\operatorname{D}^2(\mu,\tilde \mu)
					&\leq \big|\overline{r}(\mu,\overline \pi^*(\mu) ) - \overline{r}(\tilde \mu, \tilde \Lambda^\diamond)\big|+ \beta \int_{E^0\times E^0} \big|\overline V (\overline{\operatorname{F}}(\mu,\overline \pi^* (\mu),e^0))- \overline V (\overline {\operatorname{F}}(\tilde \mu,\tilde \Lambda^\diamond,\tilde e^0)) \big| \Gamma (de^0,d\tilde e^0)\nonumber\\
					&\leq 2\overline C_r {\cal W}_{{\cal P}(S\times A)}(\overline \pi^*(\mu),\tilde \Lambda^\diamond)
					\label{eq:D1_estimate}\\
					&\quad + \overline C_{{\operatorname{F}}} \overline{L}\beta  \bigg( 2 {\cal W}_{{\cal P}(S\times A)}(\overline\pi^*(\mu),\tilde \Lambda^\diamond)+
					\int_{E^0\times E^0} d_{E^0}(e^0, \tilde e^0 ) \Gamma (de^0,d\tilde e^0) \bigg). \nonumber
				\end{align}
				
				For the last line of \eqref{eq:D1_estimate}, we take infimum over all $\Gamma \in \operatorname{Cpl}_{E^0\times E^0}(\overline{p}^*(\tilde \Lambda^\diamond),\overline{p}^*(\tilde \Lambda^\diamond))$ and then use the estimate given in \eqref{eq:lips_1} to obtain  
				\begin{align}\label{eq:D1_estimate2}
					\operatorname{D}^2(\mu,\tilde \mu)
					\leq \big(2\overline C_r +  2 \overline C_{{\operatorname{F}}} \overline{L} \beta \big){\cal W}_{{\cal P}(S)}(\mu,\tilde \mu)\leq \overline{L} {\cal W}_{{\cal P}(S)}(\mu,\tilde \mu),
				\end{align}
				where the last inequality holds by the inequality $\overline{L}\geq 2\overline C_r/(1-2 \overline C_{{\operatorname{F}}}\beta)$ with $2\overline C_{\operatorname{F}}\beta<1$.
				
				By \eqref{eq:lips_0}, \eqref{eq:D_estimate} and \eqref{eq:D1_estimate2}, we have that 
				\[
				{\cal T}\overline{V}(\mu)-{\cal T}\overline{V}(\tilde \mu)=\operatorname{D}(\mu,\tilde \mu)\leq \operatorname{D}^1(\mu,\tilde \mu)\leq \operatorname{D}^2(\mu,\tilde \mu)\leq \overline{L} {\cal W}_{{\cal P}(S)}(\mu,\tilde \mu).
				\]
				Since $\mu,\tilde \mu\in {\cal P}(S)$ are chosen arbitrary, one can have that ${\cal T}\overline{V}(\cdot)$ is $\overline{L}$-Lipschitz continuous. Hence, we conclude that ${\cal T}\overline{V}\in \operatorname{Lip}_{b,\overline L}({\cal P}(S);\mathbb{R})$.

				\vspace{0.5em}
				To verify \eqref{eq:contraction}, let $\overline{V},\overline{W}\in \operatorname{Lip}_{b,\overline L}({\cal P}(S);\mathbb{R})$. By Proposition \ref{pro:dpp}\;(ii), for every $\mu\in {\cal P}(S)$
				\begin{align*}
					|{\cal T}\overline V(\mu)-{\cal T}\overline W(\mu)|
					\leq  \beta\sup_{p \in \mathfrak{P}^0} \int_{{\cal P}(S)} |\overline V(\mu')-\overline W(\mu')  | \overline p(d\mu'|\mu,\overline\pi^*(\mu),p)
					\leq \beta \|\overline V - \overline W\|_{\infty},
				\end{align*}
				which ensures \eqref{eq:contraction} to hold.
				
				\vspace{0.5em}
				Since $\beta <1$ and ${\cal T}(\operatorname{Lip}_{b,\overline L}({\cal P}(S);\mathbb{R}) )\subseteq \operatorname{Lip}_{b,\overline L}({\cal P}(S);\mathbb{R})$, ${\cal T}$ is a contraction on $\operatorname{Lip}_{b,\overline L}({\cal P}(S);\mathbb{R})$. Hence, an application of the Banach's fixed point theorem ensures the existence and uniqueness of $\overline{V}^*\in \operatorname{Lip}_{b,\overline L}({\cal P}(S);\mathbb{R})$ such that for every $\overline{V}\in \operatorname{Lip}_{b,\overline L}({\cal P}(S);\mathbb{R})$, $\overline{V}^* = {\cal T}\overline{V}^* = \lim_{n\to \infty} {\cal T}^n \overline{V}$. This completes the proof. \qed

			\section{Proof of results in Section \ref{subsec:main_thm}}\label{proof:subsec:main_thm}
			We begin by presenting an observation that plays a key role in the proof of Lemmas \ref{lem:worst_lift_dynamics} and~\ref{lem:robust_cont}. Recall the set ${\cal Q}$ given in Definition \ref{dfn:measures} and the filtration $\mathbb{G}=({\cal G}_t)_{t\geq 0}$ given in \eqref{eq:filtrations_represent}.
			\begin{lem}\label{lem:law_universal}
				Denote for every $t\geq0$ by $L^0_{{\cal G}_t}(Z)$ the set of all ${\cal G}_t$ measurable random variables $\zeta_t$  with values in a compact Polish space $Z$. Then for every $\zeta_0\in L^0_{{\cal G}_0}(Z)$ and $\mathbb{P},\widetilde{\mathbb{P}}\in {\cal Q}$, it holds that $\mathscr{L}_{\mathbb{P}}(\zeta_0)= \mathscr{L}_{\widetilde {\mathbb{P}}}(\zeta_0)$. Furthermore, for every  $t\geq 1$, $\zeta_t\in L^0_{{\cal G}_t}(Z)$, and $\mathbb{P},\widetilde{\mathbb{P}}\in {\cal Q}$, it holds that $\mathscr{L}_{\mathbb{P}}(\zeta_t\,|\,\varepsilon^0_{1:t})=\mathscr{L}_{\widetilde {\mathbb{P}}}(\zeta_t\,|\,\varepsilon^0_{1:t})$, $\mathbb{P}$-a.s..
			\end{lem}
			\begin{proof}
				Without loss of generality, we consider the case $t\geq 1$, as the case $t=0$ can be subsumed into it. Then, let $\zeta_t\in L^0_{{\cal G}_t}(Z)$ and  $\mathbb{P},\widetilde{\mathbb{P}}\in {\cal Q}$ be given. 
				
				By the same arguments presented for the proof of Lemma \ref{lem:measurability}\;(ii), $\mathscr{L}_{\mathbb{P}}(\zeta_t\,|\,\varepsilon^0_{1:t})$ and $\mathscr{L}_{\widetilde{\mathbb{P}}}(\zeta_t\,|\,\varepsilon^0_{1:t})$ are ${\cal F}_t^0$ measurable. Hence it suffices to show that for any bounded Borel measurable functions $\hat g_{t}:(E^0)^{t}\to \mathbb{R}$ and $\hat f:Z\to \mathbb{R}$,
				\begin{align*}
					\mathbb{E}^{\mathbb{P}}[\hat g_{t}(\varepsilon_{1:t}^0) \hat f(\zeta_t)]=\mathbb{E}^{\mathbb{P}}\bigg[\hat g_{t}(\varepsilon_{1:t}^0)\int_Z \hat f(\tilde z)\mathscr{L}_{\widetilde{\mathbb{P}}}(\zeta_t\,|\,\varepsilon^0_{1:t})(d\tilde z) \bigg].
				\end{align*}
			
			Note that since $\zeta_t$ is ${\cal G}_t$ measurable, there exists a Borel measurable function $\hat l: G\times \Theta^{t+1}\times E^{t}\times (E^0)^t \to Z$ such that $\zeta = \hat l(\gamma, \vartheta_{0:t}, \varepsilon_{1:t},\varepsilon_{1:t}^0).$
			
			Moreover, since $\varepsilon_{1:t}^0$ is independent of $\gamma, \vartheta_{0:t}, \varepsilon_{1:t}$ (see Remark \ref{rem:well_dfn_1}\;(i)), 
			\begin{align*}
				\mathbb{E}^{\mathbb{P}}[\hat g_{t}(\varepsilon_{1:t}^0) \hat f(\zeta_t)]&=\mathbb{E}^{\mathbb{P}}\big[\hat g_{t}(\varepsilon_{1:t}^0)  \hat f\big(\hat l(\gamma, \vartheta_{0:t}, \varepsilon_{1:t},\varepsilon_{1:t}^0)\big)\big]\\
				&=\int_{(E^0)^t} \hat g_t(e^0_{1:t})  \mathbb{E}^{\mathbb{P}}\big[\hat f\big(\hat l(\gamma, \vartheta_{0:t}, \varepsilon_{1:t},e^0_{1:t})\big)\big] \mathscr{L}_{\mathbb{P}}(\varepsilon_{1:t}^0)(de^0_{1:t})\\
				&=\int_{(E^0)^t} \hat g_t(e^0_{1:t})  \mathbb{E}^{\widetilde {\mathbb{P}}}\big[\hat f\big(\hat l(\gamma, \vartheta_{0:t}, \varepsilon_{1:t},e^0_{1:t})\big)\big] \mathscr{L}_{\mathbb{P}}(\varepsilon_{1:t}^0)(de^0_{1:t})=:\operatorname{I}_t,
			\end{align*}
			where the second equality holds by Fubini's theorem and the last equality follows from the fact that $\mathscr{L}_{\mathbb{P}}(\gamma, \vartheta_{0:t}, \varepsilon_{1:t})=\mathscr{L}_{\widetilde{\mathbb{P}}}(\gamma, \vartheta_{0:t}, \varepsilon_{1:t})$ (see Remark \ref{rem:well_dfn_1}\;(ii)).
			
			Therefore, by definition of $\mathscr{L}_{\widetilde{\mathbb{P}}}(\zeta_t\,|\,\varepsilon^0_{1:t})$ and $\mathscr{L}_{\mathbb{P}}(\varepsilon_{1:t}^0)$, 
			\begin{align*}
				\begin{aligned}
					\operatorname{I}_t&= \int_{(E^0)^t} \hat g_t(e^0_{1:t})  \mathbb{E}^{\widetilde {\mathbb{P}}}\big[\mathbb{E}^{\widetilde {\mathbb{P}}}[\hat f(\zeta ) | \varepsilon_{1:t}^0=e^0_{1:t}]\big] \mathscr{L}_{\mathbb{P}}(\varepsilon_{1:t}^0)(de^0_{1:t})\\
					&= \int_{(E^0)^t} \hat g_t(e^0_{1:t}) \bigg(\int_{Z} \hat f(z ) \mathscr{L}_{\widetilde{\mathbb{P}}}(\zeta | \varepsilon_{1:t}^0=e_{1:t})(dz)\bigg)\mathscr{L}_{\mathbb{P}}(\varepsilon_{1:t}^0)(de^0_{1:t})\\
					&=\mathbb{E}^{\mathbb{P}}\bigg[\hat g_{t}(\varepsilon_{1:t}^0)\int_Z \hat f(\tilde z)\mathscr{L}_{\widetilde{\mathbb{P}}}(\zeta_t\,|\,\varepsilon^0_{1:t})(d\tilde z) \bigg],
				\end{aligned}
			\end{align*}
			as claimed. 
			\end{proof}

			\subsection{Proof of Lemma \ref{lem:worst_lift_dynamics}}  We first prove \eqref{eq:worst_mdp}. Let $a\in {\cal A}$ be given. 
				We will construct $\underline{p}_1^{\xi,a} \in \mathfrak{P}^0$ and the sequence of kernels $\underline{p}_t^{\xi,a} :(E^{0})^{t-1}\ni e_{1:t-1}^0\mapsto \underline{p}_t^{\xi,a} (e_t^0|e_{1:t-1}^0)\in\mathfrak{P}^0$ for $t\geq 2$ to define $\underline{\mathbb{P}}^{\xi,a}\in {\cal Q}$ induced by $(\underline{p}_t^{\xi,a})_{t\geq 1}\in  \mathcal{K}^0$.
				
				\vspace{0.5em}
				\noindent{\it Step 1:} Let $\widetilde{\mathbb{P}}\in {\cal Q}$ be some arbitrary. Then set 
				\begin{align}\label{eq:const_init}
					\begin{aligned}
						\tilde s_0:=\xi,\quad 
						\tilde \Lambda_0:=\mathscr{L}_{\widetilde{\mathbb{P}}}((\tilde s_0,a_0)),
					\end{aligned}
				\end{align}
				and define by 
				\begin{align}\label{eq:msr_first}
						\underline p_1^{\xi,a}:=\overline p^*(\tilde \Lambda_0)\in \mathfrak{P}^0,
					\end{align}
					where $\overline{p}^*$ is given in Proposition \ref{pro:dpp}\;(i). 
					
					Next 
					set
					\begin{align}\label{eq:kernel_first0}
						\tilde{s}_1:= \operatorname{F}\big(\tilde s_0,a_0,\tilde \Lambda_0,\varepsilon_{1},\varepsilon_{1}^0\big),\quad  \tilde \Lambda_1:=\mathscr{L}_{\widetilde{\mathbb{P}}}((\tilde s_1,a_1)\,|\,\varepsilon_1^0),
					\end{align}
					where $(\tilde s_0,\tilde \Lambda_0)$ are given in \eqref{eq:const_init}. We see that $(\tilde{s}_1,a_1)$ are ${\cal G}_1$ measurable (because $\tilde s_0 \in L_{{\cal F}_0}^0(S)$, $a_0=\pi_0(\gamma,\vartheta_0)$, $a_1=\pi_1(\gamma,\vartheta_{0:1},\varepsilon_1,\varepsilon_{1}^0)$) and $\varepsilon_1^0$ is independent of $(\gamma,\vartheta_{0:1},\varepsilon_1)$ (see Remark \ref{rem:well_dfn_1}\;(iii)). 
					
					Moreover, an application of Lemma~\ref{lem:regular_cond}\;(ii) implies that $\tilde \Lambda_1$ is ${\cal F}_1^0$ measurable,  
					which ensures the existence of a Borel measurable function $l_1:E^0\to {\cal P}(S\times A)$ such~that 
					\begin{align}\label{eq:kernel_first}
						l_1 (\varepsilon_1^0)=\tilde \Lambda_1.
					\end{align}
					From this, define $\underline p_2^{\xi,a}:E^0\ni e_1^0\mapsto  \underline p_2^{\xi,a}(\cdot\,|\,e_1^0)\in{\cal P}(E^0)$ by
					\begin{align}\label{eq:kernel_first1}
						\underline{p}_2^{\xi,a}(\,\cdot\,|\,e_1^0):= \overline p^*\big(l_1 (e_1^0)\big)\in \mathfrak{P}^0.
					\end{align}
					
					Using the same arguments presented for \eqref{eq:kernel_first0}--\eqref{eq:kernel_first1}, for every $t\geq 1$ we inductively set
					\begin{align}\label{eq:kernel_general0}
						\tilde{s}_t:= \operatorname{F}(\tilde s_{t-1},a_{t-1},\tilde \Lambda_{t-1},\varepsilon_{t},\varepsilon_{t}^0),\quad \tilde \Lambda_t:=\mathscr{L}_{\widetilde{\mathbb{P}}}((\tilde s_{t},a_{t})\,|\,\varepsilon_{1:t}^0),
					\end{align}
					where $(\tilde s_t,a_t)$ are ${\cal G}_t$ measurable, and $\tilde \Lambda_t$ is ${\cal F}_t^0$ measurable.
					
					Hence,  there exists a Borel measurable function $l_t:(E^0)^t\to {\cal P}(S\times A)$ such that 
					\begin{align}\label{eq:kernel_general}
						l_t (\varepsilon_{1:t}^0)=\tilde \Lambda_t.
					\end{align}
					From this, define $\underline p_{t+1}^{\xi,a}:(E^0)^t\ni e_{1:t}^0\mapsto  \underline p_{t+1}^{\xi,a}(\cdot\,|\,e_{1:t}^0)\in{\cal P}(E^0)$ by
					\begin{align}\label{eq:kernel_general1}
						\underline{p}_{t+1}^{\xi,a}(\,\cdot\,|\,e_{1:t}^0):= \overline p^*\big(l_t(e_{1:t}^0)\big)\in \mathfrak{P}^0.
					\end{align}
					
					Using $(\underline{p}_t^{\xi,a})_{t\geq 1}\in  \mathcal{K}^0$, constructed via \eqref{eq:msr_first}, \eqref{eq:kernel_first1}, and \eqref{eq:kernel_general1}, we define the measure $\underline{\mathbb{P}}^{\xi,a}\in {\cal Q}$ induced by this sequence. We underline that the existence of such a measure is ensured by Ionescu--Tulcea's theorem (see Remark~\ref{rem:well_dfn_1}), and that the above inductive construction is invariant and can be carried out under any $\widetilde{\mathbb{P}}\in {\cal Q}$.
					
					\vspace{0.5em}
					\noindent{\it Step 2:} Recall for every $t\geq 0$, $\tilde \Lambda_t$ is the conditional joint law of  $(\tilde s_t,a_t)$ given $\varepsilon_{1:t}^0$ under $\widetilde{\mathbb{P}}$, as given in \eqref{eq:const_init}, \eqref{eq:kernel_first0}, and \eqref{eq:kernel_general0}. We claim that for every $t\geq 0$, $\underline{\mathbb{P}}^{\xi,a}$-a.s.
					\begin{align}\label{eq:induct_msr_const}
						s_t^{\xi,a,\underline{\mathbb{P}}^{\xi,a}}=\tilde s_t,\qquad \underline \Lambda_t^{\xi,a}=\tilde \Lambda_t,
					\end{align}
					where $\underline{\Lambda}_{t}^{\xi,a}$ is the conditional joint law of $(s_t^{\xi,a,\underline{\mathbb{P}}^{\xi,a}},a_t)$ given $\varepsilon^0_{1:t}$ under~$\underline{\mathbb{P}}^{\xi,a}$.
					
					The proof uses an induction over $t\geq 0$: For $t=0$, clearly $s_0^{\xi,a,\underline{\mathbb{P}}^{\xi,a}}=\tilde s_0=\xi\in L_{{\cal F}_0}^0(S)$. Moreover, since $a_0$ is ${\cal G}_0$ measurable (noting that ${\cal G}_0=\sigma(\gamma,\vartheta_0)$) and  $\mathscr{L}_{\underline{\mathbb{P}}^{\xi,a}}(\gamma,\vartheta_0)=\mathscr{L}_{\tilde{\mathbb{P}}}(\gamma,\vartheta_0)$ (see Remark \ref{rem:well_dfn_1}\;(ii)), it holds that $\underline \Lambda_0^{\xi,a}=\tilde \Lambda_0$.

					Assume that the induction claim holds true for some $t\geq 0$. For the case $t+1$, by the conditional McKean-Vlasov dynamics given in \eqref{eq:MKV_represent} and the induction hypothesis for $t$, it holds that $\underline{\mathbb{P}}^{\xi,a}$-a.s.,
					\begin{align}\label{eq:state_equal_as}
						\begin{aligned}
							s_{t+1}^{\xi,a,\underline{\mathbb{P}}^{\xi,a}}&=\operatorname{F}(s_{t}^{\xi,a,\underline{\mathbb{P}}^{\xi,a}},a_t,\underline \Lambda_t^{\xi,a},\varepsilon_{t+1},\varepsilon_{t+1}^0)\\
							&=\operatorname{F}(\tilde s_{t},a_t,\tilde \Lambda_t,\varepsilon_{t+1},\varepsilon_{t+1}^0)=\tilde{s}_{t+1},
						\end{aligned}
					\end{align}
					where the second equality holds by the Borel measurability of $\operatorname{F}$ (see Definition \ref{dfn:basic_element}\;(i)), and the last equality holds by definition \eqref{eq:kernel_general0}, as claimed.
					
					We now show that $\underline \Lambda_{t+1}^{\xi,a}=\tilde \Lambda_{t+1}$, $\underline{\mathbb{P}}^{\xi,a}$-a.s.. By ${\cal F}_{t+1}^0$-measurability of  $(\Lambda_{t+1}^{\xi,a},\tilde \Lambda_{t+1})$, it suffices to show that for any bounded Borel measurable functions $\hat g_{t+1}:(E^0)^{t+1}\to \mathbb{R}$ and $\hat f:S\times A\to \mathbb{R}$,
					\begin{align}\label{eq:claim3_worst}
						\mathbb{E}^{\underline{\mathbb{P}}^{\xi,a}}[\hat g_{t+1}(\varepsilon_{1:t+1}^0) \hat f(s_{t+1}^{\xi,a,\underline{\mathbb{P}}^{\xi,a}},a_{t+1}) ] = \mathbb{E}^{\underline{\mathbb{P}}^{\xi,a}}\bigg[\hat g_{t+1}(\varepsilon^0_{1:t+1})\int_{S\times A} f(\tilde s,\tilde a) \tilde \Lambda_{t+1}(d\tilde s,d\tilde a) \bigg].
					\end{align}
					
					Indeed, by \eqref{eq:state_equal_as}, 
					\begin{align*}
						&\mathbb{E}^{\underline{\mathbb{P}}^{\xi,a}}[\hat g_{t+1}(\varepsilon_{1:t+1}^0) \hat f(s_{t+1}^{\xi,a,\underline{\mathbb{P}}^{\xi,a}},a_{t+1}) ]=\mathbb{E}^{\underline{\mathbb{P}}^{\xi,a}}[\hat g_{t+1}(\varepsilon_{1:t+1}^0) \hat f(\tilde s_{t+1},a_{t+1}) ]=:\operatorname{I}^{t+1}.
					\end{align*}
					
					Moreove, since $(\tilde s_{t+1},a_{t+1})$ are ${\cal G}_{t+1}$ measurable (with ${\cal G}_{t+1}=\sigma(\gamma,\vartheta_{0:t+1},\varepsilon_{1:t+1},\varepsilon_{1:t+1}^0)$), an application of Lemma \ref{lem:law_universal} ensures that 
					$\underline{\mathbb{P}}{}^{\xi,a}$-a.s.,
					\[
						\mathscr{L}_{\underline{\mathbb{P}}{}^{\xi,a}}\big((\tilde s_{t+1},a_{t+1})\,|\,\varepsilon^0_{1:t+1}\big)=\mathscr{L}_{\widetilde {\mathbb{P}}}\big((\tilde s_{t+1},a_{t+1})\,|\,\varepsilon^0_{1:t+1}\big)=\tilde \Lambda_{t+1},
					\]
					which implies that $\operatorname{I}^{t+1}$ equals the second term given in \eqref{eq:claim3_worst}, as claimed.
					
					By induction hypothesis, the claim \eqref{eq:induct_msr_const} holds for all $t\geq 0$. 

					\vspace{0.5em}
					\noindent {\it Step 3:} Recall that $\underline{\mathbb{P}}^{\xi,a}\in {\cal Q}$ is the measure induced by $(\underline{p}_t^{\xi,a})_{t\geq 1} \in  \mathcal{K}^0$ given in \eqref{eq:msr_first}, \eqref{eq:kernel_first1}, and \eqref{eq:kernel_general1} (see Step 1).  Then from Remark \ref{rem:well_dfn_1}\;(iii), it holds that  $\underline{\mathbb{P}}^{\xi,a}$-a.s.
					\begin{align}\label{eq:equiv_1}
						\begin{aligned}
							&\mathscr{L}_{\underline{\mathbb{P}}^{\xi,a}}(\varepsilon_{1}^0)=\underline p_1^{\xi,a}\in \mathfrak{P}^0,\\\;\; &\mathscr{L}_{\underline{\mathbb{P}}^{\xi,a}}(\varepsilon_{t}^0|{\cal F}_{t-1}^0)= \underline{p}_t^{\xi,a}(\cdot|\varepsilon_{1:t-1}^0)\in \mathfrak{P}^0\;\; \mbox{for all $t\geq 2$}.
						\end{aligned}
					\end{align}
					Moreover, since $\underline \Lambda_t^{\xi,a}=\tilde \Lambda_t$ $\underline{\mathbb{P}}^{\xi,a}$-a.s.\;for all $t\geq 0$ (see \eqref{eq:induct_msr_const} in Step 2), 
					it holds that $\underline{\mathbb{P}}^{\xi,a}$-a.s.
					\begin{align}\label{eq:equiv_2}
						\underline p_1^{\xi,a}=\overline p^*(\underline \Lambda_0^{\xi,a}),\qquad \underline{p}_t^{\xi,a}(\cdot|\varepsilon_{1:t-1}^0)= \overline p^*\big(\underline \Lambda_{t-1}^{\xi,a}\big) \quad \mbox{for all $t\geq 2$},
					\end{align}
					which ensures \eqref{eq:worst_mdp} to hold, as claimed.
					
					\vspace{0.5em}
					The proof for \eqref{eq:worst_mdp2} is straightforward. Indeed, by \eqref{eq:lift_2} in Proposition \ref{pro:lift_dynamics} it holds that $\underline{\mathbb{P}}^{\xi,a}$-a.s.
					\begin{align*}
						\mathscr{L}_{\underline{\mathbb{P}}^{\xi,a}}(\underline{\mu}_{1}^{\xi,a})&=\overline p(\,\cdot\,|\,\operatorname{pj}_S(\underline{\Lambda}_{0}^{\xi,a}),\,\underline{\Lambda}_{0}^{\xi,a},\,\underline{p}_1^{\xi,a}(\cdot))\\
						\mathscr{L}_{\underline{\mathbb{P}}^{\xi,a}}(\underline{\mu}_{t+1}^{\xi,a})&=\overline p(\,\cdot\,|\,\operatorname{pj}_S(\underline{\Lambda}_{t}^{\xi,a}),\,\underline{\Lambda}_{t}^{\xi,a},\,\underline{p}_t^{\xi,a}(\,\cdot\,|\varepsilon_{1:t-1}^0))\quad \mbox{for all $t\geq 1$}. 
					\end{align*}
					Combined with \eqref{eq:equiv_2}, this ensures \eqref{eq:worst_mdp2} to hold, as claimed. This completes the proof. \qed

				\subsection{Proof of Lemma \ref{lem:robust_cont}}
					We first introduce some kernels used for constructing $a^*\in {\cal A}$. We denote by 
					\begin{align}\label{eq:universal_zero}
						{\cal K}_{S\times A}:S\times {\cal P}(S\times A)\times {\cal P}(S)\ni (s,\Lambda,\mu)\mapsto {\cal K}_{S\times A}(\,\cdot \,|\,s,\Lambda,\mu)\in {\cal P}(A)
					\end{align}
					the universal disintegration kernel (see Lemma \ref{lem:univ_disint}). Then, we define a kernel
					\begin{align}\label{eq:policy_zero}
						\psi^*:S\times {\cal P}(S)\ni (s,\mu) \mapsto \psi^*(\,\cdot\,|s,\mu):={\cal K}_{S\times A}(\,\cdot \,|\,s,\overline \pi^*(\mu),\mu)\in {\cal P}(A),
					\end{align}
					where $\overline \pi^*$ is the local maximizer given in Proposition \ref{pro:dpp}\;(ii). 
					
					Moreover, denote by  
					\begin{align}\label{eq:BlackDubins_A}
						\rho_A:{\cal P}(A)\times [0,1]\ni(\eta,u)\mapsto \rho_A(\eta,u)\in A 
					\end{align}
					the Blackwell--Dubins function of the action space $A$ (see Lemma \ref{lem:BlackDubins}).
					
					\vspace{0.5em}
					\noindent{\it Step 1.}  Let $\widetilde {\mathbb{P}}\in {\cal Q}$ be some arbitrary. We will inductively construct $a^*\in{\cal A}$ over time $t\geq 0$.
					Let 
					\begin{align}\label{eq:robust_const0}
						\begin{aligned}
						\tilde s_0&:=\xi, \quad&&\tilde \mu_0:= \mathscr{L}_{\widetilde {\mathbb{P}}}(\tilde s_0),\\
						a^*_0&:=\rho_A\big(\psi^*(\cdot\,|\,\tilde s_0,\tilde \mu_0),h_0(\vartheta_0)\big),\quad &&\tilde \Lambda_0:= \mathscr{L}_{\widetilde {\mathbb{P}}}((\tilde s_0,a_0^*)),
						\end{aligned}
					\end{align}
					where $h_0:\Theta\to [0,1]$ is given in Remark \ref{rem:randomize} (so that $h_0(\vartheta_0)\sim {\cal U}_{[0,1]}$).  In particular, since $\xi\in L_{{\cal F}_0}^0(S)$, $\tilde s_0$ is ${\cal F}_0$ measurable, and $a_0^*$ is ${\cal G}_0$ measurable. 
					
					For every $t\geq 1$ we inductively define
					\begin{align}\label{eq:robust_const_general}
						\begin{aligned}
						\tilde s_t&:=\operatorname{F}(\tilde s_{t-1},a_{t-1}^*,\tilde \Lambda_{t-1},\varepsilon_t,\varepsilon_t^0), \quad &&\tilde \mu_t:= \mathscr{L}_{\widetilde {\mathbb{P}}}(\tilde s_t\,|\,\varepsilon_{1:t}^0),\\
						a^*_t&:=\rho_A\big(\psi^*(\,\cdot\,|\,\tilde s_t,\tilde \mu_t),h_t(\vartheta_t)\big), \quad &&\tilde \Lambda_t:=\mathscr{L}_{\widetilde {\mathbb{P}}}((\tilde s_t,a_t^*)\,|\,\varepsilon_{1:t}^0),
						\end{aligned}
					\end{align}
					where $h_t:\Theta \to [0,1]$ is given in Remark \ref{rem:randomize}\;(ii) (so that $(h_u(\vartheta_u))_{0\leq u\leq t}$ is i.i.d.\;with law ${\cal U}_{[0,1]}$). Moreover, by the same arguments presented for the proof of Lemma \ref{lem:measurability}, $\tilde s_t$ is ${\cal F}_t$ measurable, while $a_t^*$ is ${\cal G}_t$ measurable. Moreover, $(\tilde \mu_t,\tilde \Lambda_t)$ are ${\cal F}_t^0$ measurable.
					
					Since $a^*=(a_t^*)_{t\geq 0}$ constructed via \eqref{eq:robust_const0} and \eqref{eq:robust_const_general} is $\mathbb{G}$ adapted, it is in ${\cal A}$.  We underline that the above inductive construction is invariant and can be carried out under any $\widetilde{\mathbb{P}}\in {\cal Q}$.
					
					\vspace{0.5em}
					\noindent {\it Step 2.} We claim that for every $\mathbb{P}\in {\cal Q}$,
					\begin{align}\label{eq:robust_const}
						s_t^{\xi,a^*,{\mathbb{P}}}=\tilde s_t,\quad \mu_t^{\xi,a^*,\mathbb{P}}=\tilde \mu_t,\quad \Lambda_t^{\xi,a^*,\mathbb{P}}=\tilde \Lambda_t, \quad \mbox{$\mathbb{P}$-a.s.},\quad \mbox{for all $t\geq 0$},
					\end{align}
					where $s_t^{\xi,a^*,{\mathbb{P}}}$, ${\mu}_{t}^{\xi,a^*,\mathbb{P}}$, and $\Lambda_t^{\xi,a^*,\mathbb{P}}$ are given in \eqref{eq:MKV_represent}, \eqref{eq:MDP_state} and \eqref{eq:MDP_action}, respectively, under $(a^*,\mathbb{P})$.
					
					Let $\mathbb{P}\in {\cal Q}$ be given. The proof uses an induction over $t\geq 0$: For $t=0$, clearly $s_0^{\xi,a,{\mathbb{P}}}=\tilde s_0=\xi\in L_{{\cal F}_0}^0(S)$. Moreover, since $a_0^*$ is ${\cal G}_0$ measurable (noting that ${\cal G}_0=\sigma(\gamma,\vartheta_0)$) and  $\mathscr{L}_{{\mathbb{P}}}(\gamma,\vartheta_0)=\mathscr{L}_{\tilde{\mathbb{P}}}(\gamma,\vartheta_0)$ (see Remark \ref{rem:well_dfn_1}\;(ii)), it holds that $\mu_0^{\xi,a^*,\mathbb{P}}=\tilde \mu_0$ and  $\Lambda_0^{\xi,a^*,\mathbb{P}}=\tilde \Lambda_0$.

					Assume that the induction claim holds true for some $t\geq 0$. For the case $t+1$, by the conditional McKean-Vlasov dynamics given in \eqref{eq:MKV_represent} and the induction hypothesis for $t$, it holds that $\mathbb{P}$-a.s.,
					\begin{align}\label{eq:state_equal_as.1}
						\begin{aligned}
							s_{t+1}^{\xi,a^*,{\mathbb{P}}}&=\operatorname{F}(s_{t}^{\xi,a^*,{\mathbb{P}}},a^*_t, \Lambda_t^{\xi,a^*,\mathbb{P}},\varepsilon_{t+1},\varepsilon_{t+1}^0)\\
							&=\operatorname{F}(\tilde s_{t},a^*_t,\tilde \Lambda_t,\varepsilon_{t+1},\varepsilon_{t+1}^0)=\tilde{s}_{t+1},
						\end{aligned}
					\end{align}
					where the second equality holds by the Borel measurability of $\operatorname{F}$ (see Definition \ref{dfn:basic_element}\;(i)), and the last equality holds by definition \eqref{eq:robust_const_general}.
					
					We now show that $\Lambda_{t+1}^{\xi,a^*,\mathbb{P}}=\tilde \Lambda_{t+1}$, ${\mathbb{P}}$-a.s.. By ${\cal F}_{t+1}^0$-measurability of  $(\Lambda_{t+1}^{\xi,a^*,\mathbb{P}},\tilde \Lambda_{t+1})$, it suffices to show that for any bounded Borel measurable functions $\hat g_{t+1}:(E^0)^{t+1}\to \mathbb{R}$ and $\hat f:S\times A\to \mathbb{R}$,
					\begin{align}\label{eq:claim3_worst.1}
						\mathbb{E}^{{\mathbb{P}}}[\hat g_{t+1}(\varepsilon_{1:t+1}^0) \hat f(s_{t+1}^{\xi,a^*,\mathbb{P}},a^*_{t+1}) ] = \mathbb{E}^{{\mathbb{P}}}\bigg[\hat g_{t+1}(\varepsilon^0_{1:t+1})\int_{S\times A} f(\tilde s,\tilde a) \tilde \Lambda_{t+1}(d\tilde s,d\tilde a) \bigg].
					\end{align}
					
					Indeed, by \eqref{eq:state_equal_as.1}, 
					\begin{align*}
						\mathbb{E}^{{\mathbb{P}}}[\hat g_{t+1}(\varepsilon_{1:t+1}^0) \hat f(s_{t+1}^{\xi,a^*,\mathbb{P}},a^*_{t+1}) ]=\mathbb{E}^{{\mathbb{P}}}[\hat g_{t+1}(\varepsilon_{1:t+1}^0) \hat f(\tilde s_{t+1},a^*_{t+1}) ]=:\operatorname{I}^{t+1}.
					\end{align*}
					
					Moreover, as $(\tilde s_{t+1},a^*_{t+1})$ is ${\cal G}_{t+1}$ measurable, an application of Lemma \ref{lem:law_universal} ensures that~${\mathbb{P}}$-a.s.
					\[
					\mathscr{L}_{\mathbb{P}}\big((\tilde s_{t+1},a^*_{t+1}) | \varepsilon^0_{1:t+1}\big)=\mathscr{L}_{\widetilde {\mathbb{P}}}\big((\tilde s_{t+1},a^*_{t+1}) | \varepsilon^0_{1:t+1}\big)=\tilde \Lambda_{t+1},
					\]
					which implies that $\operatorname{I}^{t+1}$ equals the second term in \eqref{eq:claim3_worst.1}, as claimed. 
					
					Using the same arguments presented for \eqref{eq:claim3_worst.1}, we have that $\mu_{t+1}^{\xi,a^*,\mathbb{P}}=\tilde \mu_{t+1}$ $\mathbb{P}$-a.s.. Hence, by induction hypothesis, the claim \eqref{eq:robust_const} holds.
					
					\vspace{0.5em}
					\noindent {\it Step 3.} Let $\mathbb{P}\in{\cal Q}$ be some arbitrary. Then we claim that \eqref{eq:law_control} holds. Without loss of generality, we consider the case $t\geq 1$, as the case $t=0$ can be subsumed into~it.  
					
					By the ${\cal F}_{t}^0$-measurability of  $({\Lambda}_{t}^{\xi,a^*,\mathbb{P}},{\mu}_t^{\xi,a^*,\mathbb{P}})$, it suffices to show that for any bounded Borel measurable functions $\hat g_{t}:(E^0)^{t}\to \mathbb{R}$ and $\hat f:S\times A\to \mathbb{R}$,
					\begin{align}\label{eq:claim0_optimal}
						\mathbb{E}^{{\mathbb{P}}}[\hat g_{t}(\varepsilon_{1:t}^0) \hat f(s_{t}^{\xi,a^*,{\mathbb{P}}},a^*_{t}) ] = \mathbb{E}^{{\mathbb{P}}}\bigg[\hat g_{t}(\varepsilon^0_{1:t})\int_{S\times A} f(\tilde s,\tilde a) \overline \pi^*\big({\mu}_t^{\xi,a^*,\mathbb{P}}\big)(d\tilde s,d\tilde a) \bigg],
					\end{align}
					where $\overline \pi^*$ is the local maximizer given in Proposition \ref{pro:dpp}\;(ii).
					
					Since $\hat g(\varepsilon_{1:t}^0)$ is ${\cal F}_t^0$ measurable and it holds that $s_t^{\xi,a^*,{\mathbb{P}}}=\tilde s_t$, ${\mathbb{P}}$-a.s. (see \eqref{eq:robust_const} in Step 2), 
					\begin{align*}
						\begin{aligned}
							\mathbb{E}^{{\mathbb{P}}}[\hat g_{t}(\varepsilon_{1:t}^0) \hat f(s_{t}^{\xi,a^*,{\mathbb{P}}},a^*_{t}) ] &=\mathbb{E}^{{\mathbb{P}}}[\hat g_{t}(\varepsilon_{1:t}^0) \hat f\big(\tilde s_t,a^*_{t}\big) ]\\
							&=\mathbb{E}^{{\mathbb{P}}}\Big[\hat g_{t}(\varepsilon_{1:t}^0)\mathbb{E}^{{\mathbb{P}}}\big[ \mathbb{E}^{{\mathbb{P}}}[\hat f(\tilde s_t,a^*_{t})| {\cal F}_t]\big| {\cal F}_t^0\big] \Big]=:\operatorname{I}_t,
						\end{aligned}
					\end{align*} 
					where the last equality follows from the tower property with fact that ${\cal F}_t^0\subset {\cal F}_t$. 
					
					Since $\tilde s_t$ is ${\cal F}_t$ measurable and $h_t(\vartheta_t)\sim {\cal U}_{[0,1]}$ is independent of ${\cal F}_t$ (noting that ${\cal F}_t$ does not contain the current randomization source $\vartheta_t$), 
					\begin{align}\label{eq:claim2_optimal}
						\begin{aligned}
							\operatorname{I}_t&=\mathbb{E}^{{\mathbb{P}}}\bigg[\hat g_{t}(\varepsilon_{1:t}^0)\mathbb{E}^{{\mathbb{P}}}\bigg[\mathbb{E}^{{\mathbb{P}}}\bigg[\int_A\hat f(\tilde s_t,\tilde a)\psi^*(d\tilde a\,|\,\tilde s_t,\tilde \mu_t)\,\bigg|\,{\cal F}_t\bigg]\bigg|\,{\cal F}_t^0\bigg]\bigg] \\
							&=\mathbb{E}^{{\mathbb{P}}}\bigg[\hat g_{t}(\varepsilon_{1:t}^0)\mathbb{E}^{{\mathbb{P}}}\bigg[ \int_{S\times A}\hat f(\tilde s,\tilde a){\cal K}_{S\times A}(d\tilde a \,|\,\tilde s,\overline \pi^*(\tilde \mu_t),\tilde \mu_t)\tilde \mu_t(d\tilde s)\,\Big|\,{\cal F}_t^0\bigg]\bigg]\\
							&=\mathbb{E}^{{\mathbb{P}}}\bigg[\hat g_{t}(\varepsilon_{1:t}^0) \int_{S\times A}\hat f(\tilde s,\tilde a) \overline{\pi}^*(\tilde \mu_t)(d\tilde s,d\tilde a)\bigg],
						\end{aligned}
					\end{align}
					where the first equality follows from definition of $a_t^*$ given in \eqref{eq:robust_const_general}, the second equality follows from definition of $\psi^*(\cdot|\tilde s_t,\tilde \mu_t)$ (see \eqref{eq:policy_zero}) and ${\cal F}_t^0$-measurability of $\tilde \mu_t$, and the last equality follows from definition of the universal differentiation kernel ${\cal K}_{S\times A}$ (see \eqref{eq:universal_zero}).
					
					Moreover, since $\tilde \mu_t= {\mu}_t^{\xi,a^*,\mathbb{P}}$, $\mathbb{P}$-a.s.~(see \eqref{eq:robust_const}  in Step 2), the last term in \eqref{eq:claim2_optimal} equals the second term in \eqref{eq:claim0_optimal}, as claimed. This completes the proof. \qed

			\subsection{Proof of Theorem \ref{thm:origin_MDP}} 
				For notational simplicity, set $\mu:=\mathscr{L}(\xi)$.
				
				\noindent {\it Step 1:} We claim that  for every $n\in \mathbb{N}$ 
				\begin{align}\label{eq:origin_induct_1}
					{\cal I}_n^{\xi,a^*}:=\inf_{{\mathbb{P}}\in {\cal Q}} \mathbb{E}^{{\mathbb{P}}} \bigg[\sum_{t=0}^{n-1}\beta^t\,{r}(s_{t}^{\xi,a^*,{\mathbb{P}}},a^*_t,{\Lambda}_{t}^{\xi,a^*,\mathbb{P}})+\beta^n\,\overline{V}^*({\mu}_{n}^{\xi,a^*,\mathbb{P}}) \bigg]\geq \overline{V}^*(\mu),
				\end{align}
				where for every $\mathbb{P}\in {\cal Q}$, let $({ \mu}_{t}^{\xi,a^*,\mathbb{P}})_{t\geq 0}$ and $({\Lambda}_{t}^{\xi,a^*,\mathbb{P}})_{t\geq 0}$ be given by \eqref{eq:MDP_state} and \eqref{eq:MDP_action}, respectively. 
				
				We prove \eqref{eq:origin_induct_1} via an induction over $n$. Before proceeding, note that for every $\mathbb{P}\in {\cal Q}$ and~$t\geq 0$,
				\begin{align}
					\begin{aligned}
					\mathbb{E}^{{\mathbb{P}}}\big[{r}(s_{t}^{\xi,a^*,{\mathbb{P}}},a^*_t,{\Lambda}_{t}^{\xi,a^*,\mathbb{P}})\big]
					&=\mathbb{E}^{{\mathbb{P}}}\big[\overline{r}(\operatorname{pj}_S({\Lambda}_t^{\xi,a^*,\mathbb{P}})\,,\,{\Lambda}_t^{\xi,a^*,\mathbb{P}})\big]\\
					&=\mathbb{E}^{{\mathbb{P}}}\big[\overline{r}\big({\mu}_t^{\xi,a^*,\mathbb{P}}\,,\,\overline \pi^*({\mu}_t^{\xi,a^*,\mathbb{P}})\big)\big],\label{eq:lift_reward1}
					\end{aligned}
				\end{align}
				where the first equality holds by 
				\eqref{eq:lift_reward} in Remark \ref{rem:lift_gains} and the second equality follows from \eqref{eq:law_control} in Lemma \ref{lem:robust_cont} and the fact that $\overline \pi^*(\mu)\in \mathfrak{U}(\mu)$ (see Proposition~\ref{pro:dpp}\;(ii)).
				
				Hence by the property \eqref{eq:lift_reward1}, ${\cal I}_n^{\xi,a^*}$ given in \eqref{eq:origin_induct_1} can be represented by 
				\begin{align}\label{eq:origin_induct_1_1}
					{\cal I}_n^{\xi,a^*}=\inf_{{\mathbb{P}}\in {\cal Q}}\mathbb{E}^{{\mathbb{P}}}\bigg[\sum_{t=0}^{n-1}\beta^t\,\overline{r}\big({\mu}_t^{\xi,a^*,\mathbb{P}}\,,\,\overline \pi^*({\mu}_t^{\xi,a^*,\mathbb{P}})\big)+\beta^n\,\overline{V}^*({\mu}_{n}^{\xi,a^*,\mathbb{P}}) \bigg].
				\end{align}
				\noindent{\it Step 1a:} For $n=1$, let ${\mathbb{P}}\in {\cal Q}$ be induced by some $(p_t)_{t\geq 1}\in  \mathcal{K}^0$ (see Definition~\ref{dfn:measures}). 
				
				We first note that ${\mu}_0^{\xi,a^*,\mathbb{P}}=\mu$ with trivial ${\cal F}_0^0$ and $\mathscr{L}_{\mathbb{P}}(\varepsilon_{1}^0)= p_1\in \mathfrak{P}^0$ (see Remark~\ref{rem:well_dfn_1}\;(iii)). Combined with \eqref{eq:lift_2} (see Proposition \ref{pro:lift_dynamics}), this implies that
				\begin{align}
					\mathbb{E}^{{\mathbb{P}}}\big[\overline{r}\big({\mu}_0^{\xi,a^*,\mathbb{P}}\,,\,\overline \pi^*({\mu}_0^{\xi,a^*,\mathbb{P}})\big)+\beta\,\overline{V}^*({\mu}_{1}^{\xi,a^*,\mathbb{P}})\big]&=\overline{r}(\mu,\overline \pi^*(\mu))+ \beta \int_{{\cal P}(S)} \overline{V}^*(\mu') \overline p(d\mu'\big|\mu,\overline \pi^*(\mu),p_1)\nonumber \\
					&\geq \overline{r}(\mu,\overline \pi^*(\mu) ) + \beta \inf_{p\in \mathfrak{P}^0} \int_{{\cal P}(S)} \overline{V}^*(\mu') \overline p(d\mu'\big|\mu,\overline \pi^*(\mu),p)\label{eq:induction_1_msr}\\
					&={\cal T}\overline{V}^*(\mu)=\overline{V}^*(\mu),\nonumber
				\end{align}
				where the last line follows from the optimality of $\overline \pi^*(\mu)\in \mathfrak{U}(\mu)$ for ${\cal T}\overline{V}^*(\mu)$ (see Proposition\;\ref{pro:dpp}\;(ii) for $\overline V^*\in \operatorname{Lip}_{b,\overline L}({\cal P}(S);\mathbb{R})$) and the fixed point result given in Proposition \ref{pro:fixed_point}.
				
				Since \eqref{eq:induction_1_msr} holds for any $\mathbb{P}\in {\cal Q}$, by \eqref{eq:origin_induct_1_1} we have that ${\cal I}_1^{\xi,a^*}\geq \overline V^*(\mu)$.
				
				\vspace{0.5em}
				\noindent{\it Step 1b:} Assume that \eqref{eq:origin_induct_1} holds for some $n\geq 1$. Let ${\mathbb{P}}\in {\cal Q}$ be induced by some $(p_t)_{t\geq 1}\in  \mathcal{K}^0$. 
				
				Note that ${\mu}_n^{\xi,a^*,\mathbb{P}}$ and $\mathscr{L}_{\mathbb{P}}(\varepsilon_{n+1}^0|{\cal F}_n^0)$ are ${\cal F}_n^0$ measurable and  $\mathscr{L}_{\mathbb{P}}(\varepsilon_{n+1}^0|{\cal F}_n^0)=p_{n+1}(\cdot|\varepsilon_{1:n}^0)\in \mathfrak{P}^0$ $\mathbb{P}$-a.s. (see Remark \ref{rem:well_dfn_1}\;(ii)). 
				
				From this, we can use the same arguments presented for \eqref{eq:induction_1_msr} to have that $\mathbb{P}$-a.s.
				\begin{align*}
					\begin{aligned}
						&\mathbb{E}^{{\mathbb{P}}} \big[\overline{r}\big({\mu}_t^{\xi,a^*,\mathbb{P}},\overline \pi^*({\mu}_n^{\xi,a^*,\mathbb{P}})\big)+\beta \overline{V}^*({\mu}_{n+1}^{\xi,a^*,\mathbb{P}}) \big| {\cal F}_n^0\big]\\
						&\quad =\overline{r}({\mu}_n^{\xi,a^*,\mathbb{P}},\overline \pi^*({\mu}_n^{\xi,a^*,\mathbb{P}}))+ \beta \int_{{\cal P}(S)}\overline{V}^*(\mu') \overline p(d\mu'|{\mu}_n^{\xi,a^*,\mathbb{P}},\overline \pi^*({\mu}_n^{\xi,a^*,\mathbb{P}}),p_{n+1}(\cdot|\varepsilon_{1:n}^0))\\
						&\quad \geq \overline{r}({\mu}_n^{\xi,a^*,\mathbb{P}},\overline \pi^*({\mu}_n^{\xi,a^*,\mathbb{P}}))+ \beta \inf_{p\in \mathfrak{P}^0}  \int_{{\cal P}(S)}\overline{V}^*(\mu') \overline p(d\mu'\big|{\mu}_n^{\xi,a^*,\mathbb{P}},\overline \pi^*({\mu}_n^{\xi,a^*,\mathbb{P}}),p)\\
						&\quad = {\cal T}\overline{V}^*({\mu}_n^{\xi,a^*,\mathbb{P}})=\overline{V}^*({\mu}_n^{\xi,a^*,\mathbb{P}}),
					\end{aligned}
				\end{align*}
				which ensures that
				\begin{align}\label{eq:local_n_induct_1}
					\begin{aligned}
						&\mathbb{E}^{{\mathbb{P}}} \bigg[\sum_{t=0}^{n}\beta^t\;\overline{r}({\mu}_t^{\xi,a^*,\mathbb{P}},\overline \pi^*({\mu}_t^{\xi,a^*,\mathbb{P}}))+\beta^{n+1}\overline{V}^*({\mu}_{n+1}^{\xi,a^*,\mathbb{P}}) \bigg]\\
						&\quad \geq  \mathbb{E}^{{\mathbb{P}}} \bigg[\sum_{t=0}^{n-1}\beta^t\;\overline{r}({\mu}_t^{\xi,a^*,\mathbb{P}},\overline \pi^*({\mu}_t^{\xi,a^*,\mathbb{P}}))+\beta^{n}\overline{V}^*({\mu}_n^{\xi,a^*,\mathbb{P}}) \bigg]\geq {\cal I}_n^{\xi,a^*}\geq \overline V^*(\mu),
					\end{aligned}
				\end{align}
				where the second inequality follows from definition of ${\cal I}_n^{\xi,a^*}$ given in \eqref{eq:origin_induct_1_1} and the last inequality follows from assumption of the induction for $n$  (see \eqref{eq:origin_induct_1}).
				
				As \eqref{eq:local_n_induct_1} holds for any $\mathbb{P}\in {\cal Q}$, we have ${\cal I}_{n+1}^{\xi,a^*}\geq \overline V^*(\mu)$. Therefore, by the induction hypothesis, \eqref{eq:origin_induct_1} holds for every $n\in \mathbb{N}$. We conclude that the claim for Step 1 holds.
				
				\vspace{0.5em}
				\noindent {\it Step 2:} We claim that $\overline{V}^*(\mu)\leq {V}(\xi)$. Since $\overline r$ and $\overline{V}^*$ is bounded and $\beta<1$ (see Lemma~\ref{lem:corr_S}\;(iii) and $\overline{V}^*\in\operatorname{Lip}_{b,\overline{L}}({\cal P}(S);\mathbb{R})$), the dominated convergence theorem asserts that for every $\mu\in {\cal P}(S)$
				\begin{align*}
					\begin{aligned}
						\limsup_{n\to \infty}{\cal I}_n^{\xi,a^*}&\leq \inf_{{\mathbb{P}}\in {\cal Q}} \bigg\{ \limsup_{n\to \infty}\mathbb{E}^{{\mathbb{P}}} \bigg[\sum_{t=0}^{n-1}\beta^t\;\overline{r}\big({\mu}_t^{\xi,a^*,\mathbb{P}},\overline \pi^*({\mu}_t^{\xi,a^*,\mathbb{P}})\big) \bigg]+\limsup_{n\to \infty}\mathbb{E}^{{\mathbb{P}}}\big[\beta^n\big|\overline{V}^*({\mu}_n^{\xi,a^*,\mathbb{P}}) \big| \big]\bigg\}\\
						&=\inf_{{\mathbb{P}}\in {\cal Q}}\mathbb{E}^{{\mathbb{P}}} \bigg[\sum_{t=0}^{\infty}\beta^t\;\overline{r}\big({\mu}_t^{\xi,a^*,\mathbb{P}},\overline \pi^*({\mu}_t^{\xi,a^*,\mathbb{P}})\big) \bigg]={\cal J}^{a^*}(\xi)  \leq {V}(\xi),
					\end{aligned}
				\end{align*}
				where the second equality follows from \eqref{eq:lift_reward1} and the definition of ${\cal J}^{a^*}(\xi)$ (see \eqref{eq:worst_represent}).
				
				Combining this with \eqref{eq:origin_induct_1} (as shown in Step\;1),   we conclude that 
				\begin{align}\label{eq:step2_ineq_ext}
					\overline{V}^*(\mu)\leq \limsup_{n\to \infty}{\cal I}_n^{\xi,a^*}\leq{\cal J}^{a^*}(\xi)  \leq {V}(\xi),
				\end{align}
				as claimed.
				
				\vspace{0.5em}
				\noindent{\it Step 3:} 
				We claim that ${V}(\xi) \leq \overline{V}^*(\mu)$, which ensures the statement (i) to hold. For every $a\in{\cal A}$, let $\underline{\mathbb{P}}{}^{\xi,a}\in {\cal Q}$ be induced by $(\underline{p}_t^{\xi,a})_{t\geq 1}\in \mathcal{K}^0$
				such that \eqref{eq:worst_mdp} and \eqref{eq:worst_mdp2} given in Lemma \ref{lem:worst_lift_dynamics} hold. 
				
				Then, define ${{\cal V}}^{a}(\xi)$ by 
				\begin{align*}
					{{\cal V}}^{a}(\xi):=&\mathbb{E}^{\underline{\mathbb{P}}{}^{\xi,a}}\bigg[\sum_{t=0}^\infty \beta^t r(s_t^{\xi,a,\underline{\mathbb{P}}{}^{\xi,a}},a_t,\underline{\Lambda}_t^{\xi,a}) \bigg]
					=\mathbb{E}^{\underline{\mathbb{P}}{}^{\xi,a}} \bigg[\sum_{t=0}^\infty\beta^t\;\overline{r}\big(\operatorname{pj}_S(\underline{\Lambda}_t^{\xi,a}),\underline{\Lambda}_t^{\xi,a} \big)\bigg],
				\end{align*}
				where $\underline{\Lambda}_0^{\xi,a}$ is the joint law of $(s_0^{\xi,a,\underline{\mathbb{P}}{}^{\xi,a}},a_0)$ under $\underline{\mathbb{P}}{}^{\xi,a}$, for $t\geq 1$ $\underline{\Lambda}_t^{\xi,a}$ is the conditional joint law of $(s_t^{\xi,a,\underline{\mathbb{P}}{}^{\xi,a}},a_t)$ under $\underline{\mathbb{P}}{}^{\xi,a}$ given $\varepsilon_{1:t}^0$, and the last equality follows from the same arguments presented for \eqref{eq:lift_reward1}.
				
				Then by definition of ${\cal J}^a(\xi)$ given in \eqref{eq:worst_represent} 
				\begin{align}\label{eq:local_ineq_00_origin}
					V(\xi)=\sup_{a\in {\cal A}} {\cal J}^a(\xi)\leq \sup_{a\in {\cal A}}{{\cal V}}^{a}(\xi).
				\end{align}
				
				Moreover, since $\overline r$ and $\overline{V}^*$ is bounded and $\beta<1$, by the dominated convergence theorem to the sums $\sum_{t=0}^n\beta^t\overline{r}(\operatorname{pj}_S(\underline{\Lambda}_t^{\xi,a}),\underline{\Lambda}_t^{\xi,a} )$ $n\in \mathbb{N}$, we can have that for every $a\in {\cal A}$ 
				\begin{align}\label{eq:local_ineq_0_origin}
					{{\cal V}}^{a}(\xi)=
					\sum_{t=0}^\infty\beta^t  \mathbb{E}^{\underline{\mathbb{P}}{}^{\xi,a}}\big[\overline{r}\big(\operatorname{pj}_S(\underline{\Lambda}_t^{\xi,a}),\underline{\Lambda}_t^{\xi,a} \big)+ \beta \overline V^*(\underline{\mu}_{t+1}^{\xi,a})
					-\beta  \overline V^*(\underline{\mu}_{t+1}^{\xi,a}) \big].
				\end{align}
				
				Then it follows from \eqref{eq:worst_mdp2} in Lemma \ref{lem:worst_lift_dynamics} that for every $t\geq 0$
				\begin{align*}
					\mathbb{E}^{\underline{\mathbb{P}}{}^{\xi,a}}\big[\overline{r}\big(\operatorname{pj}_S(\underline{\Lambda}_t^{\xi,a}),\underline{\Lambda}_t^{\xi,a} \big)+ \beta \overline V^*(\underline{\mu}_{t+1}^{\xi,a})
					\big]&=\mathbb{E}^{\underline{\mathbb{P}}{}^{\xi,a}} \Big[\mathbb{E}^{\underline{\mathbb{P}}{}^{\xi,a}}\big[\overline{r}\big(\operatorname{pj}_S(\underline{\Lambda}_t^{\xi,a}),\underline{\Lambda}_t^{\xi,a} \big)+ \beta \overline V^*(\underline{\mu}_{t+1}^{\xi,a})\;\big|\; {\cal F}_t^0
					\big] \Big],			\\
					&=: \mathbb{E}^{\underline{\mathbb{P}}{}^{\xi,a}}[\overline J(\underline{\Lambda}_t^{\xi,a})] 
				\end{align*}
				where $\overline J(\underline{\Lambda}_t^{\xi,a})$ is ${\cal F}_t^0$ measurable and satisfies
				\begin{align}\label{eq:local_ineq_2_origin}
					\begin{aligned}
						\overline J(\underline{\Lambda}_t^{\xi,a})&=\overline{r}\big(\operatorname{pj}_S(\underline{\Lambda}_t^{\xi,a}),\underline{\Lambda}_t^{\xi,a} \big)+ \beta \int_{{\cal P}(S)} \overline V^*(\tilde \mu)\,\overline p\big(d\tilde \mu\,\big|\,\operatorname{pj}_S(\underline{\Lambda}_{t}^{\xi,a}),\,\underline{\Lambda}_{t}^{\xi,a},\,\overline{p}^*(\underline{\Lambda}_{t}^{\xi,a})\big)\\
						&=\overline{r}\big(\operatorname{pj}_S(\underline{\Lambda}_t^{\xi,a}),\underline{\Lambda}_t^{\xi,a} \big)+ \beta  \inf_{p \in \mathfrak{P}^0}\int_{{\cal P}(S)} \overline V^*(\tilde \mu)\,\overline p\big(d\tilde \mu\,\big|\,\operatorname{pj}_S(\underline{\Lambda}_{t}^{\xi,a}),\,\underline{\Lambda}_{t}^{\xi,a},\,p\big)\\
						&\leq {\cal T}\overline V^*(\operatorname{pj}_S(\underline{\Lambda}_t^{\xi,a})),
					\end{aligned}
				\end{align}
				where the equality holds by the local optimality $\overline p^*(\underline{\Lambda}_{t}^{\xi,a})\in \mathfrak{P}^0$ (see Proposition \ref{pro:dpp}\;(i)) and the inequality holds by definition of ${\cal T}\overline V^*(\operatorname{pj}_S(\underline{\Lambda}_t^{\xi,a}))$ (see \eqref{dfn:oper_T})
				
				Combining \eqref{eq:local_ineq_00_origin}--\eqref{eq:local_ineq_2_origin} with the marginal constraint (i.e., $\operatorname{pj}_S(\underline{\Lambda}_t^{\xi,a})=\underline \mu_t^{\xi,a}$ $\underline {\mathbb{P}}^{\xi,a}$-a.s.; see \eqref{eq:marginal_const}), and the fixed point result (i.e.,  ${\cal T} \overline {V}^*=\overline V^*$; see Proposition \ref{pro:fixed_point}), we conclude that
				\begin{align*}
					V(\xi)\leq \sup_{a\in {\cal A}}\sum_{t=0}^\infty\big(\beta^t \, \mathbb{E}^{\underline{\mathbb{P}}{}^{\xi,a}} [ \overline {V}^*(\underline{\mu}_t^{\xi,a})] -\beta^{t+1} \, \mathbb{E}^{\underline{\mathbb{P}}{}^{\xi,a}}  [\overline V^*(\underline{\mu}_{t+1}^{\xi,a})] \big)= \overline V^*(\mu),
				\end{align*}
				where the last equality holds by the dominated convergence theorem and the fact that $\underline \mu_0^{\xi,a}=\mu$, as claimed.
				
				\vspace{0.5em}
				\noindent {\it Step 4:} It remains to show that \eqref{eq:verify} holds. Recall that $a^*\in {\cal A}$ is such that 
				\eqref{eq:law_control} holds for every $\mathbb{P}\in {\cal Q}$ (see Lemma \ref{lem:robust_cont}). Moreover, let $\underline{\mathbb{P}}^{\xi,a^*}\in {\cal Q}$ is induced by $(\underline {p}_t^{\xi,a^*})_{t\geq 1}\in \mathcal{K}^0$ satisfying \eqref{eq:worst_mdp} and \eqref{eq:worst_mdp2} (see Lemma \ref{lem:worst_lift_dynamics}). 
				
				By applying the dominated convergence theorem to $\sum_{t=0}^n(\beta^t\overline {V}^*(\underline{\mu}_{t}^{\xi,a^*})-\beta^{t+1}\overline {V}^*(\underline{\mu}_{t+1}^{\xi,a^*}))$ $n\in \mathbb{N}$, 
				\begin{align}\label{eq:sum_linear_origin}
					\overline V^*(\mu)=\sum_{t=0}^\infty\big(\beta^t \, \mathbb{E}^{\underline{\mathbb{P}}{}^{\xi,a^*}} \big[ \overline {V}^*(\underline{\mu}_t^{\xi,a^*}) \big] -\beta^{t+1} \, \mathbb{E}^{\underline{\mathbb{P}}{}^{\xi,a^*}}  \big[\overline V^*(\underline{\mu}_{t+1}^{\xi,a^*})\big] \big),
				\end{align}
				where $\underline{\mu}_t^{\xi,a^*}$ is the conditional law of $s_t^{\xi,a^*,\underline{\mathbb{P}}{}^{\xi,a^*}}$ given $\varepsilon_{1:t}^0$.
				
				Note that for every $\mu'\in {\cal P}(S)$
				\begin{align}\label{eq:origin_dpp}
					\overline{V}^*(\mu')={\cal T}\overline{V}^*(\mu')= \overline{r}(\mu',\overline \pi^*(\mu')) + \beta  \int_{{\cal P}(S)} \overline V^*(\tilde \mu')  \overline p\big(d\tilde \mu'|\mu',\overline \pi^*(\mu'),\overline p^*(\overline \pi^*(\mu'))\big).
				\end{align}
				where the first equality follows from Proposition \ref{pro:fixed_point} and the second equality follows from the optimality of the local optimizers $\overline \pi^*$ and $\overline{p}^*$ given in Proposition \ref{pro:dpp}.
				
				From \eqref{eq:origin_dpp}, it holds that for every $t\geq 0$ 
				\begin{align*}
					&\mathbb{E}^{\underline{\mathbb{P}}{}^{\xi,a^*}} \big[ \overline {V}^*(\underline{\mu}_t^{\xi,a^*}) \big]=\mathbb{E}^{\underline{\mathbb{P}}{}^{\xi,a^*}} \big[ {\cal T}\overline {V}^*(\underline{\mu}_t^{\xi,a^*}) \big]\nonumber\\
					&\hspace{2.em} =\mathbb{E}^{\underline{\mathbb{P}}{}^{\xi,a^*}} \bigg[ \overline{r}(\underline{\mu}_t^{\xi,a^*},\overline \pi^*(\underline{\mu}_t^{\xi,a^*})) + \beta  \int_{{\cal P}(S)} \overline V^*(\tilde \mu')  \overline p\big(d\tilde \mu'|\underline{\mu}_t^{\xi,a^*},\overline \pi^*(\underline{\mu}_t^{\xi,a^*}),\overline p^*(\overline \pi^*(\underline{\mu}_t^{\xi,a^*}))\big)\bigg]\\
					&\hspace{2.em} =\mathbb{E}^{\underline{\mathbb{P}}{}^{\xi,a^*}} \bigg[\overline{r}\big(\operatorname{pj}_S(\underline{\Lambda}_t^{\xi,a^*}),\underline{\Lambda}_t^{\xi,a^*} \big)+\beta  \int_{{\cal P}(S)} \overline V^*(\tilde \mu')  \overline p\big(d\tilde \mu'|\operatorname{pj}_S(\underline{\Lambda}_t^{\xi,a^*}),\underline{\Lambda}_t^{\xi,a^*} ,\overline p^*(\underline{\Lambda}_t^{\xi,a^*} )\big) \bigg]=:\operatorname{I}_t,\nonumber
				\end{align*}
				where $\underline{\Lambda}_0^{\xi,a^*}$ is the joint law of $(s_0^{\xi,a^*,\underline{\mathbb{P}}{}^{\xi,a^*}},a^*_0)$ under $\underline{\mathbb{P}}{}^{\xi,a^*}$, for $t\geq 1$ $\underline{\Lambda}_t^{\xi,a^*}$ is the conditional joint law of $(s_t^{\xi,a^*,\underline{\mathbb{P}}{}^{\xi,a^*}},a_t^*)$ under $\underline{\mathbb{P}}{}^{\xi,a^*}$ given $\varepsilon_{1:t}^0$, and the last equality follows from the fact that $\underline{\Lambda}_t^{\xi,a^*}=\overline \pi^*(\underline{\mu}_t^{\xi,a^*})$ $\underline{\mathbb{P}}{}^{\xi,a^*}$-a.s.; see Lemma \ref{lem:robust_cont}, and the marginal constraint that $\operatorname{pj}_S(\underline{\Lambda}_t^{\xi,a^*})=\underline \mu_t^{\xi,a^*}$ $\underline {\mathbb{P}}^{\xi,a^*}$-a.s.; see \eqref{eq:marginal_const}.
				
				Furthermore, by \eqref{eq:worst_mdp2} in Lemma \ref{lem:worst_lift_dynamics} for $(a^*,\underline{\mathbb{P}}^{\xi,a^*})$, it holds that for every $t\geq 0$
				\begin{align*}
					\operatorname{I}_t=\mathbb{E}^{\underline{\mathbb{P}}{}^{\xi,a^*}} \big[\overline{r}\big(\operatorname{pj}_S(\underline{\Lambda}_t^{\xi,a^*}),\underline{\Lambda}_t^{\xi,a^*} \big)+\beta   \overline V^*(\underline \mu_{t+1}^{\xi,a^*}) \big].
				\end{align*}
				
				Combined with \eqref{eq:sum_linear_origin}, this ensures that
				\[
				\overline V^*(\mu)=\sum_{t=0}^\infty \beta^t \mathbb{E}^{\underline{\mathbb{P}}{}^{\xi,a^*}} \big[\overline{r}\big(\operatorname{pj}_S(\underline{\Lambda}_t^{\xi,a^*}),\underline{\Lambda}_t^{\xi,a^*} \big) \big]= \mathbb{E}^{\underline{\mathbb{P}}{}^{\xi,a^*}} \bigg[\sum_{t=0}^\infty \beta^t\overline{r}\big(\operatorname{pj}_S(\underline{\Lambda}_t^{\xi,a^*}),\underline{\Lambda}_t^{\xi,a^*} \big) \bigg].
				\]
				
				Therefore, by the equality $\overline{V}^*(\mu)={V}(\xi)$ (from Step 2 and Step 3), we conclude that 
				\begin{align*}
					\begin{aligned}
						\overline{V}^*(\mu)=V(\xi)
						=\sup_{a\in {\cal A}}{\cal J}^a(\xi) 
						&=\mathbb{E}^{\underline{\mathbb{P}}{}^{\xi,a^*}} \bigg[\sum_{t=0}^\infty \beta^t\overline{r}\big(\operatorname{pj}_S(\underline{\Lambda}_t^{\xi,a^*}),\underline{\Lambda}_t^{\xi,a^*} \big) \bigg]\\
						&=\mathbb{E}^{\underline{\mathbb{P}}{}^{\xi,a^*}} \bigg[\sum_{t=0}^\infty \beta^t r\big(s_t^{\xi,a^*,\underline{\mathbb{P}}{}^{\xi,a^*}},a^*_t,\underline{\Lambda}_t^{\xi,a^*} \big) \bigg] ={\cal J}^{a^*}(\xi) ,
					\end{aligned}
				\end{align*}
				where the last line follows from the same arguments presented for \eqref{eq:lift_reward1}, and the inequality \eqref{eq:step2_ineq_ext} given in Step 2. This completes the proof. \qed

			\section{Proof of results in Section \ref{subsec:closed_Markov}}\label{proof:subsec:closed_Markov}
			\subsection{Proof of Lemma \ref{lem:MDP_closed_loop}}
				We first prove \eqref{eq:kernel_prod}. For simplicity, denote for every $t\geq 0$ by
				\begin{align*}
					\mu_t:=\mu_t^{\xi,\pi^c,\mathbb{P}},\qquad \Lambda_t:=\Lambda_t^{\xi,\pi^c,\mathbb{P}},\qquad \nu_{t+1}:=\mathscr{L}_{\mathbb{P}}(\varepsilon_{t+1}^0|{\cal F}_t^0).
				\end{align*}
				
				As the case for $t= 0$ can be subsumed into the others for $t\geq 1$, we consider the case $t\geq 1$. 
				Since $\Lambda_t$ and $\mu_t$ are ${\cal F}_t^0$ measurable, it is sufficient to show that  for any bounded Borel measurable functions $g:(E^0)^{t}\to \mathbb{R}$ and $f:S\times A\to \mathbb{R}$,
				\begin{align*}
					\mathbb{E}^{\mathbb{P}}\big[g(\varepsilon_{1:t}^0)f\big(s_t^{\xi,\pi^c,\mathbb{P}},a_t^{\pi^c,\mathbb{P}})\big]= \mathbb{E}^{\mathbb{P}}\bigg[g(\varepsilon_{1:t}^0)\int_{S\times A} f(s',a') \pi^c_t\big(d\tilde a|\tilde s,\mu_t\big)\mu_t(d\tilde s)\bigg].
				\end{align*}
				
				Note that $g(\varepsilon_{1:t}^0)$ is ${\cal F}_t^0$ measurable and $s_t^{\xi,\pi^c,\mathbb{P}}$ is ${\cal F}_t$ measurable. Hence, by the distributional constraint that $\mathscr{L}_{\mathbb{P}}(a_t^{\pi^c,\mathbb{P}}|{\cal F}_t)=\pi_t^c(\cdot|s^{\xi,\pi^c,\mathbb{P}}_t,\mu_t)$ $\mathbb{P}$-a.s. (see~\eqref{eq:MKV_represent_closed}) and the tower property, 
				\begin{align*}
					\mathbb{E}^{\mathbb{P}}\big[g(\varepsilon_{1:t}^0)f(s_t^{\xi,\pi^c,\mathbb{P}},a_t^{\pi^c,\mathbb{P}})\big]&=\mathbb{E}^{\mathbb{P}}\Big[g(\varepsilon_{1:t}^0)\mathbb{E}^{\mathbb{P}}\big[\mathbb{E}^{\mathbb{P}}[f(s_t^{\xi,\pi^c,\mathbb{P}},a_t^{\pi^c,\mathbb{P}})\,|\,{\cal F}_t^c]\big|\,{\cal F}_t^0\big]\Big]\\
					&=\mathbb{E}^{\mathbb{P}}\bigg[g(\varepsilon_{1:t}^0)\int_{A}f(s_t^{\xi,\pi^c,\mathbb{P}},\tilde a)\pi^c_t(da'|s_t^{\xi,\pi^c,\mathbb{P}},\mu_t)\bigg]=:\operatorname{I}_t.
				\end{align*}
				Moreover by the definition of $\mu_t$ and its ${\cal F}_t^0$-measurability,
				\begin{align*}
					\operatorname{I}_t&=\mathbb{E}^{\mathbb{P}}\bigg[g(\varepsilon_{1:t}^0) \mathbb{E}^{\mathbb{P}}\bigg[\int_{A}f(s_t^{\xi,\pi^c,\mathbb{P}},\tilde a)\pi^c_t(da'|s_t^{\xi,\pi^c,\mathbb{P}},\mu_t)\,\Big|\,{\cal F}_t^0\bigg]\bigg]\\
					& = \mathbb{E}^{\mathbb{P}}\bigg[g(\varepsilon_{1:t}^0)\int_{S\times A} f(s',a') \pi^c_t\big(d\tilde a|\tilde s,\mu_t\big)\mu_t(d\tilde s)\bigg],
				\end{align*}
				as~claimed.
				
				Moreover, since $\operatorname{pj}_S(\mu_t^{\xi,\pi^c,\mathbb{P}}\mathbin{\hat \otimes} \pi_t^c(\cdot\,|\,\cdot,\mu_t^{\xi,\pi^c,\mathbb{P}}))=\mu_{t}^{\xi,\pi^c,\mathbb{P}}$,  we can the same arguments as in the proof of Proposition \ref{pro:lift_dynamics} (we refer to Section \ref{proof:subsec:lift_MDP}) to get that \eqref{eq:lift_2_closed} holds $\mathbb{P}$-a.s..\qed

				\subsection{Proof of Lemma \ref{lem:worst_lift_dynamics_closed}}
					We first prove \eqref{eq:worst_mdp_closed}. 
					\noindent{\it Step 1:}  Let $\pi^c\in \Pi^c$ be given, and let $\widetilde{\mathbb{P}}\in {\cal Q}$ be some arbitrary. Then set 
					\begin{align}\label{eq:const_init_closed}
						\begin{aligned}
							&\tilde s_0:=\xi,\quad 
							\quad 
							&&\tilde \mu_0:=\mathscr{L}_{\widetilde{\mathbb{P}}}(\tilde s_0),\\
							&\tilde a_0:=\rho_A\big(\pi_0^c(\cdot\,|\,\tilde s_0,\tilde \mu_0),h_0(\vartheta_0)\big),\quad 
						\end{aligned}
					\end{align}
					where $\rho_A$ is the Blackwell-Dubins function on $A$ (see Lemma \ref{lem:BlackDubins}) and $h_0$ is given in Remark \ref{rem:randomize}. Here we note that 
					$\tilde s_0$ is ${\cal F}_0$ measurable (as $\xi \in L_{{\cal F}_0}^0(S)$) and $\tilde a_0$ is ${\cal G}_0$ measurable.

					Then we define by 
					\begin{align}\label{eq:msr_first_closed}
							\underline p_1^{\xi,\pi^c}:=\overline p^*\big(\tilde \mu_0\mathbin{\hat \otimes} \pi_0^c(\cdot\,|\,\cdot,\tilde \mu_0)\big)\in \mathfrak{P}^0,
						\end{align}
						where $\overline{p}^*$ is given in Proposition \ref{pro:dpp}\;(i). 
						
						Next, for every $t\geq 1$ we inductively set
						\begin{align}\label{eq:kernel_general0_closed}
							\begin{aligned}
								&\tilde{s}_t:= \operatorname{F}\big(\tilde s_{t-1},\tilde a_{t-1},\tilde \mu_{t-1}\mathbin{\hat \otimes} \pi_{t-1}^c(\cdot\,|\,\cdot,\tilde \mu_{t-1}) ,\varepsilon_{t},\varepsilon_{t}^0\big),\quad &&\tilde\mu_t:=\mathscr{L}_{\widetilde{\mathbb{P}}}(\tilde s_{t}\,|\,\varepsilon_{1:t}^0),\\
								&\tilde{a}_t:=\rho_A\big(\pi_t^c(\cdot\,|\,\tilde s_t,\tilde \mu_t),h_t(\vartheta_t)\big),\quad  
							\end{aligned}
						\end{align}
						Here, by using the same arguments presented for the proof of Lemma \ref{lem:measurability}\;(ii), we can deduce that $\tilde s_t$ is ${\cal F}_t$ measurable and $\tilde a_t$ is ${\cal G}_t$ measurable. Moreover, $(\tilde \mu_t,\tilde \Lambda_t)$ are ${\cal F}_t^0$ measurable. 
						
						From this, we can consider a Borel measurable function $l_t:(E^0)^t\to {\cal P}(S\times A)$ such that 
						\begin{align}\label{eq:kernel_general_closed}
							l_t (\varepsilon_{1:t}^0)=\tilde \mu_t\mathbin{\hat \otimes} \pi_t^c(\cdot\,|\,\cdot,\tilde \mu_t).
						\end{align}
						Then, define $\underline p_{t+1}^{\xi,\pi^c}:(E^0)^t\ni e_{1:t}^0\mapsto  \underline p_{t+1}^{\xi,\pi^c}(\cdot\,|\,e_{1:t}^0)\in{\cal P}(E^0)$ by
						\begin{align}\label{eq:kernel_general1_closed}
							\underline{p}_{t+1}^{\xi,\pi^c}(\,\cdot\,|\,e_{1:t}^0):= \overline p^*\big(l_t(e_{1:t}^0)\big)\in \mathfrak{P}^0.
						\end{align}
						
						Therefore we can define by $\underline{\mathbb{P}}^{\xi,\pi^c}\in {\cal Q}$ the measure induced by $(\underline{p}_t^{\xi,\pi^c})_{t\geq 1}\in  \mathcal{K}^0$ given in \eqref{eq:msr_first_closed} and \eqref{eq:kernel_general1_closed}.
						
						\vspace{0.5em}
						\noindent {\it Step 2:}  Recall $(\tilde \mu_t)_{t\geq 0}$ given in \eqref{eq:const_init_closed} and \eqref{eq:kernel_general0_closed}. We claim that $\underline{\mathbb{P}}^{\xi,\pi^c}$-a.s.
						\begin{align}\label{eq:induct_msr_const_closed}
							\underline \mu_t^{\xi,\pi^c}=\tilde \mu_t,\quad \mbox{for all $t\geq 0$,}
						\end{align}
						where $\underline{\mu}_{0}^{\xi,a}$ is the law of $s_0^{\xi,a,\underline{\mathbb{P}}^{\xi,a}}$ under $\underline{\mathbb{P}}^{\xi,a}$, and for $t\geq 1$ $\underline{\mu}_{t}^{\xi,a}$ is the conditional law of $s_t^{\xi,a,\underline{\mathbb{P}}^{\xi,a}}$ under $\underline{\mathbb{P}}^{\xi,a}$ given $\varepsilon^0_{1:t}$.
						
						The proof uses an induction over $t\geq 0$: For $t=0$, clearly $s_0^{\xi,\pi^c,\underline{\mathbb{P}}^{\xi,a}}=\tilde s_0=\xi\in L_{{\cal F}_0}^0(S)$. Moreover, since  $\mathscr{L}_{\underline{\mathbb{P}}^{\xi,\pi^c}}(\gamma)=\mathscr{L}_{\tilde{\mathbb{P}}}(\gamma)$ (see Remark \ref{rem:well_dfn_1}\;(ii)), it holds that $\underline \mu_0^{\xi,\pi^c}=\tilde \mu_0$.

						Assume that the induction claim holds for some $t\geq 0$. 
						By ${\cal F}_{t+1}^0$-measurability of  $(\mu_{t+1}^{\xi,\pi^c},\tilde \mu_{t+1})$, it suffices to show that for any bounded Borel measurable functions $\hat g_{t+1}:(E^0)^{t+1}\to \mathbb{R}$ and $\hat f:S\to \mathbb{R}$,
						\begin{align}\label{eq:claim3_worst_closed}
							\mathbb{E}^{\underline{\mathbb{P}}^{\xi,\pi^c}}\big[\hat g_{t+1}(\varepsilon_{1:t+1}^0) \hat f(s_{t+1}^{\xi,\pi^c,\underline{\mathbb{P}}^{\xi,\pi^c}}) \big] = \mathbb{E}^{\underline{\mathbb{P}}^{\xi,\pi^c}}\bigg[\hat g_{t+1}(\varepsilon^0_{1:t+1})\int_{S\times A} f(\tilde s) \tilde \mu_{t+1}(d\tilde s) \bigg].
						\end{align}
						
						Indeed, by the conditional McKean-Vlasov dynamics given in \eqref{eq:MKV_represent_closed} and Fubini's theorem
						\begin{align}
							\mathbb{E}^{\underline{\mathbb{P}}^{\xi,\pi^c}}\big[&\hat g_{t+1}(\varepsilon_{1:t+1}^0) \hat f(s_{t+1}^{\xi,\pi^c,\underline{\mathbb{P}}^{\xi,\pi^c}}) \big]
							=\mathbb{E}^{\underline{\mathbb{P}}^{\xi,\pi^c}}\Big[\hat g_{t+1}(\varepsilon_{1:t+1}^0) \hat f\big(\operatorname{F}(s^{\xi,\pi^c,\mathbb{P}}_t,a_t^{\pi^c,\mathbb{P}},\underline{ \Lambda}_{t}^{\xi,\pi^c},\varepsilon_{t+1},\varepsilon_{t+1}^0)\big) \Big] \nonumber\\
							&\hspace{1.em} =\int_E\mathbb{E}^{\underline{\mathbb{P}}^{\xi,\pi^c}}\Big[\hat g_{t+1}(\varepsilon_{1:t+1}^0) \hat f\big(\operatorname{F}(s^{\xi,\pi^c,\mathbb{P}}_t,a_t^{\pi^c,\mathbb{P}},\underline{ \Lambda}_{t}^{\xi,\pi^c},e,\varepsilon_{t+1}^0)\big) \Big]\lambda_{\varepsilon}(de)=:\operatorname{I}_t,\label{eq:BD_law0}
						\end{align}
						where the second equality holds since $\varepsilon_{t+1}$ is independent of ${\cal G}_t\vee\sigma(\varepsilon_{t+1}^0)$ with $\mathscr{L}_{\underline{\mathbb{P}}^{\xi,\pi^c}}(\varepsilon_{t+1})=\lambda_\varepsilon$ (see Remark \ref{rem:well_dfn_1}\;(i), (ii)).
						
						Moreover, since $\varepsilon_{t+1}^0$ is conditionally independent of ${\cal G}_t$ given ${\cal F}_t^0$ (see Remark \ref{rem:well_dfn_1}\;(iii)) with $\mathscr{L}_{\underline{\mathbb{P}}^{\xi,\pi^c}}(\varepsilon_{t+1}^0|{\cal F}_t^0)=\underline{p}_{t+1}^{\xi,\pi^c}(de^0\,|\,\varepsilon_{1:t}^0)$ (by definition of $\underline{\mathbb{P}}^{\xi,\pi^c}$), and $s^{\xi,\pi^c,\mathbb{P}}_t,a_t^{\pi^c,\mathbb{P}}$, and $\underline{ \Lambda}_{t}^{\xi,\pi^c}$ are all ${\cal G}_t$ measurable, we have
						\begin{align}\label{eq:BD_law1}
							\operatorname{I}_t=\int_E\mathbb{E}^{\underline{\mathbb{P}}^{\xi,\pi^c}}\bigg[\int_{E^0} \Big(\hat g_{t+1}(\varepsilon_{1:t}^0,e^0) \operatorname{D}_{{\cal F}_t^0}(e,e^0) \Big) \underline{p}_{t+1}^{\xi,\pi^c}(de^0\,|\,\varepsilon_{1:t}^0) \bigg]\lambda_{\varepsilon}(de)
						\end{align}
						where for every $(e,e^0)\in E\times E^0$
						\begin{align*}
							\operatorname{D}_{{\cal F}_t^0}(e,e^0):=&\mathbb{E}^{\underline{\mathbb{P}}^{\xi,\pi^c}}\Big[\hat f\big( \operatorname{F}(s^{\xi,\pi^c,\mathbb{P}}_t,a_t^{\pi^c,\mathbb{P}},\underline{ \Lambda}_{t}^{\xi,\pi^c},e,e^0)\big)\,\Big|\,{\cal F}_t^0\Big]\,\\
							=&\int_{S\times A} \hat f\big(\operatorname{F}(s,a,\underline{ \Lambda}_{t}^{\xi,\pi^c},e,e^0)\big)\underline{ \Lambda}_{t}^{\xi,\pi^c}(d s, d a).
						\end{align*}
						
						Moreover, from \eqref{eq:kernel_prod} in Lemma \ref{lem:MDP_closed_loop} it holds for every $(e,e^0)\in E\times E^0$ that $\underline{\mathbb{P}}^{\xi,\pi^c}$-a.s.,
						\begin{align*}
							\begin{aligned}
								\operatorname{D}_{{\cal F}_t^0}(e,e^0)&=\int_{S\times A} \hat f\Big(\operatorname{F}\big(s,a,\big(\underline \mu_t^{\xi,\pi^c}\mathbin{\hat \otimes} \pi_t^c(\cdot\,|\,\cdot,\underline \mu_t^{\xi,\pi^c})\big),e,e^0\big)\Big)\big(\underline \mu_t^{\xi,\pi^c}\mathbin{\hat \otimes} \pi_t^c(\cdot\,|\,\cdot,\underline \mu_t^{\xi,\pi^c})\big)(ds, da)\\
								&=\int_{S\times A} \hat f\Big(\operatorname{F}\big(s,a,\big(\tilde \mu_t\mathbin{\hat \otimes} \pi_t^c(\cdot\,|\,\cdot,\tilde \mu_t)\big),e,e^0\big)\Big)\big(\tilde \mu_t\mathbin{\hat \otimes} \pi_t^c(\cdot\,|\,\cdot,\tilde \mu_t)\big)(ds, da)
							\end{aligned}
						\end{align*}
						where the second inequality follows from the induction assumption at $t$.
						
						Furthermore, since $\tilde s_t$ is ${\cal G}_t$ measurable (noting that ${\cal F}_t\subset {\cal G}_t$), an application of Lemma \ref{lem:law_universal} ensures that $\tilde \mu_t = \mathscr{L}_{\underline{\mathbb{P}}^{\xi,\pi^c}}(\tilde s_t |{\cal F}_t^0)$ $\underline{\mathbb{P}}^{\xi,\pi^c}$-a.s.. This implies that $\underline{\mathbb{P}}^{\xi,\pi^c}$-a.s.
						\begin{align}
								\operatorname{D}_{{\cal F}_t^0}(e,e^0)&=\int_{S\times A} \hat f\Big(\operatorname{F}\big(s,a,\big(\tilde \mu_t\mathbin{\hat \otimes} \pi_t^c(\cdot\,|\,\cdot,\tilde \mu_t)\big),e,e^0\big)\Big)\big(\mathscr{L}_{\underline{\mathbb{P}}^{\xi,\pi^c}}(\tilde s_t |{\cal F}_t^0)\mathbin{\hat \otimes} \pi_t^c(\cdot\,|\,\cdot,\tilde \mu_t)\big)(ds, da)\nonumber\\
								&=\mathbb{E}^{\underline{\mathbb{P}}^{\xi,\pi^c}}\bigg[\int_A \hat f\Big(\operatorname{F}\big(\tilde s_t,a,\big(\tilde \mu_t\mathbin{\hat \otimes} \pi_t^c(\cdot\,|\,\cdot,\tilde \mu_t)\big),e,e^0\big)\Big)\pi_t^c(da\,|\,\tilde s_t,\tilde \mu_t)\,\Big|{\cal F}_t^0\bigg]\label{eq:BD_law2}\\
								&=\mathbb{E}^{\underline{\mathbb{P}}^{\xi,\pi^c}}\Big[ \hat f\Big(\operatorname{F}\big(\tilde s_t,\tilde a_t,\big(\tilde \mu_t\mathbin{\hat \otimes} \pi_t^c(\cdot\,|\,\cdot,\tilde \mu_t)\big),e,e^0\big)\Big)\,\Big|{\cal F}_t^0\Big],\nonumber
						\end{align}
						where the last equality holds by definition of $\tilde a_t$ given in \eqref{eq:kernel_general0_closed} (which follows from the property of the Blackwell-Dubins function and the fact that $\mathscr{L}_{\underline{\mathbb{P}}^{\xi,\pi^c}}(h_t(\vartheta_t))={\cal U}_{[0,1]}$; see Remark \ref{rem:randomize}).
						
						Combining \eqref{eq:BD_law1} with \eqref{eq:BD_law1} and \eqref{eq:BD_law0}, we hence have 
						\begin{align*}
							\mathbb{E}^{\underline{\mathbb{P}}^{\xi,\pi^c}}\big[\hat g_{t+1}(\varepsilon_{1:t+1}^0) \hat f(s_{t+1}^{\xi,\pi^c,\underline{\mathbb{P}}^{\xi,\pi^c}}) \big]&=\mathbb{E}^{\underline{\mathbb{P}}^{\xi,\pi^c}}\Big[\hat g_{t+1}(\varepsilon_{1:t+1}^0)  \hat f\Big(\operatorname{F}\big(\tilde s_t,\tilde a_t,\big(\tilde \mu_t\mathbin{\hat \otimes} \pi_t^c(\cdot\,|\,\cdot,\tilde \mu_t)\big),\varepsilon_{t+1},\varepsilon_{t+1}^0\big)\Big) \Big]\\
							&=\mathbb{E}^{\underline{\mathbb{P}}^{\xi,\pi^c}}\big[\hat g_{t+1}(\varepsilon_{1:t+1}^0) \hat{f}(\tilde s_{t+1}) \big]\\
							&=\mathbb{E}^{\underline{\mathbb{P}}^{\xi,\pi^c}}\bigg[\hat g_{t+1}(\varepsilon_{1:t+1}^0) \int_S\hat{f}(s)\mathscr{L}_{\underline{\mathbb{P}}^{\xi,\pi^c}}(\tilde s_{t+1}|\varepsilon^0_{1:t+1})(ds) \bigg],
						\end{align*}
						where the last line holds by definition of $\tilde s_{t+1}$ given in \eqref{eq:kernel_general0_closed}. 
						
						Moreover, since $\tilde s_{t+1}$ is ${\cal G}_{t+1}$ measurable, another application of Lemma \ref{lem:law_universal} ensures that 
						\[
						\mathscr{L}_{\underline{\mathbb{P}}^{\xi,\pi^c}}(\tilde s_{t+1}|\varepsilon^0_{1:t+1})=\tilde \mu_{t+1},\quad \mbox{$\underline{\mathbb{P}}^{\xi,\pi^c}$-a.s.},
						\]
						which ensures \eqref{eq:claim3_worst_closed} to hold, as claimed. 
						
						By the induction hypothesis, \eqref{eq:induct_msr_const_closed} holds for all $t\geq 0$.
						
						\vspace{0.5em}
						\noindent {\it Step 3:} Recall that $\underline{\mathbb{P}}^{\xi,\pi^c}\in {\cal Q}$ is the measure induced by $(\underline{p}_t^{\xi,\pi^c})_{t\geq 1}\in  \mathcal{K}^0$ given in \eqref{eq:msr_first_closed} and \eqref{eq:kernel_general1_closed} (see Step 1).  Then from Remark \ref{rem:well_dfn_1}\;(iii), it holds that  $\underline{\mathbb{P}}^{\xi,a}$-a.s.
						\begin{align}\label{eq:equiv_1_closed}
							\begin{aligned}
								&\mathscr{L}_{\underline{\mathbb{P}}^{\xi,\pi^c}}(\varepsilon_{1}^0)=\underline p_1^{\xi,\pi^c}\in \mathfrak{P}^0,\\\;\; &\mathscr{L}_{\underline{\mathbb{P}}^{\xi,\pi^c}}(\varepsilon_{t}^0|{\cal F}_{t-1}^0)= \underline{p}_t^{\xi,\pi^c}(\cdot|\varepsilon_{1:t-1}^0)\in \mathfrak{P}^0\;\; \mbox{for all $t\geq 2$}.
							\end{aligned}
						\end{align}
						Moreover, by \eqref{eq:induct_msr_const_closed} in Step 2 and \eqref{eq:kernel_prod}  in Lemma \ref{lem:MDP_closed_loop}, it holds that $\underline{\mathbb{P}}^{\xi,\pi^c}$-a.s.
						\begin{align}\label{eq:equiv_2_closed}
							\underline p_1^{\xi,\pi^c}=\overline p^*(\underline \Lambda_0^{\xi,\pi^c}),\qquad \underline{p}_t^{\xi,\pi^c}(\cdot|\varepsilon_{1:t-1}^0)= \overline p^*\big(\underline \Lambda_{t-1}^{\xi,\pi^c}\big) \quad \mbox{for all $t\geq 2$},
						\end{align}
						which ensures \eqref{eq:worst_mdp_closed} to hold, as claimed. 
						
						A direct consequence of \eqref{eq:lift_2_closed} ensures \eqref{eq:worst_mdp2_closed} to hold, as claimed. This completes the proof. \qed

			\subsection{Proof of Corollary \ref{cor:origin_MDP_closed}}
				As the essential arguments of the proof closely follow those of Theorem~\ref{thm:origin_MDP}, we provide the outline of the proof and omit some details here.
				
				\noindent {\it Step 1.} For notational simplicity, set $\mu:=\mathscr{L}(\xi)$. We first consider for every $n\in \mathbb{N}$ 
				\begin{align*}
					{\cal I}_n^{\xi,\pi^{c,*}}:=\inf_{{\mathbb{P}}\in {\cal Q}} \mathbb{E}^{{\mathbb{P}}} \bigg[\sum_{t=0}^{n-1}\beta^t\,{r}(s_{t}^{\xi,\pi^{c,*},{\mathbb{P}}},a^{\pi^{c,*},\mathbb{P}}_t,{\Lambda}_{t}^{\xi,\pi^{c,*},\mathbb{P}})+\beta^n\,\overline{V}^*({\mu}_{n}^{\xi,\pi^{c,*},\mathbb{P}}) \bigg],
				\end{align*}
				where for each $\mathbb{P}\in {\cal Q}$, $({ \mu}_{t}^{\xi,\pi^{c,*},\mathbb{P}})_{t\geq 0}$ and $({\Lambda}_{t}^{\xi,\pi^{c,*},\mathbb{P}})_{t\geq 0}$ are given in \eqref{eq:MDP_state_closed}. 
				
				Note that by \eqref{eq:kernel_prod} in Lemma \ref{lem:MDP_closed_loop} and definition of $\pi^{c,*}_t=\pi_{\operatorname{loc}}^{c,*}$ given in \eqref{eq:local_max_ext} together with the property \eqref{eq:local_max_ext_integ}, it holds for every $\mathbb{P}\in {\cal Q}$ that  $\mathbb{P}$-a.s.,
				\[
				\overline \pi^*({\mu}_t^{\xi,\pi^{c,*},\mathbb{P}})={\Lambda}_{t}^{\xi,\pi^{c,*},\mathbb{P}}\quad \mbox{for all $t\geq 0$.}
				\]
				From this, using the same arguments presented for \eqref{eq:lift_reward1}, we have that  for every $n\in \mathbb{N}$
				\begin{align*}
					{\cal I}_n^{\xi,\pi^{c,*}}=\inf_{{\mathbb{P}}\in {\cal Q}} \mathbb{E}^{{\mathbb{P}}} \bigg[\sum_{t=0}^{n-1}\beta^t\,\overline{r}({\mu}_t^{\xi,\pi^{c,*},\mathbb{P}}\,,\,{\Lambda}_{t}^{\xi,\pi^{c,*},\mathbb{P}})+\beta^n\,\overline{V}^*({\mu}_{n}^{\xi,\pi^{c,*},\mathbb{P}}) \bigg].
				\end{align*}
				Hence, from the representation of the Markov decision process of the lifted state process in  \eqref{eq:lift_2_closed} (see Lemma \ref{lem:MDP_closed_loop}), we can use the same arguments presented for Steps 1 and 2 in the proof of Theorem \ref{thm:origin_MDP} (that relies on the local optimality of $\overline \pi^*({\mu}_t^{\xi,\pi^{c,*},\mathbb{P}})$ to ${\cal T}\overline V^*({\mu}_t^{\xi,\pi^{c,*},\mathbb{P}})$ in Proposition~\ref{pro:dpp}\;(ii) and the fixed point theorem in Proposition \ref{pro:fixed_point}; see Section \ref{proof:subsec:main_thm}) to have  
				\begin{align*}
					\overline{V}^*(\mu)\leq \limsup_{n\to \infty}{\cal I}_n^{\xi,\pi^{c,*}}\leq{\cal J}^{\pi^{c,*}}(\xi)  \leq {V}^c(\xi).
				\end{align*}
				
				\noindent {\it Step 2.} For every $\pi^c\in \Pi^c$, let $\underline {\mathbb{P}}^{\xi,\pi^c}\in {\cal Q}$ be induced by some $(\underline{p}_t^{\xi,\pi^c})_{t\geq 1}\in \mathcal{K}^0$ satisfying \eqref{eq:worst_mdp} and~\eqref{eq:worst_mdp_closed} (see Lemma \ref{lem:worst_lift_dynamics_closed}).  Then define ${\cal V}^{\pi^c}(\xi)$ by
				\begin{align*}
					{{\cal V}}^{\pi^c}(\xi):=\mathbb{E}^{\underline{\mathbb{P}}{}^{\xi,\pi^c}}\bigg[\sum_{t=0}^\infty \beta^t r(s_t^{\xi,\pi^c,\underline{\mathbb{P}}{}^{\xi,\pi^c}},a^{\pi^c,\underline{\mathbb{P}}{}^{\xi,\pi^c}}_t,\underline{\Lambda}_t^{\xi,\pi^c}) \bigg]
					=\mathbb{E}^{\underline{\mathbb{P}}{}^{\xi,\pi^c}} \bigg[\sum_{t=0}^\infty\beta^t\;\overline{r}\big(\operatorname{pj}_S(\underline{\Lambda}_t^{\xi,\pi^c}),\underline{\Lambda}_t^{\xi,\pi^c} \big)\bigg],
				\end{align*}
				where $\underline{ \Lambda}_{t}^{\xi,\pi^c}$ is the conditional joint law of $(s_{t}^{\xi,\pi^c,\underline{\mathbb{P}}^{\xi,\pi^c}},a^{\pi^c,\mathbb{P}}_t)$ under $\underline {\mathbb{P}}^{\xi,\pi^c}$ given $\varepsilon^0_{1:t}$.
				
				By the local optimality of $\overline{p}^*(\underline{ \Lambda}_{t}^{\xi,\pi^c})$ to ${\cal T}\overline V^*(\operatorname{pj}_S(\underline{ \Lambda}_{t}^{\xi,\pi^c}))$ (see Proposition~\ref{pro:dpp}\;(i)), we can use the same arguments presented for Step 3 in the proof of Theorem \ref{thm:origin_MDP} to have
				\[
				V^c(\xi)\leq \sup_{\pi^c\in \Pi^c}{\cal V}^{\pi^c}(\xi)\leq \sup_{\pi^c\in \Pi^c}\sum_{t=0}^\infty\Big(\beta^t \, \mathbb{E}^{\underline{\mathbb{P}}{}^{\xi,\pi^c}} [ \overline {V}^*(\underline{\mu}_t^{\xi,\pi^c}) ] -\beta^{t+1} \, \mathbb{E}^{\underline{\mathbb{P}}{}^{\xi,\pi^c}}  [\overline V^*(\underline{\mu}_{t+1}^{\xi,\pi^c})] \Big)= \overline V^*(\mu),
				\] 
				where $\underline{ \mu}_{t}^{\xi,\pi^c}$ is the conditional law of $s_{t}^{\xi,\pi^c,\underline{\mathbb{P}}^{\xi,\pi^c}}$ under $\underline {\mathbb{P}}^{\xi,\pi^c}$ given $\varepsilon^0_{1:t}$.
				
				Therefore, we have obtained that $\overline{V}^*(\mu)=V^{c}(\xi)$, as claimed. In fact, $\overline{V}^*(\mu)=V(\xi)$ follows from Theorem \ref{thm:origin_MDP}. Hence the statement (i) holds.
				
				\vspace{0.5em}
				\noindent 
				{\it Step 3.} Lastly, we consider $\underline{\mathbb{P}}^{\xi,\pi^{c,*}}\in {\cal Q}$ which is induced by $(\underline{p}_t^{\xi,\pi^{c,*}})_{t\geq 1}\in \mathcal{K}^0$ satisfying \eqref{eq:worst_mdp_closed} and~\eqref{eq:worst_mdp2_closed} (see Lemma \ref{lem:worst_lift_dynamics_closed}). Then 
				by definition of $\pi^{c,*}$ and of $\underline{\mathbb{P}}^{\xi,\pi^{c,*}}$ (noting that both satisfy the local optimality given in Proposition \ref{pro:dpp}),  it holds that for every $t\geq 0$ 
				\begin{align*}
					&\mathbb{E}^{\underline{\mathbb{P}}{}^{\xi,\pi^{c,*}}} [ \overline {V}^*(\underline{\mu}_t^{\xi,\pi^{c,*}}) ]=\mathbb{E}^{\underline{\mathbb{P}}{}^{\xi,\pi^{c,*}}} [ {\cal T}\overline {V}^*(\underline{\mu}_t^{\xi,\pi^{c,*}}) ]\nonumber\\
					&\hspace{2.em} =\mathbb{E}^{\underline{\mathbb{P}}{}^{\xi,\pi^{c,*}}} \bigg[\overline{r}\big(\operatorname{pj}_S(\underline{\Lambda}_t^{\xi,\pi^{c,*}}),\underline{\Lambda}_t^{\xi,\pi^{c,*}} \big)+\beta  \int_{{\cal P}(S)} \overline V^*(\tilde \mu')  \overline p\big(d\tilde \mu'|\operatorname{pj}_S(\underline{\Lambda}_t^{\xi,\pi^{c,*}}),\underline{\Lambda}_t^{\xi,\pi^{c,*}} ,\overline p^*(\underline{\Lambda}_t^{\xi,\pi^{c,*}} )\big) \bigg].
				\end{align*}
				Hence by using the same arguments presented for Step 4 of the proof of Theorem \ref{thm:origin_MDP}, we deduce that \eqref{eq:verify_closed} holds. This completes the proof. \qed

	\appendix
	\section{Supplementary statements}
		Let us provide some elementary observations on conditional laws. 
		\begin{lem}\label{lem:regular_cond}
			Fix a probability space $(\tilde \Omega, \tilde {\cal F}, \tilde {\mathbb P})$. Let $X$ be Borel space and $Y$ be measurable space. For every random elements ${\cal X}$ and ${\cal Y}$ with values in $X$ and $Y$, respectively, the following~hold:
			\begin{enumerate}[leftmargin=3.em]
				\item [(i)] There exists a kernel $k^{{\cal X}|{\cal Y}}:Y\ni y\mapsto  k^{{\cal X}|{\cal Y}}(dx|y)\in {\cal P}(X)$ such that for~every $B\in {\cal B}(X)$, $\tilde {\mathbb{P}}({\cal X}\in B|{\cal Y}) = k^{{\cal X}|{\cal Y}}(B|{\cal Y})$ $\tilde{\mathbb{P}}$-a.s., and $k^{{\cal X}|{\cal Y}}$ is unique $\mathscr{L}_{\tilde{\mathbb{P}}}({\cal Y})$-a.e.. As a consequence, $k^{{\cal X}|{\cal Y}}(\cdot\,|\,{\cal Y})$ is $\sigma({\cal Y})$ measurable and we denote for every $\tilde \omega\in \tilde \Omega$
				\[
				\mathscr{L}_{\tilde{\mathbb{P}}}({\cal X}|{\cal Y})(\tilde \omega):=k^{{\cal X}|{\cal Y}}(\cdot|{\cal Y})(\tilde \omega),
				\]
				i.e., a conditional law of ${\cal X}$ given ${\cal Y}$; see, e.g., \cite[Section 6, p.106--107]{kallenberg2002foundations}. 
				\item [(ii)] If ${\cal X}$ is given by ${\cal X}=\varphi({\cal Y}, {\cal Z})$, where $\varphi:Y\times Z\to X$ is a measurable function and ${\cal Z}$ is a random element in $Z$ and independent of~${\cal Y}$, then $\mathscr{L}_{\tilde{\mathbb{P}}}({\cal X}|{\cal Y})=\mathscr{L}_{\tilde{\mathbb{P}}}(\varphi(y,{\cal Z}))|_{y={\cal Y}}$ and $\mathscr{L}_{\tilde{\mathbb{P}}}({\cal X}|{\cal Y})$ is $\sigma({\cal Y})$ measurable. 
			\end{enumerate}
		\end{lem}
		\begin{proof}
			Part\;(i) is shown in \cite[Theorem 6.3]{kallenberg2002foundations}. We proceed to prove (ii), which is a consequence of~(i) with an application of Fubini's theorem. Clearly, it is sufficient to show that for any bounded measurable function $g:Y\to \mathbb{R}$ and bounded Borel measurable function $f:X\to \mathbb{R}$,
			\begin{align*}
				\mathbb{E}^{\tilde {\mathbb{P}}}\bigg[g({\cal Y})\int_X f(x')\mathscr{L}_{\tilde{\mathbb{P}}}({\cal X}|{\cal Y})(dx')\bigg]= \mathbb{E}^{\tilde{\mathbb{P}}}\bigg[g({\cal Y})\int_X f(x')\mathscr{L}_{\tilde{\mathbb{P}}}(\varphi(y,{\cal Z}))|_{y={\cal Y}}(dx')\bigg].
			\end{align*}
			
			Indeed, by definition of the conditional law $\mathscr{L}_{\tilde{\mathbb{P}}}({\cal X}|{\cal Y})$ (given in (i)) it holds that 
			\[
			\mathbb{E}^{\tilde {\mathbb{P}}}\bigg[g({\cal Y})\int_X f(x')\mathscr{L}_{\tilde{\mathbb{P}}}({\cal X}|{\cal Y})(dx')\bigg]= \mathbb{E}^{\tilde {\mathbb{P}}}\big[g({\cal Y})\mathbb{E}^{\tilde {\mathbb{P}}}[f({\cal X})|{\cal Y}]\big]= \mathbb{E}^{\tilde {\mathbb{P}}}[g({\cal Y})f({\cal X})]=:\operatorname{I},
			\]
			where the second equality follows from the $\sigma({\cal Y})$-measurability of $g({\cal Y})$ and the tower property.
			
			Moreover since ${\cal X}=\varphi({\cal Y}, {\cal Z})$, and 
			${\cal Y}$ and ${\cal Z}$ are independent,
			\begin{align*}
				\operatorname{I}= \mathbb{E}^{\tilde {\mathbb{P}}}\Big[g({\cal Y})\mathbb{E}^{\tilde {\mathbb{P}}}\big[f(\varphi({\cal Y},{\cal Z}))|{\cal Y}\big]\Big]&= \int_{Y} g(y) \mathbb{E}^{\tilde {\mathbb{P}}}\Big[f(\varphi(y,{\cal Z}))\Big] \mathscr{L}_{\tilde{\mathbb{P}}}({\cal Y})(dy)\\
				&=\int_{Y} g(y)\mathbb{E}^{\tilde {\mathbb{P}}}\bigg[\int_Xf(x')\mathscr{L}_{\tilde{\mathbb{P}}}(\varphi(y,{\cal Z}))(dx')\bigg] \mathscr{L}_{\tilde{\mathbb{P}}}({\cal Y})(dy)\\
				&=\mathbb{E}^{\tilde{\mathbb{P}}}\bigg[g({\cal Y})\int_Xf(x')\mathscr{L}_{\tilde{\mathbb{P}}}(\varphi(y,{\cal Z}))|_{y={\cal Y}}(dx')\bigg],
			\end{align*}
			where the second equality follows from definition of $\mathscr{L}_{\tilde{\mathbb{P}}}(\varphi(y,{\cal Z}))$ and the last one follows from Fubini's theorem (since both $f$ and $g$ are bounded). The $\sigma({\cal Y})$-measurability of $\mathscr{L}_{\tilde{\mathbb{P}}}({\cal X}|{\cal Y})$ follows from (i). This concludes the proof.
	\end{proof}

	\begin{lem}[Blackwell and Dubins \cite{blackwell1983extension}]\label{lem:BlackDubins}
		For any Polish space $X$, there exists a Borel measurable function $\rho_X:{\cal P}(X)\times [0,1]\to X$ satisfying the following conditions:
		\begin{itemize}
			\item [(i)] for every $\lambda\in {\cal P}(X)$ and every uniform random variable $U\sim {\cal U}_{[0,1]}$, $\rho_X(\lambda,U)$ is distributed according to $\lambda$; 
			\item [(ii)] for almost every $u$, the map $\lambda\mapsto \rho_X(\lambda,u)$ is continuous w.r.t.\;the weak topology of~${\cal P}(X)$.
		\end{itemize}
		We call $\rho_X$ the Blackwell--Dubins function of the space $X$.
	\end{lem}

	\begin{lem}[Universal disintegration; {see, e.g., \cite[Corollarly 1.26]{kallenberg2017random}}]\label{lem:univ_disint} For any Borel spaces $X$ and $Y$, there exists a kernel ${\cal K}_{X\times Y}:X\times {\cal P}(X\times Y)\times {\cal P}(X)\ni(x,\lambda,\eta)\mapsto {\cal K}_{X\times Y}(\cdot|x,\lambda,\eta)\in {\cal P}(Y)$  
		such that for every $\lambda \in {\cal P}(X\times Y)$ and $\eta \in {\cal P}(X)$ satisfying $\operatorname{pj}_X(\lambda)\ll \eta$, it holds that
		\[
		\lambda = \eta \mathbin{\hat \otimes} {\cal K}_{X\times Y}(\,\cdot\,|\,\cdot,\lambda, \eta),
		\]
		Moreover, ${\cal K}_{X\times Y}(\,\cdot\,|,\cdot,\lambda, \eta)$ is unique $\eta$-a.e. for fixed $\lambda$ and $\eta$.
	\end{lem}

\bibliographystyle{abbrv}
\bibliography{references}

\end{document}